\documentclass{article}
\usepackage[top=1in, bottom=1in, left=1in, right=1in]{geometry}
\usepackage[T1]{fontenc}
\usepackage{lmodern}
\usepackage[UKenglish]{babel}
\usepackage{amsmath,amsthm,amssymb,mathtools,graphicx,float}
\usepackage{microtype,mathrsfs,slashed}
\usepackage{caption,longtable,booktabs,needspace,multicol}
\usepackage[colorlinks=true, pdfstartview=Fit,linkcolor=blue,citecolor=blue,urlcolor=blue]{hyperref}

\let\LinkSHodgeSin\sin
\let\LinkSHodgeMathOpen\(
\let\LinkSHodgeMathClose\)

\newcommand{\rG}{{\rm G}}

\newcommand{\bw}{{\bf w}}

\newcommand{\bF}{{\bf F}}

\newcommand{\cB}{\mathcal{B}}

\newcommand{\cE}{\mathcal{E}}
\newcommand{\cF}{\mathcal{F}}

\newcommand{\cJ}{\mathcal{J}}

\newcommand{\cL}{\mathcal{L}}

\newcommand{\cO}{\mathcal{O}}

\newcommand{\cX}{\mathcal{X}}

\newcommand{\fm}{{\mathfrak m}}

\newcommand{\Z}{\mathbb{Z}}
\newcommand{\Q}{\mathbb{Q}}
\newcommand{\R}{\mathbb{R}}
\newcommand{\C}{\mathbb{C}}

\renewcommand{\bF}{\mathbb{F}}

\newcommand{\Spin}{{\rm Spin}}
\newcommand{\SU}{{\rm SU}}

\newcommand{\SL}{\mathrm{SL}}

\DeclareMathOperator{\rank}{rank}

\renewcommand{\det}{\operatorname{det}}

\renewcommand{\epsilon}{\varepsilon}
\DeclareMathOperator{\sen}{sen}
\renewcommand{\sin}{\sen}

\def\({\mathopen{}\left(}
\def\){\right)\mathclose{}}
\def\<{\mathopen{}\left<}
\def\>{\right>\mathclose{}}

\newcommand{\qandq}{\quad\text{e}\quad}

\newcommand{\qforq}{\quad\text{para}\quad}

\usepackage{multicol,color}

\definecolor{gold}{rgb}{0.85,.66,0}
\definecolor{cherry}{rgb}{0.9,.1,.2}
\definecolor{burgundy}{rgb}{0.8,.2,.2}
\definecolor{orangered}{rgb}{0.85,.3,0}
\definecolor{orange}{rgb}{0.85,.4,0}
\definecolor{olive}{rgb}{.45,.4,0}
\definecolor{lime}{rgb}{.6,.9,0}
\definecolor{green}{rgb}{.2,.7,0}
\definecolor{grey}{rgb}{.4,.4,.2}
\definecolor{brown}{rgb}{.4,.3,.1}

\usepackage{xparse}

\makeatletter
\NewDocumentCommand{\NewExampleTag}{mm}{%
  \expandafter\newcommand\csname tag#1\endcsname[1]{%
    \ifdefined\mmriemexercisereplayactive\else
      \addcontentsline{explist#1}{#1}{##1}%
    \fi
    \ifdefined\mmriemexerciseactive
      \textnormal{(cf.~\ref{exemplos:#2}).}\space
    \fi
  }%
  \expandafter\newcommand\csname showtag#1\endcsname{%
    \section{#2}\label{exemplos:#2}%
    \@starttoc{explist#1}%
  }%
  \expandafter\newcommand\csname l@#1\endcsname[2]{%
    \textbf{##1} -- pág.~##2\par
  }%
}
\makeatother

\NewExampleTag{euclid}{Espaço euclidiano}
\NewExampleTag{esfera}{Esfera}
\NewExampleTag{reproj}{Espaço projetivo real}
\NewExampleTag{comproj}{Espaço projetivo complexo}
\NewExampleTag{regval}{Variedades de nível}
\NewExampleTag{superficie}{Superfícies regulares}
\NewExampleTag{homog}{Espaços homogêneos}
\NewExampleTag{lie}{Grupos de Lie}

\renewcommand{\qforq}{\quad\text{for}\quad}
\renewcommand{\qandq}{\quad\text{and}\quad}
\let\sin\LinkSHodgeSin
\let\(\LinkSHodgeMathOpen
\let\)\LinkSHodgeMathClose
\newtheorem{theorem}{Theorem}[section]
\newtheorem{maintheorem}{Theorem}

\renewcommand{\thefootnote}{\fnsymbol{footnote}}

\newcommand{\ghlogo}{\raisebox{-0.2em}{\includegraphics[height=1em]{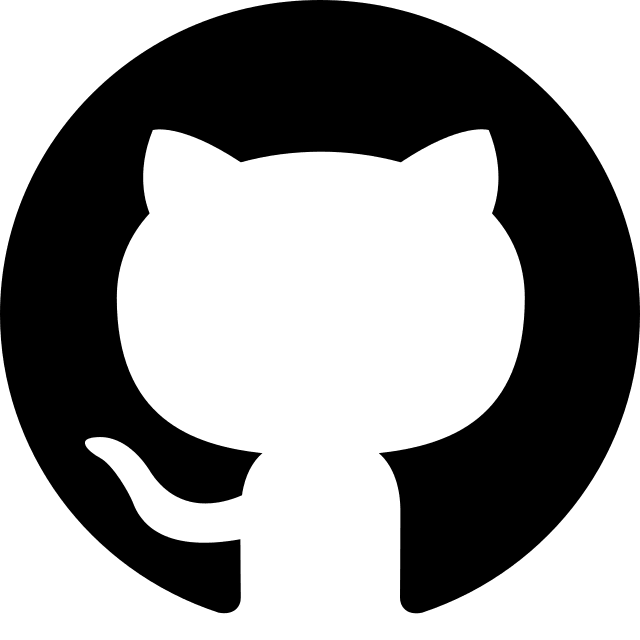}}}

\newtheorem{lemma}[theorem]{Lemma}
\newtheorem{proposition}[theorem]{Proposition}
\newtheorem{corollary}[theorem]{Corollary}
\theoremstyle{definition}
\newtheorem{remark}[theorem]{Remark}
\newtheorem{definition}[theorem]{Definition}
\newtheorem{example}[theorem]{Example}
\theoremstyle{plain}
\newcommand{\WP}{\mathbb{P}}
\newcommand{\MM}{\mathbb{M}}
\DeclareMathOperator{\age}{age}
\providecommand{\coloneq}{\coloneqq}
\newtheorem{introconjecture}{Conjecture}

\title{Sasakian Hodge Numbers of Calabi--Yau Links} 
\author{Edward Hirst\thanks{\texttt{ehirst@unicamp.br}}, Henrique N. S\'a Earp\thanks{\texttt{hqsaearp@unicamp.br}}, Tom\'as S. R. Silva\thanks{\texttt{tomas@ime.unicamp.br}}}
\date{\small\today \\
\textit{\small Instituto de Matem\'atica, Estat\'istica e Computa\c{c}\~ao Cient\'ifica (IMECC), \\ Universidade Estadual de Campinas (UNICAMP), 13083-859, Brazil}\\
}

\begin{document}
\maketitle
\renewcommand{\thefootnote}{\arabic{footnote}}
\setcounter{footnote}{0}

\begin{abstract}
The link of an isolated weighted homogeneous hypersurface singularity whose degree equals the sum of its weights is a Sasakian manifold that admits a contact Calabi--Yau structure, and its Sasakian Hodge numbers are dimensions of graded pieces of the Milnor algebra of the defining polynomial.
In earlier work these numbers were computed with Gr\"obner bases for the $7$-dimensional links of $7549$ of the $7555$ weighted Calabi--Yau threefold links, and two conjectures were raised from the data.
The first asserts that polynomials with the same weights have Milnor algebras whose relevant graded pieces have equal dimensions; the second asserts that the Sasakian Hodge number $h^{2,1}_S$ of the link is bounded above by the Hodge number $h^{2,1}$ of the base Calabi--Yau threefold.
Both conjectures are proved here.
The first holds for every polynomial with an isolated singularity, because its partial derivatives form a regular sequence; an independent derivation of the Sasakian Hodge numbers by sheaf cohomology on the leaf-space orbifold is also given.
For the second, the stringy Hodge numbers of the Calabi--Yau orbifold are shown to bound the middle Sasakian Hodge numbers in every dimension $n\ge2$, and for $n\ge3$ the gap in $h^{n-1,1}$ is shown to be a sum of geometric genera of loci of orbifold points; for threefolds the bound is attained exactly when all curves of orbifold points are rational.
The resulting closed formulae give the Sasakian Hodge numbers of the six $7$-dimensional links left open in the earlier work and complete datasets for all $95$ five-dimensional K3 links and all $1\,100\,055$ nine-dimensional Calabi--Yau fourfold links; for the K3 links they recover the previously known second Betti numbers.

\small Code and data: \ghlogo\ \href{https://github.com/GeometrIALab/Link_SasakianHodgeNumbers}{\texttt{github.com/GeometrIALab/Link\_SasakianHodgeNumbers}}.
\end{abstract}
\numberwithin{equation}{section}
\numberwithin{figure}{section}
\numberwithin{table}{section}
\newpage
\begingroup
\small
\section*{Contents}
\makeatletter
\@starttoc{toc}
\makeatother
\endgroup

\bigskip \bigskip
\paragraph{Acknowledgements.}
\noindent Edward Hirst acknowledges support from the S\~ao Paulo Research Foundation (FAPESP) under
grant \mbox{[2024/18994-7]}.

\noindent Henrique N. S\'a Earp acknowledges support from FAPESP grants \mbox{[2021/04065-6]}, associated with
\emph{BRIDGES: Brazil--France interplays in Gauge Theory, extremal structures
and stability}; \mbox{[2024/00923-6]}, associated with \emph{CBG:
Brazilian Centre for Geometry}; and \mbox{[2020/09838-0]}, associated with
\emph{BI0S: Brazilian Institute of Data Science}. He also acknowledges support from the Brazilian
National Council for Scientific and Technological Development (CNPq),
through Research Productivity Fellowship \mbox{[307145/2025-5]}, level PQ-A.

\noindent Tom\'as S. R. Silva acknowledges support from FAPESP grant \mbox{[2022/09891-4]}.

This work was supported in part by a grant of access to OpenAI models through the ChatGPT for Academic Researchers program.

\paragraph{AI usage statement.}
\noindent The proofs of the central results of this paper, Conjectures~\ref{conj:wd} and~\ref{conj:cyb}, were developed with substantial assistance from the Claude Fable (Anthropic).
Claude was used to explore proof strategies and to produce the initial versions of the arguments, which were then reviewed by the authors.
Some parts of the proofs are still being refined for clarity of exposition.
Claude also assisted with writing the first draft of the manuscript, which was subsequently heavily reviewed and edited by the authors with the further assistance of both Claude and GPT Astra (OpenAI).
The authors take responsibility for the mathematical content of this paper and the correctness of its results.
\newpage

\section{Introduction}
\label{sec:intro}

Weighted hypersurface singularities provide a common setting for Sasakian geometry, Hodge theory and the topology of links.
The problem considered here is to determine which Hodge-theoretic invariants of a Calabi--Yau link are fixed by the weights of its defining equation, and how they compare with those of the corresponding Calabi--Yau hypersurface.
A weighted homogeneous polynomial $f$ in $N$ complex variables with an isolated singularity at the origin determines two closely related geometric objects.
The first is the weighted projective hypersurface 
$$X_f\coloneq \{f=0\}\subset\WP^{N-1}(\bw),$$ 
a projective variety with cyclic quotient singularities, which is the coarse space of a compact complex orbifold $\cX_f$ of dimension $n=N-2$.
The second is the \emph{link} $$K_f\coloneq f^{-1}(0)\cap S^{2N-1},$$ a smooth compact manifold of dimension $2n+1$, which is a Seifert circle orbibundle over $\cX_f$.
The link inherits a quasi-regular Sasakian structure whose Reeb flow is the weighted circle action.
When the degree of $f$ equals the sum of the weights, the orbifold $\cX_f$ has trivial canonical bundle, and after a deformation of its transverse K\"ahler metric the link carries a \emph{contact Calabi--Yau} structure \cite{TomassiniVezzoni2008,HabibVezzoni2015,CalvoAndradeRodriguezSaEarp2020}.
In real dimension $7$ such structures induce $\rG_2$-structures \cite{HabibVezzoni2015,CalvoAndradeRodriguezSaEarp2020}, which is one of the reasons for interest in them \cite{LotaySaEarp2023}.

The \emph{Sasakian Hodge numbers} $h^{p,q}_S$ of the link are the dimensions of the spaces of harmonic forms of type $(p,q)$ in degrees up to the middle dimension.
Apart from $h^{0,0}_S=1$ they vanish below the middle degree, because the link is highly connected (Remark~\ref{rem:diamond}).
The remaining ones give a Hodge decomposition of the middle cohomology of $K_f$, and in particular determine its Betti number, but not finer invariants such as torsion or the diffeomorphism type.
By a theorem going back to Steenbrink and Itoh \cite{Steenbrink1977,Itoh2004}, they are the dimensions of certain graded pieces of the Milnor algebra $$\MM_f\coloneq \C[z]/(\partial_1f,\dots,\partial_Nf).$$
This quotient imposes the partial derivatives as relations among polynomials; its grading turns the Hodge-number calculation into a count in specified weighted degrees (Section~\ref{sec:weights}).
In \cite{AHHHSS2024} this description was turned into an algorithm based on Gr\"obner bases, to compute the Sasakian Hodge numbers and the Crowley--Nordstr\"om invariant \cite{CrowleyNordstrom2015} of the $7$-dimensional links for $7549$ of the $7555$ weighted Calabi--Yau threefolds \cite{CandelasLynkerSchimmrigk1990,KlemmSchimmrigk1994,KreuzerSkarke1992}; the resulting dataset was then studied with machine learning.
The invariant $\nu\in\Z/48$ is attached to the homotopy class of the induced $\rG_2$-structure; its definition and the Milnor-fibre formula used here are recalled in Section~\ref{sec:nu}.
The interpretation of these experiments led to the Weight Determinacy and Calabi--Yau Bound conjectures of \cite{AHHHSS2024}, stated below as Conjectures~\ref{conj:wd} and~\ref{conj:cyb} after fixing notation.

\begin{definition}[Admissible weighted hypersurfaces]
\label{def:intro-setting}
Let $N=n+2\ge4$ and let $\bw\coloneq (w_1,\dots,w_N)$ be positive integer weights with $\gcd(w_1,\dots,w_N)=1$.
Give $S\coloneq \C[z_1,\dots,z_N]$ the grading $\deg z_i=w_i$.
A polynomial $f\in S_d$ is \emph{admissible} if it has an isolated critical point at the origin.
For such an $f$, retain the notation $K_f$, $\cX_f$ and $X_f$ for its link, quotient orbifold and coarse hypersurface, respectively, and write
\[
Y\coloneq f^{-1}(0)\setminus0,\qquad \MM_f\coloneq S/(\partial_1f,\dots,\partial_Nf).
\]
The data are \emph{Calabi--Yau} when $d=\sum_iw_i$.
Their \emph{monomial counts} and \emph{formal coefficients} $a_k$ are defined by
\[
p_{\bw}(k)\coloneq \dim S_k,\qquad
P_{\bw}(t)\coloneq \prod_{i=1}^N\frac{1-t^{d-w_i}}{1-t^{w_i}}=\sum_{k\ge0}a_kt^k,
\]
respectively, with $p_{\bw}(k)\coloneq0$ for $k<0$.
\end{definition}

The normalisation $\gcd(\bw)=1$ removes a common ineffective stabiliser from the weighted action.
It is essential when comparing the quotient orbifold with the stringy invariants of its coarse space in Theorems~\ref{thm:defect}--\ref{prop:lw} below. 
The Hilbert-series identity itself also holds for unnormalised weights, with the degrees regraded accordingly.

\begin{definition}[Weak R-equivalence]
\label{def:weak-r-equivalence}
In the setting of Definition~\ref{def:intro-setting}, let $f,g\in S_d$ be admissible.
Following \cite{AHHHSS2024}, they are \emph{weakly R-equivalent} if
\[
\dim(\MM_f)_{\ell_p}=\dim(\MM_g)_{\ell_p}
\qquad(0\le p\le n),\qquad
\ell_p\coloneq(p+1)d-\sum_iw_i.
\]
For weight vectors agreeing up to permutation, the coordinates are first relabelled to identify the gradings.
These are precisely the degrees computing the Sasakian Hodge numbers in Theorem~\ref{thm:itoh}.
\end{definition}

\begin{introconjecture}[Weight Determinacy (WD); Conjecture 7 of \cite{AHHHSS2024}]
\label{conj:wd}
Two admissible weighted homogeneous polynomials of the same degree whose weight vectors agree up to permutation are weakly R-equivalent in the sense of Definition~\ref{def:weak-r-equivalence}.
\end{introconjecture}

The admissibility hypothesis makes explicit the isolated-singularity assumption in the discussion surrounding the original conjecture; it is absent from its displayed formulation in \cite{AHHHSS2024}.
The same discussion predicts equality of the Sasakian Hodge numbers and, for seven-dimensional Calabi--Yau links, of the Crowley--Nordstr\"om $\nu$-invariants of the induced $\rG_2$-structures.
Theorem~\ref{cor:conj7} establishes these consequences as well as equality of the Hilbert functions in every degree.

\begin{introconjecture}[Calabi--Yau Bound (CYB); Conjecture 8 of \cite{AHHHSS2024}]
\label{conj:cyb}
In the Calabi--Yau context of Definition~\ref{def:intro-setting} with $n=3$, let $K_f$ and $X_f$ denote the link and the coarse hypersurface, respectively.
Then
\[
h^{2,1}_S(K_f)\le h^{2,1}_{CY},
\]
where $h^{2,1}_{CY}$ is the Calabi--Yau Hodge number recorded for $X_f$ in the Kreuzer--Skarke list, interpreted here as its stringy Hodge number $h^{2,1}_{st}(X_f)$.
\end{introconjecture}

Conjecture~\ref{conj:cyb} is numbered Conjecture~11 in the first arXiv version of \cite{AHHHSS2024}.
Equality was observed in $4198$ of the $7549$ computed cases.
The distinction between the ordinary orbifold and stringy Hodge numbers is explained in Table~\ref{tab:hodge-notation} and Section~\ref{sec:conj8}.

\paragraph{Weight determinacy and the sheaf-cohomological formula.}
We prove Conjectures~\ref{conj:wd} and~\ref{conj:cyb} and relate them to the classical descriptions of the invariants involved.
The algebraic input is the Hilbert-series calculation for a weighted complete intersection, underlying the Milnor--Orlik formula \cite{MilnorOrlik1970,Stanley1978}; a recent account of the Koszul-complex proof is given by Soler \cite{Soler2026}.
Conjecture~\ref{conj:wd} is a consequence of this classical calculation.
When $f$ has an isolated singularity its partial derivatives form a regular sequence, the condition on successive equations recalled in Section~\ref{sec:hilbert}.
The resulting resolution of the Milnor algebra gives, with the notation of Definition~\ref{def:intro-setting},
\begin{equation}
\label{eq:intro-series}
\sum_{\ell\ge0}\dim(\MM_f)_\ell\,t^\ell=\prod_{i=1}^N\frac{1-t^{\,d-w_i}}{1-t^{\,w_i}}\,.
\end{equation}
This gives the statement in every degree, not only in those relevant to the Hodge numbers, once the isolated-singularity hypothesis implicit in \cite{AHHHSS2024} is made explicit; Example~\ref{ex:nonisolated} shows that the hypothesis cannot be dropped.

\begin{maintheorem}[Weight determinacy]
\label{cor:conj7}
Let $f$ and $g$ be admissible polynomials in $N$ variables, in the sense of Definition~\ref{def:intro-setting}, of the same degree $d$ and with weight vectors agreeing up to permutation.
Then their graded Milnor algebras satisfy
\[
\dim(\MM_f)_\ell=\dim(\MM_g)_\ell\qquad(\ell\in\Z).
\]
In particular $f$ and $g$ are weakly R-equivalent in the sense of Definition~\ref{def:weak-r-equivalence}, and they have the same Sasakian Hodge numbers and the same Milnor number.
If moreover $N=5$ and $d=\sum_iw_i$, they have the same value of \eqref{eq:nuSt}, and hence the same Crowley--Nordstr\"om $\nu$-invariant of their induced $\rG_2$-structures.
\end{maintheorem}

Theorem~\ref{cor:conj7} also answers the question posed by Itoh \cite[p.~123, Problem following Cor.~1.3]{Itoh2004}: whether the Sasakian Hodge numbers depend only on the weights and degree, at least for five-dimensional links.
The answer holds in every dimension $2n+1\ge5$ considered here, without the Calabi--Yau condition.
At the algebraic level, \eqref{eq:intro-series} refines the Milnor--Orlik product formula by retaining the grading; evaluating at $t=1$ recovers the Milnor number \cite{MilnorOrlik1970,Soler2026}.

The next theorem computes the basic Hodge numbers of the link by sheaf cohomology on the leaf-space orbifold and reproduces the coefficients of \eqref{eq:intro-series}.

\begin{maintheorem}[Sasakian Hodge numbers by sheaf cohomology]
\label{thm:sheaf}
In the setting of Definition~\ref{def:intro-setting}, $h^{0,0}_S(K_f)=1$, $h^{p,q}_S(K_f)=0$ for $0<p+q<n$, and for $p+q=n$
\begin{equation}
\label{eq:sheaf}
h^{p,q}_S(K_f)=\sum_{\substack{I\subseteq\{1,\dots,N\}\\|I|\le p}}(-1)^{|I|}\,p_{\bw}\Big((p+1)d-\sum_iw_i-\sum_{i\in I}(d-w_i)\Big)=a_{(p+1)d-\sum_iw_i}.
\end{equation}
\end{maintheorem}

Combining Theorem~\ref{thm:sheaf} with the classical Hilbert-series identity \eqref{eq:intro-series} recovers Itoh's formula for the Sasakian Hodge numbers in terms of the Milnor algebra \cite[Thm.~1.2]{Itoh2004}.
Itoh's proof uses the residue-theoretic description of primitive cohomology due to Griffiths and Steenbrink \cite{Griffiths1969a,Griffiths1969b,Steenbrink1977}.
The contribution here is an alternative sheaf-cohomological derivation under the same assumptions, using neither residues nor the Jacobian-ring description of the Hodge structure.

In the Calabi--Yau case the outcome is a closed formula in terms of the numbers $p_{\bw}(k)$ of monomials of weighted degree $k$; for $7$-dimensional links, for instance,
\begin{equation}
\label{eq:intro-h21}
h^{2,1}_S=p_{\bw}(d)-\sum_{i=1}^5p_{\bw}(w_i)\,.
\end{equation}

\paragraph{The Calabi--Yau bound.}
Conjecture~\ref{conj:cyb} concerns the relation between ordinary orbifold cohomology and stringy Hodge theory.
The relevant framework is provided by Chen--Ruan cohomology \cite{ChenRuan2004}, Batyrev's stringy invariants \cite{Batyrev1998} and their comparison by Yasuda \cite{Yasuda2004}.
The proof compares the Hodge numbers in Table~\ref{tab:hodge-notation}; in particular, the ordinary Hodge numbers of the quotient orbifold $\cX_f$ must be distinguished from its Chen--Ruan Hodge numbers.
\begin{table}[htbp]
\caption{Hodge-number conventions for the link and its quotient.}
\label{tab:hodge-notation}
\centering
\begin{tabular}{@{}lp{0.67\textwidth}@{}}
\toprule
Notation & Meaning\\
\midrule
$h^{p,q}_B(K_f)=h^{p,q}(\cX_f)$ & Basic Hodge numbers of the link, equal to the ordinary Hodge numbers of the quotient orbifold.\\
$h^{p,q}_S(K_f)$ & Sasakian Hodge numbers, the primitive part of basic cohomology in degrees up to the middle.\\
$h^{p,q}_{\rm orb}(\cX_f)$ & Chen--Ruan Hodge numbers, including contributions from twisted sectors.\\
$h^{p,q}_{CY}=h^{p,q}_{st}(X_f)$ & Stringy Hodge numbers of the coarse hypersurface, equal to $h^{p,q}_{\rm orb}(\cX_f)$.\\
\bottomrule
\end{tabular}
\end{table}
The Kreuzer--Skarke lists record the stringy Hodge numbers.
In the middle degree, the comparison is
\begin{equation}
\label{eq:intro-chain}
h^{p,q}_{CY}=h^{p,q}_B+\big(\text{twisted contributions}\big),\qquad h^{p,q}_B=h^{p,q}_S+\delta_{pq}\qquad(p+q=n).
\end{equation}
The non-negativity of the twisted contributions gives the bound in every dimension $n\ge2$.
The following theorem identifies the gap in $h^{n-1,1}$ for $n\ge3$.
Here $\mu_r$ is the group of $r$-th roots of unity, and for an index set $F$ the notation $p_{\bw_F}$ counts monomials in the variables indexed by $F$, with their inherited weights.

\begin{maintheorem}[The Calabi--Yau bound and its defect]
\label{thm:defect}
Let $f$ be an admissible Calabi--Yau polynomial as in Definition~\ref{def:intro-setting}, and put $X\coloneq X_f$ and $h^{p,q}_{CY}\coloneq h^{p,q}_{st}(X)$.
For $F\subseteq\{1,\dots,N\}$ write $\bw_F\coloneq (w_i)_{i\in F}$ and $g_F\coloneq \gcd(w_i:i\in F)$, and set $p_{\bw_F}(k)\coloneq0$ for $k<0$.
Let $p+q=n$.
\begin{enumerate}
\item[(i)] For $n\ge2$, $h^{p,q}_S(K_f)\le h^{p,q}_{CY}-\delta_{pq}$.
\item[(ii)] For $n\ge3$,
\begin{equation}
\label{eq:defect}
h^{n-1,1}_{CY}-h^{n-1,1}_S(K_f)=\sum_{\substack{|F|=n\\ g_F>1}}(g_F-1)\;p_{\bw_F}\Big(d-\sum_{i\in F}w_i\Big).
\end{equation}
The summand indexed by $F$ is $g_F-1$ times the geometric genus $p_g=h^{n-2,0}$ of the closed suborbifold $\mathcal C_F=[(Y\cap\C^F)/\C^*]$ of points whose stabiliser contains $\mu_{g_F}$, whose coarse space is $X\cap\WP(\bw_F)$.
\end{enumerate}
\end{maintheorem}

The loci $\mathcal C_F$ in Theorem~\ref{thm:defect} have dimension $n-2$ and codimension two in $X_f$.
For threefolds they are curves, and the bound is attained exactly when all of them are rational.

\Needspace{4\baselineskip}
\begin{maintheorem}[The threefold bound]
\label{thm:conj8}
Under the assumptions and notation of Theorem~\ref{thm:defect}, let $n=3$, so that $X_f$ is a threefold and $K_f$ is seven-dimensional.
For each index set $F$ with $|F|=3$ and $g_F>1$, the coarse space $C_F\coloneq X_f\cap\WP(\bw_F)$ of $\mathcal C_F$ is a smooth connected projective curve.
Its generic points have stabiliser $\mu_{g_F}$ and transverse singularity of type $A_{g_F-1}$ in $X_f$: the normal slice is the cyclic quotient $\C^2/\mu_{g_F}$ with action $(u,v)\mapsto(\zeta u,\zeta^{-1}v)$.
Finitely many special points may have larger stabilisers.
With the geometric genus $p_g(C_F)=p_{\bw_F}(d-\sum_{i\in F}w_i)$ of the whole curve,
\[h^{2,1}_{CY}=h^{2,1}_S(K_f)+\sum_{|F|=3,\ g_F>1}(g_F-1)\,p_g(C_F)\ \ge\ h^{2,1}_S(K_f).\]
Equality holds if and only if every such curve is rational, in particular when there is none.
\end{maintheorem}

Theorems~\ref{thm:defect} and~\ref{thm:conj8} express the classical distinction between polynomial and non-polynomial deformations \cite{GreenHubsch1987,CdOK1995} in terms of the Sasakian Hodge numbers of the link.
They give a geometric interpretation and an exact error formula for the $h^{n-1,1}$ approximation of \cite[\S5]{Hirst:2023kdl}, recovering its lower-bound property in dimensions three and four and giving the corresponding statement for every $n\ge3$.
Remark~\ref{rem:classical} makes these comparisons explicit.
Equality holds for $4202$ of the $7555$ weight systems, and for $4198$ of the $7549$ considered in \cite{AHHHSS2024}.
The same decomposition explains why the Sasakian Hodge numbers are not mirror symmetric, as the link sees only the untwisted sector, while mirror symmetry exchanges untwisted with twisted contributions (Remark~\ref{rem:mirror}).

\paragraph{Computations for fourfolds.}
For fourfolds the same decomposition gives a Libgober--Wood identity, which explains the near-linear relation between the two Sasakian Hodge numbers seen in the data.
For the punctured cone $Y$ of Definition~\ref{def:intro-setting} and a finite-order element $t\in\C^*$, let $Y^t\coloneq \{z\in Y:t\cdot z=z\}$ and $\cX_t\coloneq [Y^t/\C^*]$.
The non-empty sectors with $t\ne1$ are the twisted sectors.
Their \emph{age} is the sum of the fractional arguments of the eigenvalues of $t$ on the tangent space of $\cX_f$, with arguments taken in $[0,1)$ (Lemma~\ref{lem:age}).

\Needspace{4\baselineskip}
\begin{maintheorem}[Libgober--Wood identity for nine-dimensional links]
\label{prop:lw}
Under the assumptions of Theorem~\ref{thm:defect}, let $n=4$ and use the sectors $\cX_t$ defined above.
The middle Sasakian Hodge numbers of the nine-dimensional link $K_f$ satisfy
\begin{equation}
\label{eq:link-lw}
h^{2,2}_S=4h^{3,1}_S+47+R,\qquad R\coloneq \sum_{t\neq1}\rho(t),
\end{equation}
where $t$ runs over the non-trivial sectors $\cX_t$ defined above.
If $\cX_t$ has age $a$ and dimension $m$, then $m\le2$, and its contribution is given by
\begin{center}
\begin{tabular}{ccc}
\toprule
Dimension $m$ & Age $a$ & Contribution $\rho(t)$\\
\midrule
$2$ & $1$ & $6\chi(\cO_{\cX_t})-e(\cX_t)$\\
$1$ & $1$ & $e(\cX_t)+1$\\
$1$ & $2$ & $-1$\\
$0$ & $1$ & $4c$\\
$0$ & $2$ & $-c$\\
$0$ & $3$ & $0$\\
\bottomrule
\end{tabular}
\end{center}
where $\chi(\cO_{\cX_t})$ is the holomorphic Euler characteristic, $e$ is the Euler characteristic of the coarse space and $c$ the number of points of a sector of dimension zero.
The sectors of dimension one occur in pairs $\{t,t^{-1}\}$ of ages $1$ and $2$, and each pair contributes $e(\cX_t)$.
In particular $h^{2,2}_S=4h^{3,1}_S+47$ when $X_f$ has no orbifold points.
\end{maintheorem}

Theorem~\ref{prop:lw} is a consequence of Batyrev's stringy Libgober--Wood identity \cite[Cor.~3.10, Cor.~3.11]{Batyrev2000MRL}, after separating the untwisted and twisted contributions.
It gives the resulting relation on the link together with the correction contributed by each sector.

The closed formulae express each Sasakian Hodge number considered here as an alternating sum of at most $22$ denumerants $p_{\bw}(k)$, which are computed exactly and on the order of milliseconds (far faster than Gr\"obner basis methods).
In the first place they complete the $7$-dimensional dataset of \cite{AHHHSS2024}, via the invariants of the six weight systems whose Gr\"obner bases could not be computed as given in Section~\ref{sec:six}, where they are also certified exactly for explicit polynomials.
They also give two further datasets.
The first consists of the $5$-dimensional links of the $95$ weighted K3 surfaces.
Through $b_2(K_f)=2+h^{1,1}_S(K_f)$, it recovers the Betti numbers computed by Boyer, Galicki and Matzeu \cite[Ex.~41, Cor.~42]{BGM2006} and gives a uniform Hodge-theoretic calculation with independent Gr\"obner-basis checks; the full list is given in Table~\ref{tab:k3}.
The second consists of the $9$-dimensional links of all $1\,100\,055$ weighted Calabi--Yau fourfolds \cite{LSW1999}.
The Gr\"obner-basis approach of \cite{AHHHSS2024} enumerates a monomial basis of the whole Milnor algebra, whose dimension reaches $3\times10^{12}$ in this dataset, and it originally did not complete under months of CPU compute for these links; with the closed formulae all of them are processed in about five minutes on a laptop.
The results are checked by four verification methods, including exact certificates for explicit polynomials and a comparison with the Kreuzer--Skarke data, and they are summarised in Tables~\ref{tab:9d-stats} and~\ref{tab:9d-examples}.

\paragraph{From experiment to proof.}
The passage from \cite{AHHHSS2024} to the present paper completes an experiment--conjecture--proof cycle in AI-assisted mathematics: computational data and their machine-learning-assisted interpretation suggested Conjectures~\ref{conj:wd} and~\ref{conj:cyb}, and the arguments below establish them deductively.
The experiments identify patterns to explain; the proofs supply hypotheses, mechanisms and conclusions that hold beyond the sampled cases.
This interaction has prominent precedents, including the machine-learning-guided discovery of a relation between geometric and algebraic knot invariants by Davies and collaborators \cite{Davies2021}.
This project style illustrates the complementary roles for computational discovery and subsequent mathematical justification; where AI agents can support the deduction of the rigorous mathematics behind numerical observations.

\paragraph{Organisation.}
This paper is organised as follows.
Section~\ref{sec:background} recalls the background, fixes the conventions for the orbifold and its coarse space, and recalls equivalence of singularities.
Section~\ref{sec:conj7} proves Theorems~\ref{cor:conj7} and~\ref{thm:sheaf} and derives the closed formulae.
Section~\ref{sec:conj8} proves Theorems~\ref{thm:defect}, \ref{thm:conj8}, and~\ref{prop:lw}.
Section~\ref{sec:applications} discusses the data for $5$- and $9$-dimensional links in Subsections~\ref{sec:5d} and~\ref{sec:9d}; the \hyperref[sec:conclusions]{afterword} returns to the experiments and further questions.
Appendix~\ref{app:comp} records the computational certificates, the computing environment and the complete K3 table.
The data, code and certificates are available in the corresponding \href{https://github.com/GeometrIALab/Link_SasakianHodgeNumbers}{GitHub} repository.

\newpage
\section{Background}
\label{sec:background}

We fix the weighted grading and the conventions for the link $K_f$, the quotient orbifold $\cX_f$, and its coarse space $X_f$.
The distinction between basic, Sasakian, and stringy Hodge numbers is essential for the comparison outlined in Section~\ref{sec:intro}.
Subsection~\ref{sec:algebra-background} states the algebraic facts used for the cone, and Subsection~\ref{sec:orbifold-background} assembles the quotient-singularity criterion, sector formula and stringy comparison needed in Section~\ref{sec:conj8}.
The Milnor-algebra description identifies the degrees used in Definition~\ref{def:weak-r-equivalence}; we also recall how weak R-equivalence differs from equivalence under a change of coordinates.

\subsection{Weighted homogeneous polynomials and Milnor algebras}
\label{sec:weights}
Throughout, $N=n+2\ge4$, and $\bw\coloneq (w_1,\dots,w_N)$ is a vector of positive integer weights with $\gcd(\bw)=1$; when weight systems are listed or counted, their weights are ordered nondecreasingly.
The polynomial ring $S\coloneq \C[z_1,\dots,z_N]$ is graded by $\deg z_i=w_i$, and $S_k$ denotes its degree-$k$ part.
The dimension of $S_k$ is the \emph{denumerant}
\begin{equation}
\label{eq:denumerant}
p_{\bw}(k)\coloneq\#\Big\{m\in\Z_{\ge0}^N:\ \sum_{i=1}^Nm_iw_i=k\Big\},\qquad p_{\bw}(k)\coloneq0\qforq k<0,
\end{equation}
which counts the monomials of weighted degree $k$.
A polynomial $f\in S_d$ is called \emph{$\bw$-homogeneous of degree $d$}; equivalently,
\begin{equation}
\label{eq:weighted-homogeneity}
f(t^{w_1}z_1,\dots,t^{w_N}z_N)=t^df(z)\qquad(t\in\C^*).
\end{equation}
Its \emph{Jacobian ideal} and graded \emph{Milnor algebra} are, respectively,
\begin{equation}
\label{eq:milnor-algebra}
\cJ_f\coloneq (\partial_1f,\dots,\partial_Nf),\qquad
\MM_f\coloneq S/\cJ_f=\bigoplus_{\ell\ge0}(\MM_f)_\ell.
\end{equation}

\paragraph{Admissibility and quasi-smoothness.}
Recall from Definition~\ref{def:intro-setting} that $f$ is \emph{admissible} if it has an isolated critical point at the origin; write $U_{\bw,d}\subseteq S_d$ for the set of such polynomials of degree $d$.
By Lemma~\ref{lem:VJ} below, $f$ is admissible exactly when the common zero set of its partial derivatives is $\{0\}$.
Euler's identity is
\begin{equation}
\label{eq:euler}
d\,f=\sum_iw_iz_i\partial_if.
\end{equation}
Thus every critical point of $f$ lies on $\{f=0\}$, so this means that the affine cone $\{f=0\}\subset\C^N$ is smooth away from the origin, that is, that the hypersurface $X_f$ is \emph{quasi-smooth} \cite[\S3.1]{Dolgachev1982}.
Quasi-smoothness concerns the punctured affine cone; the coarse hypersurface can still have quotient singularities.
In particular, the isolated hypersurface singularity at the cone vertex and the cyclic quotient singularities of $X_f$ are different objects.
The pair $(\bw,d)$ is called \emph{non-degenerate} when $U_{\bw,d}\neq\emptyset$, and \emph{Calabi--Yau} when moreover $d=\sum_iw_i$.
The non-degenerate Calabi--Yau weight systems are classified for small $n$: there are $95$ for $n=2$ \cite{Yonemura1990,IanoFletcher2000}, $7555$ for $n=3$ \cite{CandelasLynkerSchimmrigk1990,KlemmSchimmrigk1994,KreuzerSkarke1992} and $1\,100\,055$ for $n=4$ \cite{LSW1999}, all available from \cite{KSwebsite,KSfourfolds}.

\paragraph{Hilbert series and monomial counts.}
The degree $d$ is fixed throughout and is suppressed from the notation $P_{\bw}$ and $\mu(\bw)$ below; in the Calabi--Yau case it is determined by $\bw$, and $U_{\bw}=U_{\bw,d}$.
The following power series and its value at $t=1$ will appear repeatedly:
\begin{equation}
\label{eq:Pw}
P_{\bw}(t)\coloneq \prod_{i=1}^N\frac{1-t^{\,d-w_i}}{1-t^{\,w_i}}=\sum_{k\ge0}a_kt^k,\qquad\mu(\bw)\coloneq \lim_{t\to1}P_{\bw}(t)=\prod_{i=1}^N\Big(\frac d{w_i}-1\Big).
\end{equation}
The limit exists factor by factor, and Theorem~\ref{thm:hilbert} shows that $P_{\bw}$ is a polynomial when $(\bw,d)$ is non-degenerate.
Expanding the numerator and the geometric series in the denominator expresses the coefficients, for every $(\bw,d)$, as
\begin{equation}
\label{eq:incl-excl}
a_\ell=\sum_{I\subseteq\{1,\dots,N\}}(-1)^{|I|}\,p_{\bw}\Big(\ell-\sum_{i\in I}(d-w_i)\Big),
\end{equation}
an identity of formal power series that involves no polynomial $f$.
The \emph{Hilbert function} of $\MM_f$ is $\ell\mapsto\dim(\MM_f)_\ell$, and its \emph{Hilbert series} is the generating series of these dimensions.
Theorem~\ref{thm:hilbert} identifies this series with $P_{\bw}$.
For a chosen monomial order, a \emph{Gr\"obner basis} of $\cJ_f$ permits polynomial division with a unique remainder.
The monomials not divisible by any leading monomial of the basis are called \emph{standard monomials}; their residue classes form a vector-space basis of $\MM_f$ \cite[Chs.~2,5]{CoxLittleOShea2015}.
For the Fermat quintic in Example~\ref{ex:nonisolated}, these are
\[
z_1^{m_1}\cdots z_5^{m_5},\qquad 0\le m_i\le3,
\]
because $\cJ_f=(z_1^4,\dots,z_5^4)$.
Counting standard monomials degree by degree gives the Hilbert function.
The product formula computes the same counts without constructing a basis.

Finally, $\mu_r\subset\C^*$ denotes the group of $r$-th roots of unity, empty products and least common multiples equal $1$, and Hodge numbers with an index outside the allowed range are $0$.

\subsection{Commutative algebra for the cone}
\label{sec:algebra-background}
The algebraic results below replace transversality at the singular vertex.
We state the versions needed later, keeping their geometric interpretation in view.

\begin{remark}[regular sequences and the Cohen--Macaulay property]
\label{rem:cohen-macaulay}
A sequence $q_1,\dots,q_N$ is \emph{regular} in $S$ if each $q_i$ is a non-zero-divisor in $S/(q_1,\dots,q_{i-1})$ and the final quotient is nonzero.
Thus multiplication by the next equation is injective after imposing the earlier ones.
A homogeneous \emph{system of parameters} consists here of $N$ positive-degree elements whose quotient is finite-dimensional.
A polynomial ring is \emph{Cohen--Macaulay}: its local \emph{depth}, the maximal length of a regular sequence in the maximal ideal, equals its dimension.
The consequence used below is that every homogeneous system of parameters is regular \cite[\S\S1.1,2.1]{BrunsHerzog1993}.
The quotient by a regular sequence is a \emph{complete intersection}.
This condition makes the Koszul complex exact, allowing the Hilbert series to be computed by an alternating sum.
Geometrically, a regular sequence imposes successive equations with the expected codimension, but does not require their zero sets to be smooth.
This distinction matters at the cone vertex: the partial derivatives all vanish there, and the quotient retains their algebraic multiplicities.
The Cohen--Macaulay property allows the dimension condition to control these equations through an exact complex, even though the differential-geometric regular-value theorem does not apply.
\end{remark}

\begin{example}[regularity without transversality]
\label{ex:regular-sequence}
In $R\coloneq\C[x,y]$ with the ordinary grading, the sequence $x^2,y^3$ is regular.
Indeed, $x^2$ is a non-zero-divisor in $R$, and multiplication by $y^3$ is injective in $R/(x^2)$: write each class uniquely as $a(y)+xb(y)$.
Nevertheless, the differentials of both equations vanish at the origin.
The quotient has basis $1,x,y,xy,y^2,xy^2$, so its Hilbert series is
\[
1+2t+2t^2+t^3=(1+t)(1+t+t^2)
=\frac{(1-t^2)(1-t^3)}{(1-t)^2}.
\]
This is the product mechanism of Theorem~\ref{thm:hilbert} in two variables.
By contrast, $x^2,xy$ is not regular: the nonzero class of $x$ in $R/(x^2)$ annihilates $xy$, and the common zero set is a line.
\end{example}

\begin{proposition}[algebraic inputs]
\label{prop:algebra-inputs}
Let $S=\C[z_1,\dots,z_N]$ with positive weights, and let $\fm=(z_1,\dots,z_N)$.
The following facts will be used for the Jacobian ideal and the coordinate ring of the cone.
\begin{enumerate}
\item[(i)] If an ideal of $S$ is generated by $c$ elements, every irreducible component of its zero set has dimension at least $N-c$ (Krull's height theorem).
\item[(ii)] The ring $S$ is Cohen--Macaulay. If homogeneous elements $q_1,\dots,q_N$ of positive degrees have finite-dimensional quotient, they form a regular sequence. Its Koszul complex is a free resolution of $S/(q_1,\dots,q_N)$.
\item[(iii)] Let $A$ be a finitely generated positively graded Cohen--Macaulay $\C$-algebra of dimension $m$, with $A_0=\C$ and homogeneous maximal ideal $\fm_A=A_{>0}$. Its quotient by a homogeneous non-zero-divisor of positive degree is Cohen--Macaulay of dimension $m-1$. Moreover,
\[
H^j_{\fm_A}(A)=0\qforq j\ne m.
\]
\end{enumerate}
These are the height, parameter, Koszul-exactness and depth-vanishing results of \cite[\S\S1.1,1.2,1.6,2.1 and Ch.~3]{BrunsHerzog1993}.
\end{proposition}

Here $H^j_{\fm_A}(A)$ is \emph{local cohomology supported at the vertex}: it is the cohomology of the augmented \v{C}ech complex formed by inverting homogeneous generators of $\fm_A$.
The explicit complex and its comparison with the punctured cone are written out before Lemma~\ref{lem:lines}.
Part~(iii) says that a Cohen--Macaulay cone has only one potentially nonzero local-cohomology degree; this is the source of the line-bundle vanishings in that lemma.
Part~(i) will also exclude overly large fixed loci in Lemma~\ref{lem:twisted}.

For a finite-dimensional connected graded algebra $R$ ($R_0=\C$), the \emph{socle} is the subspace annihilated by $R_{>0}$.
Such an algebra is \emph{Artinian Gorenstein} precisely when its socle is one-dimensional.

\begin{proposition}[duality for a graded complete intersection]
\label{prop:ci-duality}
Let $q_1,\dots,q_N$ be a homogeneous regular sequence in $S$, with $\deg q_i=e_i>0$.
Then $R\coloneq S/(q_1,\dots,q_N)$ is Artinian Gorenstein, its socle lies in degree
\[
s\coloneq\sum_i e_i-\sum_iw_i,
\]
and multiplication gives perfect pairings
\[
R_k\otimes R_{s-k}\longrightarrow R_s\cong\C.
\]
In particular $\dim R_k=\dim R_{s-k}$; see \cite[\S4, especially Thm.~4.1 and its proof]{Stanley1978} and \cite[\S3.1]{BrunsHerzog1993}.
\end{proposition}

``Perfect'' means that every nonzero element has a partner with nonzero product in the socle, or equivalently $R_k\cong R_{s-k}^*$ after choosing a generator of $R_s$.
This is an algebraic analogue of a top-degree pairing.
For the Milnor algebra the degrees are $e_i=d-w_i$, giving the symmetry used in Corollary~\ref{cor:closed}.

\subsection{Links and Sasakian Hodge numbers}
\label{sec:links}
For $f\in U_{\bw,d}$ let $V_f\coloneq f^{-1}(0)\subset\C^N$ and $Y\coloneq V_f\setminus0$.
The link $K_f\coloneq V_f\cap S^{2N-1}$ is a smooth closed manifold of dimension $2n+1$, and it is $(n-1)$-connected \cite{Milnor1968}; its Betti numbers are $b_k(K_f)\coloneq \dim_\Q H^k(K_f;\Q)$.
The local smoothing of the cone is the \emph{Milnor fibre}
\[
W_f\coloneq f^{-1}(\delta)\cap B_\epsilon^{2N},\qquad 0<|\delta|\ll\epsilon\ll1.
\]
Its boundary is diffeomorphic to the link, and $W_f$ has the homotopy type of a wedge of $(N-1)$-spheres \cite{Milnor1968}.
Their number, the \emph{Milnor number}, is $\dim_\C\MM_f$; thus the finite-dimensional algebra records the middle topology of the smoothing.
Here the polynomial Jacobian quotient agrees with the local analytic quotient, since it is supported only at the origin.
Moving $\delta$ once around zero gives the \emph{monodromy} of the Milnor fibre, whose action on cohomology enters the Milnor--Orlik formula used in Subsection~\ref{sec:9d} \cite{MilnorOrlik1970}.

\paragraph{The quotient orbifold and its coarse space.}
The weighted action of $\C^*$ on $\C^N$ and its Euler vector field are
\begin{equation}
\label{eq:weighted-action}
t\cdot z\coloneq (t^{w_1}z_1,\dots,t^{w_N}z_N),\qquad E_{\bw}\coloneq \sum_iw_iz_i\partial_{z_i}.
\end{equation}
The action preserves $Y$, its circle subgroup preserves $K_f$, and its infinitesimal action is generated by $E_{\bw}$.
The stabiliser of $z\in Y$ is $\mu_{m(z)}$ with $m(z)\coloneq \gcd\{w_i:z_i\ne0\}$; since $\gcd(\bw)=1$ the generic stabiliser is trivial.

The orbifold and its coarse space are distinguished throughout.
The quotient stack $\cX_f\coloneq [Y/\C^*]$ is a smooth proper Deligne--Mumford stack, that is, a compact complex orbifold, of dimension $n$.
Concretely, a transverse slice to an orbit gives a smooth chart $U$ with a finite cyclic stabiliser group $G$.
The orbifold retains the chart and its group action, whereas its coarse space retains only the quotient $U/G$; see \cite[\S4]{Yasuda2004} for the stack terminology.
A local action is \emph{effective} when only the identity acts trivially.
We allow ineffective stabilisers on the sectors below, and prove effectivity of $\cX_f$ in the Calabi--Yau case in Section~\ref{sec:conj8}.

Coherent sheaves on $\cX_f$ correspond to coherent sheaves on $Y$ with a compatible $\C^*$-action; this is the meaning of \emph{equivariant} here.
We write $\Omega^p_{\cX_f}$ for the sheaf of holomorphic $p$-forms.
The orbifold line bundle $\cO_{\cX_f}(k)$ is specified by requiring its sections to lift to functions $\phi$ on $Y$ satisfying
\[
\phi(t\cdot z)=t^k\phi(z).
\]
Thus weighted homogeneous polynomials define sections.
These are line bundles on the orbifold; their pushforwards to the coarse space need not be locally free, as already occurs in weighted projective space \cite[\S\S1.4--1.5]{Dolgachev1982}.
The coarse space $X_f\coloneq Y/\C^*\subset\WP^{N-1}(\bw)$ is a projective variety with cyclic quotient singularities.
The Hodge numbers of the orbifold are
\begin{equation}
\label{eq:orbifold-ordinary-hodge}
h^{p,q}(\cX_f)\coloneq \dim H^q(\cX_f,\Omega^p_{\cX_f}).
\end{equation}
Locally, pushing a sheaf to $X_f$ means taking the invariant sections in a chart $[U/G]$.
Over $\C$, averaging over the finite group $G$ makes this operation exact: an invariant section in a quotient has an invariant lift.
Consequently coherent cohomology on $\cX_f$ is the cohomology of its pushforward to $X_f$.
Serre's GAGA comparison theorem identifies algebraic and analytic coherent cohomology on a complex projective variety \cite[\S3]{Serre1956}.
It therefore allows us to compute the holomorphic cohomology in \eqref{eq:orbifold-ordinary-hodge} algebraically.

\begin{theorem}[Hodge decomposition for projective V-manifolds]
\label{thm:vmanifold-hodge}
Let $\mathcal Z$ be a smooth proper complex orbifold with projective coarse space $Z$.
The ordinary rational cohomology of $Z$ is pure, and
\[
H^k(Z;\C)\cong\bigoplus_{p+q=k}H^q(\mathcal Z,\Omega^p_{\mathcal Z}).
\]
If a finite generic stabiliser acts trivially, removing it does not alter these groups.
For weighted projective space of dimension $m$, in particular,
\[
h^{p,q}(\WP(\mathbf v))=
\begin{cases}1,&0\le p=q\le m,\\0,&\text{otherwise}.\end{cases}
\]
The decomposition is Steenbrink's V-manifold Hodge theory \cite[\S1]{SteenbrinkOslo1977}, \cite[\S2]{Steenbrink1977}; the weighted projective calculation is \cite[Cor.~2.3.6]{Dolgachev1982}.
\end{theorem}

A \emph{V-manifold} is a complex space locally given by finite quotients of smooth charts.
Its sheaf $\tilde\Omega^p_Z$ consists of invariant holomorphic forms in these charts, and is the pushforward of $\Omega^p_{\mathcal Z}$.
Averaging over the finite groups explains the cohomology comparison, and a group acting trivially on a chart acts trivially on its forms.
Thus the theorem applies to $\cX_f$ and to the closed suborbifolds used later; its Hodge numbers are the ordinary ones of their coarse spaces.
The age-shifted orbifold Hodge numbers introduced in Subsection~\ref{sec:orbifold-background} are a different invariant.

\paragraph{Sasakian and contact Calabi--Yau structures.}
A Sasakian structure on a manifold $M$ of dimension $2n+1$ is a tuple $(\xi,\eta,\Phi,g)$ consisting of a contact form $\eta$ with Reeb field $\xi$, an endomorphism $\Phi$ of $TM$ and a Riemannian metric $g$, such that the metric cone $(M\times\R_{>0},dr^2+r^2g)$ is K\"ahler with a complex structure that maps $r\partial_r$ to $\xi$ and restricts to $\Phi$ on $\ker\eta$ \cite[Ch.~6]{BoyerGalicki2008}.
It is \emph{quasi-regular} if all orbits of $\xi$ are closed; then $\xi$ generates a locally free circle action, and the space of leaves is a compact K\"ahler orbifold.
The weighted Sasakian structure of $S^{2N-1}$, whose Reeb field generates the weighted circle action, restricts to a quasi-regular Sasakian structure $(\xi,\eta,\Phi,g)$ on $K_f$ \cite[Ch.~9]{BoyerGalicki2008}.
Since the real weighted dilations move $Y$ transversally to the sphere, its space of leaves is $\cX_f$, and $K_f\to\cX_f$ is a Seifert circle orbibundle.
In what follows the structure is abbreviated to $(\xi,\eta)$.

In the Calabi--Yau case $d=\sum_iw_i$ the orbifold $\cX_f$ has a holomorphic volume form.
The residue 
\[
\Omega_Y\coloneq \mathrm{Res}_{f=0}\big(dz_1\wedge\dots\wedge dz_N/f\big)
\]
 is a nowhere-vanishing holomorphic $(n+1,0)$-form on $Y$ of weight $\sum_iw_i-d=0$, so that $\cL_{E_{\bw}}\Omega_Y=0$.
Its contraction $\varepsilon\coloneq \iota_{E_{\bw}}\Omega_Y$ is a holomorphic $(n,0)$-form with $\iota_{E_{\bw}}\varepsilon=0$, and it is closed by Cartan's formula:
\[
d\varepsilon=\cL_{E_{\bw}}\Omega_Y-\iota_{E_{\bw}}d\Omega_Y=0.
\]
Being invariant and horizontal, $\varepsilon$ descends to a nowhere-vanishing section of the canonical bundle of $\cX_f$, and on $K_f$ it restricts to a closed basic $(n,0)$-form that is nowhere vanishing on $\ker\eta$.
Following \cite[Def.~3.1]{TomassiniVezzoni2008}, a \emph{contact Calabi--Yau} structure is a Sasakian structure together with a closed basic $(n,0)$-form $\varepsilon$ with
\[\varepsilon\wedge\bar\varepsilon=(-1)^{n(n+1)/2}\,i^n\,(d\eta)^n;\]
such a structure is null Sasakian and transversally Ricci-flat \cite[Cor.~3.7]{TomassiniVezzoni2008}.
The structure restricted from the weighted sphere is null Sasakian, but its transverse metric is in general not Ricci-flat.
By El Kacimi's transverse Calabi--Yau theorem \cite{ElKacimi1990} it can be deformed to $\eta'=\eta+\zeta$, with $\zeta$ a basic $1$-form, so that the Reeb field, the transverse holomorphic structure and the basic class $[d\eta]$ are unchanged and the transverse metric is Ricci-flat \cite[Prop.~6.6, Prop.~6.7]{HabibVezzoni2015}, \cite[Prop.~13]{CalvoAndradeRodriguezSaEarp2020}.
For the deformed structure $\log|\varepsilon|^2$ is a basic pluriharmonic function, hence constant, so after rescaling $\varepsilon$ by a positive constant $K_f$ is contact Calabi--Yau.
The Hodge-theoretic invariants studied below depend only on the Reeb foliation, its transverse holomorphic structure and the class $[d\eta]$, so they are the same for both Sasakian structures.

\paragraph{Basic and Sasakian Hodge numbers.}
A differential form $\alpha$ on $K_f$ is \emph{basic} if $\iota_\xi\alpha=0$ and $\cL_\xi\alpha=0$.
Basic forms make up a subcomplex of the de Rham complex, whose cohomology with real coefficients is the basic cohomology $H^*_B(K_f)$.
The transverse complex structure splits complex-valued basic forms by type, and the resulting basic Dolbeault cohomology groups $H^{p,q}_B(K_f)$ have dimensions $h^{p,q}_B$.
The form $d\eta$ is a closed basic $(1,1)$-form, the transverse K\"ahler form, and $L\coloneq [d\eta]\wedge\,\cdot\,$ denotes the corresponding Lefschetz operator.
A class in $H^k_B$, or in $H^{p,q}_B$ with $p+q=k$, is \emph{primitive} if $k\le n$ and it is annihilated by $L^{n-k+1}$.

\begin{definition}[Sasakian Hodge numbers]
\label{def:sasakian}
On a compact Sasakian $(2n+1)$-manifold $M$ the harmonic $k$-forms with $k\le n$ are basic and primitive \cite{Tachibana1965}, and they split by type \cite{Tanaka1975,Itoh2004}.
The \emph{Sasakian Hodge number} $h^{p,q}_S(M)$, for $p+q\le n$, is the dimension of the space of harmonic forms of type $(p,q)$.
By basic Hodge theory for transversely K\"ahler foliations \cite[Prop.~3.4.5, Thm.~3.4.6, \S3.4.7]{ElKacimi1990}, this space is isomorphic to the primitive part of $H^{p,q}_B(M)$, so that
\begin{equation}
\label{eq:hS-basic}
h^{p,q}_S(M)=h^{p,q}_B(M)-h^{p-1,q-1}_B(M)\qquad(p+q\le n).
\end{equation}
\end{definition}

The following description underlies the computations of \cite{AHHHSS2024}.

\begin{theorem}[\cite{Steenbrink1977,Itoh2004}]
\label{thm:itoh}
Let $f$ be a weighted homogeneous polynomial of degree $d$ in $N=n+2$ variables, with an isolated critical point at the origin.
Write $K_f$ for its link and $\MM_f\coloneq S/\cJ_f$ for its graded Milnor algebra, as in Section~\ref{sec:weights}, and let $p+q=n$.
Then
\begin{equation}
\label{eq:itoh}
h^{p,q}_S(K_f)=\dim(\MM_f)_\ell,\qquad \ell\coloneq (p+1)d-\sum_iw_i.
\end{equation}
In the Calabi--Yau case $\ell=pd$.
\end{theorem}

Itoh's proof identifies the harmonic forms on $K_f$ with primitive harmonic forms on the orbifold $\cX_f$ \cite{Itoh2004}, and then appeals to the Jacobian-ring description of the primitive cohomology of quasi-smooth hypersurfaces \cite{Griffiths1969a,Griffiths1969b,Steenbrink1977,Dolgachev1982}.
In the convention of \cite[Thm.~4.3.2]{Dolgachev1982} the space $H^{p,q}_{\rm prim}(X_f)$ corresponds to the degree $(q+1)d-\sum_iw_i$, with $p$ and $q$ exchanged; both conventions give the same numbers by Hodge symmetry.
The alternative derivation in Theorem~\ref{thm:sheaf} is proved in Section~\ref{sec:sheaf}.

\subsection{Quotient singularities and stringy Hodge theory}
\label{sec:orbifold-background}
The proof of the Calabi--Yau Bound uses two descriptions of the same correction to ordinary Hodge numbers: stabiliser sectors of the orbifold, and exceptional strata of a resolution of its coarse space.
We recall the definitions and the comparison theorem here; Section~\ref{sec:conj8} checks their hypotheses and computes the sectors from the weights.
Unless stated otherwise, we use $\cX\coloneq\cX_f$, $X\coloneq X_f$ and $Y=f^{-1}(0)\setminus0$, with $f\in U_{\bw}$ and $d=\sum_iw_i$.

\paragraph{Stabilisers, sectors and age.}
For $t\in\C^*$ of finite order put
\begin{equation}
\label{eq:fixed-loci}
S(t)\coloneq \{i:t^{w_i}=1\},\qquad Y^t\coloneq \{z\in Y:t\cdot z=z\}=Y\cap\C^{S(t)},
\end{equation}
and write $\bw_S\coloneq (w_i)_{i\in S}$ for $S\subseteq\{1,\dots,N\}$.
The \emph{inertia orbifold} records a point together with an element of its stabiliser, so the same geometric locus may occur once for each such element \cite[\S3.1]{ChenRuan2004}.
Since $\C^*$ is abelian, the part labelled by $t$ is the \emph{sector} $\cX_t\coloneq [Y^t/\C^*]$.
It is a closed suborbifold of $\cX$, and the inertia orbifold is the disjoint union of these sectors over the $t$ of finite order with $Y^t\ne\emptyset$.
The coarse space of $\cX_t$ is the locus $$X^t\coloneq X\cap\WP(\bw_{S(t)}),$$ but $\cX_t$ keeps its stabilisers.
In particular, if $e=\gcd(\bw_{S(t)})>1$, then $\mu_e$ acts trivially on $Y^t$, and $\cX_t$ is a $\mu_e$-gerbe over the quotient of $Y^t$ by $\C^*/\mu_e\cong\C^*$, which acts with the weights $\bw_{S(t)}/e$.
Here \emph{gerbe} records the retained subgroup $\mu_e$ acting trivially on the entire fixed locus; passing to the indicated quotient removes that subgroup.
It does not remove the possible larger stabilisers at special points.

For a finite-order linear transformation $g$ of a complex tangent space, write its eigenvalues as $e^{2\pi i\vartheta_j}$ with $0\le\vartheta_j<1$.
Its \emph{age} is
\[
\age(g)\coloneq\sum_j\vartheta_j.
\]
It is locally constant on each connected fixed sector.
The invariant directions have $\vartheta_j=0$; in each moving direction the angles for $g$ and $g^{-1}$ add to $1$.
Consequently
\begin{equation}
\label{eq:age-codimension}
\age(g)+\age(g^{-1})=\operatorname{codim}_{\C}\operatorname{Fix}(g),
\qquad \det g=e^{2\pi i\age(g)}.
\end{equation}
These are the degree-shifting identities of \cite[Lem.~3.2.1]{ChenRuan2004}.
In particular, an effective determinant-one action has positive integral age on every nonidentity element; if its fixed locus has codimension two, both inverse ages must be $1$.
Lemma~\ref{lem:age} will compute these angles directly for a weighted hypersurface.

\paragraph{Canonical quotient singularities.}
We recall the singularity terminology needed to compare the two Hodge theories.
A \emph{cyclic quotient singularity} is locally the quotient of a smooth chart by a finite cyclic linear action.
In this setting, \emph{Gorenstein} means that the canonical sheaf is a line bundle; the determinant-one condition ensures that a local holomorphic volume form descends.
When the normal variety $X$ has Cartier canonical divisor, as in the Gorenstein case, a resolution $\pi:\widetilde X\to X$ gives
\[
K_{\widetilde X}=\pi^*K_X+\sum_j b_jD_j
\]
along the exceptional divisors $D_j$.
Here $K$ denotes the canonical divisor class, and the \emph{exceptional divisors} are the hypersurfaces contracted by $\pi$.
The coefficients $b_j$, called \emph{discrepancies}, measure the change of canonical class under resolution.
The singularities are \emph{canonical} when these discrepancies are non-negative on every resolution.
A resolution is \emph{crepant} when all $b_j=0$, so it preserves the canonical class.

\begin{theorem}[Reid--Tai criterion and determinant-one quotients]
\label{thm:quotient-criterion}
Let $G\subset\mathrm{GL}(m,\C)$ be finite and contain no \emph{quasireflections}, that is, no nonidentity element fixing a hyperplane pointwise.
Then $\C^m/G$ has canonical singularities if and only if $\age(g)\ge1$ for every $g\ne1$ in $G$ \cite[Thm.~4.11]{Reid1987}.
If $G\subset\mathrm{SL}(m,\C)$, its invariant ring is Gorenstein \cite[Cor.~5.6]{Stanley1978}; the quotient is therefore Gorenstein canonical.
\end{theorem}

For the last assertion, a finite-order matrix is diagonalisable, and a nonidentity determinant-one matrix cannot be a quasireflection.
Its age is a positive integer by \eqref{eq:age-codimension}, so the criterion applies.
This explains how a calculation of tangent-space eigenvalues supplies the singularity hypotheses required below.

\begin{example}[the cyclic surface quotient]
\label{ex:cyclic-surface}
Let $\zeta$ be a primitive $r$-th root of unity, $r\ge2$, acting by
\[
(u,v)\longmapsto(\zeta u,\zeta^{-1}v).
\]
The invariant ring is generated by $U=u^r$, $V=v^r$, $Z=uv$, with relation $UV=Z^r$.
Indeed, an invariant monomial $u^av^b$ has $a-b$ divisible by $r$; removing $(uv)^{\min(a,b)}$ leaves a power of $u^r$ or $v^r$.
This is the surface singularity of type $A_{r-1}$.
Its minimal resolution is crepant and has a chain of $r-1$ smooth rational curves of self-intersection $-2$ \cite[\S2]{Reid1987}.
The $r-1$ nonidentity group elements all have age $k/r+(r-k)/r=1$.
Thus the number of age-one sectors equals the number of exceptional curves; the transverse version of this example explains the multiplicity in Theorem~\ref{thm:conj8}.
\end{example}

\paragraph{The sector sum.}
Assume now that all local tangent actions have determinant one, as will be proved for our Calabi--Yau hypersurfaces in Section~\ref{sec:conj8}.
Then the ages are integers.
Following Chen and Ruan \cite[Def.~3.2.4]{ChenRuan2004}, the \emph{orbifold Hodge numbers} of $\cX$ are
\begin{equation}
\label{eq:orbifold}
h^{p,q}_{\rm orb}(\cX)\coloneq \sum_th^{\,p-\age(t),\,q-\age(t)}(\cX_t),
\end{equation}
where $t$ runs over the elements of finite order with $Y^t\ne\emptyset$, and $h^{i,j}(\cX_t)\coloneq \dim H^j(\cX_t,\Omega^i_{\cX_t})$ are the ordinary Hodge numbers of the sector, which by Theorem~\ref{thm:vmanifold-hodge} are those of its coarse space $X^t$.
A sector that is disconnected contributes the sum over its components, on which the age is the same by Lemma~\ref{lem:age}.
The term $t=1$ is the \emph{untwisted sector} $\cX_1=\cX$, and the others are the \emph{twisted sectors}.
The gerbe structure of a sector does not change its Hodge numbers, since the generic stabiliser acts trivially on its differential forms; how often a given locus occurs is recorded only by the number of elements $t$ with that fixed locus.

An age-one sector contributes its $(p-1,q-1)$-classes to bidegree $(p,q)$.
The shift is determined by the tangent action normal to the fixed locus.

\paragraph{Stringy invariants and the comparison theorem.}
To specify the invariant on the coarse space independently, let $X$ be a normal projective variety with Gorenstein canonical singularities and choose a resolution whose exceptional divisors $D_1,\dots,D_s$ have \emph{simple normal crossings}: locally they are coordinate hyperplanes.
Write their discrepancies as $b_j$, and set
\[
D_J\coloneq\bigcap_{j\in J}D_j,\qquad
D_J^\circ\coloneq D_J\setminus\bigcup_{j\notin J}D_j,
\qandq D_\emptyset\coloneq\widetilde X.
\]
Batyrev defines
\begin{equation}
\label{eq:stringy-resolution}
E_{st}(X;u,v)\coloneq
\sum_{J\subseteq\{1,\dots,s\}}E(D_J^\circ;u,v)
\prod_{j\in J}\frac{uv-1}{(uv)^{b_j+1}-1}
\end{equation}
\cite[Def.~3.1]{Batyrev1998}.
Here the compactly supported Hodge--Deligne polynomial is
\[
E(Z;u,v)\coloneq\sum_{k,p,q}(-1)^k h^{p,q}(H_c^k(Z;\C))u^pv^q;
\]
compact supports account for the possibly noncompact strata, and the Hodge numbers refer to the mixed Hodge structure on $H_c^k$.
For a smooth projective $Z$ it reduces to $\sum_{p,q}(-1)^{p+q}h^{p,q}(Z)u^pv^q$.
The expression \eqref{eq:stringy-resolution} is independent of the resolution \cite[Thm.~3.4]{Batyrev1998}.
When it is a polynomial, its signed coefficients define the \emph{stringy Hodge numbers} \cite[Def.~3.8]{Batyrev1998}:
\[
E_{st}(X;u,v)=\sum_{p,q}(-1)^{p+q}h^{p,q}_{st}(X)u^pv^q.
\]

\begin{theorem}[orbifold--stringy comparison and crepant resolutions]
\label{thm:orbifold-stringy}
Let $\mathcal Z$ be a smooth proper effective complex algebraic orbifold with projective coarse space $Z$, whose local tangent actions lie in $\mathrm{SL}(m,\C)$, where $m=\dim\mathcal Z$.
For each connected component $\mathcal T$ of its inertia orbifold, let $a_{\mathcal T}$ be its age.
Then the stringy $E$-function of $Z$ is a polynomial, and
\begin{equation}
\label{eq:stringy-sector-comparison}
h^{p,q}_{st}(Z)=\sum_{\mathcal T}h^{p-a_{\mathcal T},q-a_{\mathcal T}}(\mathcal T)
=h^{p,q}_{\rm orb}(\mathcal Z).
\end{equation}
This is the Hodge-theoretic consequence of \cite[Thm.~3.15 and proof of Cor.~3.16]{Yasuda2004}.
If $Z$ admits a projective crepant resolution $\widetilde Z$, then also
$h^{p,q}_{st}(Z)=h^{p,q}(\widetilde Z)$ \cite[Thm.~3.12]{Batyrev1998}.
\end{theorem}

Yasuda uses the age of the inverse stabiliser element as his shift number.
Inversion pairs the corresponding fixed loci, so reindexing gives precisely the convention in \eqref{eq:stringy-sector-comparison}.
The first equality requires no crepant resolution.
It also makes non-negativity transparent in this setting: every summand is the dimension of an ordinary Hodge space.
When a crepant resolution does exist, all factors in \eqref{eq:stringy-resolution} are $1$, and additivity over its strata explains the last equality.
The stringy Hodge numbers have the familiar symmetries
\[
h^{p,q}_{st}(Z)=h^{q,p}_{st}(Z)=h^{m-p,m-q}_{st}(Z)
\]
\cite[Thm.~3.7]{Batyrev1998}; these are also obtained from sectorwise Hodge symmetry and duality, pairing inverse sectors using \eqref{eq:age-codimension}.

\begin{theorem}[stringy Libgober--Wood identity, Calabi--Yau case]
\label{thm:background-lw}
Let $Z$ be a projective $m$-dimensional variety with Gorenstein canonical singularities, numerically trivial canonical class and polynomial stringy $E$-function.
Then \cite[Cor.~3.10, Cor.~3.11]{Batyrev2000MRL}
\begin{equation}
\label{eq:background-lw}
\sum_{p,q}(-1)^{p+q}\Big(p-\frac m2\Big)^2h^{p,q}_{st}(Z)
=\frac m{12}e_{st}(Z),\qquad
e_{st}(Z)\coloneq\sum_{p,q}(-1)^{p+q}h^{p,q}_{st}(Z).
\end{equation}
\end{theorem}

Numerical triviality means that $K_Z$ has zero intersection with every curve; a trivial canonical bundle is sufficient.
The general identity has an additional stringy analogue of the Chern number $c_1c_{m-1}$, which vanishes under this hypothesis.
For $m=4$, \eqref{eq:background-lw}, together with Hodge symmetry, duality and the vanishing of $h^{p,0}_{st}$ for $0<p<4$, gives the relation used in Corollary~\ref{cor:n4}.
For comparison, a smooth K3 surface has only the four corner entries $h^{0,0}=h^{2,0}=h^{0,2}=h^{2,2}=1$ and $h^{1,1}=20$; the left-hand side is $4$, equal to $\frac2{12}\cdot24$.

\subsection{The Crowley--Nordstr\"om invariant}
\label{sec:nu}
A $\rG_2$-structure on an oriented seven-manifold is a reduction of its frame bundle to the compact exceptional group $\rG_2$, equivalently a positive $3$-form $\varphi$; it determines a Riemannian metric and a spin structure.
Here \emph{positive} means that the form is pointwise equivalent to the standard $3$-form on $\R^7$ whose stabiliser is $\rG_2$.
A homotopy is a continuous path through such forms, with the induced spin structure fixed.
The Crowley--Nordstr\"om invariant $\nu$ is constant along these paths and is preserved by spin diffeomorphisms \cite[Thm.~1.3]{CrowleyNordstrom2015}.
Thus different values distinguish homotopy classes of structures, even on the same underlying manifold; we use only this invariant, not the further invariants of \cite{CrowleyNordstrom2015}.

For $7$-dimensional Calabi--Yau links, that is, for $N=5$ and $d=\sum_iw_i$, the Crowley--Nordstr\"om invariant $\nu\in\Z/48$ of the $\rG_2$-structure $\varphi\coloneq \eta'\wedge d\eta'+\mathrm{Im}\,\varepsilon$ induced by the contact Calabi--Yau structure $(\eta',\varepsilon)$ above \cite[Prop.~11]{CalvoAndradeRodriguezSaEarp2020} was also computed in \cite{AHHHSS2024}.
For a closed $7$-manifold $M$ with a $\rG_2$-structure $\varphi$ it is defined by \cite{CrowleyNordstrom2015}
\begin{equation}
\label{eq:nu-definition}
\nu(\varphi)\coloneq \chi(W)-3\sigma(W)\bmod 48.
\end{equation}
Here $W$ is any compact $8$-manifold with boundary $M$, oriented compatibly with the boundary orientation of $M$, that carries a $\Spin(7)$-structure inducing $\varphi$; this is a reduction of the structure group of $TW$, and no torsion-free metric is required.
The numbers $\chi(W)$ and $\sigma(W)$ are the Euler characteristic and the signature of $W$.
Crowley and Nordstr\"om prove that such a coboundary exists and that the residue in \eqref{eq:nu-definition} is independent of its choice \cite[Def.~1.2, Thm.~1.3, Lem.~3.4]{CrowleyNordstrom2015}.
For a link, the Milnor-fibre construction of \cite[\S3]{CalvoAndradeRodriguezSaEarp2020} supplies such a $\Spin(7)$ coboundary, with the complex orientation on $W$ and the induced boundary orientation on $K_f$.
The identification of its boundary $\rG_2$-structure with the contact Calabi--Yau structure, up to homotopy, is part of that geometric construction; it is not a consequence of the Hilbert-series calculation.
It has the homotopy type of a bouquet of $\mu(\bw)$ four-spheres \cite{Milnor1968,MilnorOrlik1970}, so $\chi(W)=1+\mu(\bw)$.
If $\{z^\alpha\}$ is a monomial basis of $\MM_f$ and $l(\alpha)\coloneq \sum_i(\alpha_i+1)w_i/d$, Steenbrink's signature formula \cite{Steenbrink1977} gives $\sigma(W)=\mu_+-\mu_-$, and hence \cite[Thm.~2]{CalvoAndradeRodriguezSaEarp2020}
\begin{equation}
\label{eq:nuSt}
\nu=1+\mu-3(\mu_+-\mu_-)\bmod48,\qquad \mu_\pm\coloneq \#\big\{\alpha:\ l(\alpha)\notin\Z,\ \lfloor l(\alpha)\rfloor\ \text{even (}+\text{) or odd (}-\text{)}\big\}.
\end{equation}
The intersection form on $H^4(W;\R)$ can be degenerate because $W$ has boundary.
In Steenbrink's formula the monomials with integral $l(\alpha)$ count its nullspace, so they contribute to $\mu$ but not to $\mu_+-\mu_-$ \cite[Thm.~2]{Steenbrink1977}.
The quantities $l(\alpha)-1$ are the spectral numbers in the convention with spectrum in $(-1,N-1)$; \eqref{eq:nuSt} uses the shifted quantities $l(\alpha)$ throughout.

\paragraph{Equivalence of singularities.}
\label{sec:conjectures}
Definition~\ref{def:weak-r-equivalence} records the numerical equivalence used in Conjecture~\ref{conj:wd}.
It is useful to compare it with equivalence under a change of coordinates.
Two admissible polynomials $f$ and $g$ are \emph{R-equivalent} if $g=f\circ\phi$ for a germ of biholomorphism $\phi:(\C^N,0)\to(\C^N,0)$; then $\phi^*$ maps $\cJ_f$ to $\cJ_g$, so their local analytic Milnor algebras, and hence their finite-dimensional polynomial Milnor algebras, are isomorphic as algebras.
The coordinate change need not preserve the grading.
Weak R-equivalence compares only the dimensions in the degrees of \eqref{eq:itoh}.
It does not require an isomorphism of Milnor algebras or a coordinate change; Example~\ref{ex:dwork} illustrates this distinction.

\section{Proof of the Weight Determinacy Conjecture}
\label{sec:conj7}

We first compute the Hilbert series of the Milnor algebra $\MM_f\coloneq S/\cJ_f$ for an admissible polynomial $f\in U_{\bw,d}$.
Together with Theorem~\ref{thm:itoh}, this proves weight determinacy and gives explicit formulae in the denumerants $p_{\bw}$ of \eqref{eq:denumerant}.
We then derive the same Hodge numbers directly from sheaf cohomology on $\cX_f$, before applying the formulae to the six outstanding threefold weight systems.

\subsection{The Hilbert series of the Milnor algebra}
\label{sec:hilbert}
The proof rests on two elementary facts: an admissible polynomial has no critical points outside the origin, and its partial derivatives form a homogeneous system of parameters in the Cohen--Macaulay ring $S$.
The terminology and algebraic inputs are recalled in Subsection~\ref{sec:algebra-background}.
The first lemma supplies the isolatedness and positivity needed for that argument.

\begin{lemma}
\label{lem:VJ}
Let $f\in\C[z_1,\dots,z_N]$ be weighted homogeneous of degree $d$ for positive weights $\bw$, with an isolated critical point at the origin; equivalently, $f\in U_{\bw,d}$ in the notation of Section~\ref{sec:weights}.
Then the common zero set of $\partial_1f,\dots,\partial_Nf$ in $\C^N$ is $\{0\}$, and $d>w_i$ for every $i$.
\end{lemma}
\begin{proof}
Since $\partial_if(t\cdot z)=t^{d-w_i}\partial_if(z)$, the critical set of $f$ is invariant under the $\C^*$-action.
All weights are positive, so the closure of every orbit contains the origin.
A critical point $z\neq0$ would therefore produce critical points $t\cdot z$ accumulating at $0$ as $t\to0$, which contradicts isolatedness.
If $d\le w_i$ for some $i$, the only monomial of degree at most $d$ involving $z_i$ is $z_i$ itself, and since $f$ has a critical point at $0$ it does not contain $z_i$; but then the whole $z_i$-axis is critical, again a contradiction.
\end{proof}

\begin{theorem}
\label{thm:hilbert}
Let $S\coloneq \C[z_1,\dots,z_N]$ have grading $\deg z_i=w_i$, and let $f\in S_d$ have an isolated critical point at the origin.
Write $\MM_f\coloneq S/(\partial_1f,\dots,\partial_Nf)$, and let $P_{\bw}$, $a_\ell$ and $\mu(\bw)$ be as in \eqref{eq:Pw}.
Then $\partial_1f,\dots,\partial_Nf$ form a regular sequence in $S$, and
\[\sum_{\ell\ge0}\dim(\MM_f)_\ell\,t^\ell=P_{\bw}(t)=\prod_{i=1}^N\frac{1-t^{\,d-w_i}}{1-t^{\,w_i}}\,.\]
In particular $P_{\bw}$ is a polynomial, $\dim(\MM_f)_\ell=a_\ell$ depends only on $(\bw,d)$ for every $\ell$, and $\dim\MM_f=\mu(\bw)$.
\end{theorem}
\begin{proof}
By Lemma~\ref{lem:VJ} and the Nullstellensatz, the radical of $\cJ_f$ is the maximal ideal $\fm\coloneq (z_1,\dots,z_N)$.
Here the \emph{radical} consists of the polynomials some power of which lies in $\cJ_f$; the Nullstellensatz identifies it with the ideal of the common zero set \cite[Ch.~4]{CoxLittleOShea2015}.
In particular, for each $i$ there is an integer $b_i>0$ with $z_i^{b_i}\in\cJ_f$.
Monomials with exponents less than the corresponding $b_i$ therefore span $S/\cJ_f$.
This proves the finite-dimensionality needed below; it does not assert that these spanning monomials are linearly independent.
Hence the $N$ homogeneous generators $\partial_if$, of positive degrees $d-w_i$, form a homogeneous system of parameters of the graded ring $S$.
Since $S$ is Cohen--Macaulay, they form a regular sequence by Proposition~\ref{prop:algebra-inputs}(ii), using the Cohen--Macaulay property recalled in Remark~\ref{rem:cohen-macaulay}.
The \emph{Koszul complex} packages multiplication by the generators and the relations among them.
Its first differential sends $(g_i)$ to $\sum_i g_i\partial_i f$, and its second sends $e_i\wedge e_j$ to $\partial_i f\,e_j-\partial_j f\,e_i$.
By Proposition~\ref{prop:algebra-inputs}(ii) it is exact in positive degrees, so it is a graded free resolution
\[0\to S\Big(-\sum_i(d-w_i)\Big)\to\cdots\to\bigoplus_{i<j}S(-(d-w_i)-(d-w_j))\to\bigoplus_iS(-(d-w_i))\to S\to\MM_f\to0 .\]
Here $S(-a)_k\coloneq S_{k-a}$, so the degree shift contributes the factor $t^a$ to a Hilbert series.
Taking the alternating sum of Hilbert series, and using $\sum_k\dim S_k\,t^k=\prod_i(1-t^{w_i})^{-1}$, gives the stated product.
Its value at $t=1$ is $\prod_i(d-w_i)/w_i=\mu(\bw)$.
The individual factors can also be checked without writing the whole resolution.
Put $R_0\coloneq S$ and $R_j\coloneq S/(\partial_1f,\dots,\partial_jf)$.
Regularity gives the short exact sequence
\begin{equation}
\label{eq:successive-hilbert-quotients}
0\longrightarrow R_{j-1}(-(d-w_j))
\xrightarrow{\ \partial_jf\ }R_{j-1}
\longrightarrow R_j\longrightarrow0.
\end{equation}
Taking dimensions in each degree multiplies the Hilbert series by $1-t^{d-w_j}$.
Iterating from $R_0=S$ to $R_N=\MM_f$ gives the same product; injectivity of these multiplication maps is precisely where regularity is used.
\end{proof}

The identity of Theorem~\ref{thm:hilbert} is classical, and is implicit in the computation of the Milnor number $\mu(\bw)$ by Milnor and Orlik \cite{MilnorOrlik1970}; see also \cite{Stanley1978,Soler2026} for the algebraic argument and \cite[eq.~(2)]{KreuzerSkarke1992} for its use in the classification of Calabi--Yau weight systems.
Theorem~\ref{cor:conj7} follows at once, proving Conjecture~\ref{conj:wd} in every graded degree and answering Itoh's weight-dependence question \cite[p.~123]{Itoh2004}.

\begin{proof}[Proof of Theorem~\ref{cor:conj7}]
A permutation of the coordinates identifies $\MM_g$ with the Milnor algebra of a polynomial with the same weights as $f$, and Theorem~\ref{thm:hilbert} applies to both.
The Sasakian Hodge numbers are dimensions of graded pieces by Theorem~\ref{thm:itoh}.
For \eqref{eq:nuSt}, a monomial $z^\alpha$ of degree $k$ has $l(\alpha)=(k+\sum_iw_i)/d$, which depends only on $k$, so $\mu_+$ and $\mu_-$ are determined by the Hilbert function.
\end{proof}

The isolated-singularity hypothesis is part of the setting of \cite{AHHHSS2024}, but it is not contained in the original displayed statement of the Weight Determinacy Conjecture, and it cannot be omitted.

\begin{example}[the hypothesis is necessary]
\label{ex:nonisolated}
Take $\bw=(1,1,1,1,1)$ and $d=5$.
The Fermat quintic $f\coloneq \sum_{i=1}^5z_i^5$ is admissible, $\cJ_f=(z_1^4,\dots,z_5^4)$ and $\dim(\MM_f)_{10}=101$.
The polynomial $g\coloneq \sum_{i=1}^4z_i^5$ has the same weights and degree but is singular along the $z_5$-axis; here $\cJ_g=(z_1^4,\dots,z_4^4)$, $\dim(\MM_g)_{10}=251$, and $\MM_g$ is infinite-dimensional.
\end{example}

The conclusion of Theorem~\ref{cor:conj7} concerns dimensions only, and this is as much as can be expected: the Milnor algebras themselves are in general not isomorphic, so weak R-equivalence is strictly weaker than R-equivalence \cite{Ahmed2012,MatherYau1982}.

\begin{example}[the Milnor algebras differ]
\label{ex:dwork}
In five variables, let $F\coloneq \sum_iz_i^5$ be the Fermat quintic and let $f_\psi\coloneq \sum_iz_i^5-5\psi z_1\cdots z_5$ be a Dwork quintic, with $\psi\neq0$ and $\psi^5\ne1$.
By Theorem~\ref{thm:hilbert} the two Milnor algebras have the same Hilbert function.
In $\MM_F$ the linear form $z_1$ satisfies $z_1^4=0$.
Suppose that some $x\in\MM_{f_\psi}$ of degree one satisfied $x^4=0$, and write $x=\sum_ia_iz_i\neq0$.
Then $x^4$ would lie in
\[
(\cJ_{f_\psi})_4=\mathrm{span}\{z_i^4-\psi\prod_{j\ne i}z_j\}.
\]
The coefficient of $z_i^3z_j$ in $x^4$ is $4a_i^3a_j$, so at most one $a_i$ is nonzero; but then $x^4=a_i^4z_i^4$, which is not in that span since $\psi\neq0$.
Hence $\MM_F\not\cong\MM_{f_\psi}$ as graded algebras.
Both algebras are generated in degree one, so an isomorphism of local algebras would induce an isomorphism of their associated graded algebras with respect to the maximal ideal, which are $\MM_F$ and $\MM_{f_\psi}$ themselves.
Hence the Milnor algebras are not isomorphic even as ungraded algebras, and $F$ and $f_\psi$ are not R-equivalent.
\end{example}

\begin{remark}[the Crowley--Nordstr\"om invariant]
\label{rem:nu}
The statement of Theorem~\ref{cor:conj7} about $\nu$ concerns only $N=5$ and $d=\sum_iw_i$, where the link carries the $\rG_2$-structure of Section~\ref{sec:nu}.
For every admissible $f$ the invariant $\nu$ of that structure is given by \eqref{eq:nuSt}, whose ingredients $\chi(W)=1+\mu(\bw)$ and $\sigma(W)=\mu_+-\mu_-$ are determined by the Hilbert function of $\MM_f$; by Theorem~\ref{cor:conj7}, $\nu$ therefore depends only on $\bw$.
The rest of Theorem~\ref{cor:conj7} holds for every non-degenerate $(\bw,d)$, including links that carry no such $\rG_2$-structure.
\end{remark}

\subsection{Closed formulae}
Theorem~\ref{thm:hilbert} reduces the Sasakian Hodge numbers to coefficients of $P_{\bw}$, which are given by the short alternating sums \eqref{eq:incl-excl} of denumerants.
Proposition~\ref{prop:ci-duality} supplies the Gorenstein pairing and explains the symmetry used to shorten these formulae.
For the Fermat quintic, the socle is spanned by $z_1^3\cdots z_5^3$, of degree $15$.
The monomial $z^\alpha$ pairs with $z_1^{3-\alpha_1}\cdots z_5^{3-\alpha_5}$ to give this generator.
Taking the coefficient of that generator identifies the complementary monomial bases as dual bases, so the symmetry $a_k=a_{15-k}$ can be checked directly.

\begin{corollary}
\label{cor:closed}
Let $f$ be a weighted homogeneous polynomial of degree $d$ with an isolated critical point at the origin, in $N=n+2$ variables of weights $\bw$.
Let $a_k\coloneq \dim(\MM_f)_k$ be the Hilbert coefficients of Theorem~\ref{thm:hilbert}, and let $p_{\bw}$ denote the monomial count \eqref{eq:denumerant}.
Then $\MM_f$ is a graded Artinian Gorenstein algebra with socle in degree $s\coloneq \sum_i(d-2w_i)$, so $a_k=a_{s-k}$.
In the Calabi--Yau case $s=nd$, and $h^{p,q}_S(K_f)=a_{pd}=a_{qd}$ for $p+q=n$.
Moreover, with $m\coloneq \min(p,q)$, the coefficient $a_{md}$ is the sum of the terms of \eqref{eq:incl-excl} with $|I|\le m$, all other terms being zero.
Explicitly,
\begin{align*}
n=2:&\quad h^{2,0}_S=1,\quad h^{1,1}_S=p_{\bw}(d)-\textstyle\sum_ip_{\bw}(w_i);\\
n=3:&\quad h^{3,0}_S=1,\quad h^{2,1}_S=p_{\bw}(d)-\textstyle\sum_ip_{\bw}(w_i);\\
n=4:&\quad h^{4,0}_S=1,\quad h^{3,1}_S=p_{\bw}(d)-\textstyle\sum_ip_{\bw}(w_i),\quad
h^{2,2}_S=p_{\bw}(2d)-\textstyle\sum_ip_{\bw}(d+w_i)+\sum_{i<j}p_{\bw}(w_i+w_j).
\end{align*}
\end{corollary}
\begin{proof}
By Theorem~\ref{thm:hilbert}, $\MM_f$ is a graded Artinian complete intersection.
Proposition~\ref{prop:ci-duality} makes it Gorenstein with the stated socle degree and symmetric Hilbert function; directly, $t^sP_{\bw}(1/t)=P_{\bw}(t)$.
In the Calabi--Yau case, for $\ell=md$ the argument of $p_{\bw}$ in \eqref{eq:incl-excl} is $$(m-|I|)d+\sum_{i\in I}w_i.$$
If $|I|=m+1\le N-1$ this is negative because $\sum_{i\in I}w_i<d$, and if $|I|\ge m+2$ it is at most $-d$.
Taking $m\coloneq \min(p,q)$, which is allowed since $a_{pd}=a_{qd}$, gives the stated truncation, and Theorem~\ref{thm:itoh} identifies $a_{pd}$ with $h^{p,q}_S$.
\end{proof}

The truncation holds in the degree $md$ with $m\coloneq \min(p,q)$, but not in the symmetric degree.
For the quintic,
\[
h^{2,1}_S=a_5=p_{\bw}(5)-5p_{\bw}(1)=126-25=101,
\]
whereas in degree $10$ the terms with $|I|=2$ do not vanish:
\[
\begin{aligned}
a_{10}&=p_{\bw}(10)-5p_{\bw}(6)+10p_{\bw}(2)\\
&=1001-1050+150=101.
\end{aligned}
\]

\begin{remark}[polynomial deformations]
\label{rem:moduli}
In the Calabi--Yau case the number $a_d=p_{\bw}(d)-\sum_ip_{\bw}(w_i)$ has a deformation-theoretic meaning.
The group $G$ of weighted coordinate changes $z_i\mapsto g_i(z)$, with $g_i\in S_{w_i}$, has Lie algebra $\bigoplus_iS_{w_i}$, and the differential of its orbit map at $f$ is
\[\bigoplus_iS_{w_i}\longrightarrow S_d,\qquad(g_i)\longmapsto\sum_ig_i\,\partial_if,\]
with image $(\cJ_f)_d$.
This map is injective by the regular-sequence statement of Theorem~\ref{thm:hilbert}: a relation $\sum_ig_i\partial_if=0$, called a \emph{syzygy}, is a relation among the generators of the regular sequence $\partial_1f,\dots,\partial_Nf$, hence a combination of Koszul syzygies, whose coefficients would have degree $w_i+w_j-d<0$.
By Euler's identity $f\in(\cJ_f)_d$, so the rescaling of the equation is already contained in the image.
Hence $a_d=\dim S_d/(\cJ_f)_d$ is the dimension of the space of first-order deformations of $f$ modulo weighted coordinate changes and rescaling, that is, the number of polynomial deformations of $X_f$; equivalently, since the stabiliser of $f$ in $G$ is finite, the $G$-orbit of $f$ has codimension $a_d$ in $S_d$.
For the quintic it is $126-25=101$, and for the example $\bw=(22,29,49,50,75)$ of \cite{AHHHSS2024} it is $7-5=2$.
For $n\ge3$, this also identifies the first-order complex deformations of the smooth quotient orbifold $\cX_f$.
Indeed, the dual Euler and normal sequences of Lemma~\ref{lem:forms}, together with the vanishings of Lemma~\ref{lem:lines}, give
\[
H^1(\cX_f,T_{\cX_f})\cong
\frac{S_d/\C f}{\mathrm{im}\big(\bigoplus_i S_{w_i}\longrightarrow S_d/\C f\big)}
\cong(\MM_f)_d.
\]
The required vanishing is $H^1(T_{\WP}|_{\cX_f})=0$, obtained from $H^1(\cO(w_i))=H^2(\cO)=0$ in the Euler sequence; here $T_{\WP}$ denotes the tangent bundle of the weighted projective orbifold.
Contraction with a holomorphic volume form identifies this space with $H^{n-1,1}(\cX_f)$.
This concerns deformations of the orbifold; the additional stringy contribution is studied in Section~\ref{sec:conj8}.
For K3 surfaces $a_d$ counts deformations of $X_f$ as a weighted hypersurface, rather than the full deformation space of its smooth K3 resolution, whose complex dimension is $20$.
The Gr\"obner-basis approach of \cite{AHHHSS2024} obtains the same numbers after enumerating a monomial basis of the whole Milnor algebra, of size $\mu(\bw)$, whereas Corollary~\ref{cor:closed} needs $\sum_{j=0}^{m}\binom Nj$ denumerants for $h^{p,q}_S$ with $m\coloneq \min(p,q)$.
This is $N+1$ for $h^{n-1,1}_S$ in every dimension, and at most $1+N+\binom N2$, which is $22$ for $N=6$, for every Hodge number in dimensions $n\le5$.
\end{remark}

\subsection{A sheaf-cohomological proof}
\label{sec:sheaf}
We prove Theorem~\ref{thm:sheaf} by computing the Hodge numbers directly on the quotient orbifold.
The five lemmas below organise the argument into three stages.
First, Lemma~\ref{lem:basic} identifies basic cohomology with coherent cohomology on the leaf space.
Next, the Gysin sequence of Lemma~\ref{lem:gysin} and the connectivity calculation of Lemma~\ref{lem:offmiddle} isolate the primitive contribution in the middle degree.
Finally, Lemmas~\ref{lem:lines} and~\ref{lem:forms} compute the Euler characteristics from which the remaining Hodge numbers are recovered.
Combining the resulting formula with Theorem~\ref{thm:hilbert} recovers Itoh's formula, as explained in Section~\ref{sec:intro}.

The proof uses the notation of Section~\ref{sec:links}, with $\cX\coloneq \cX_f$.
Let $A\coloneq S/(f)$, a graded ring with Hilbert function $h_A(k)\coloneq \dim A_k=p_{\bw}(k)-p_{\bw}(k-d)$, and put $\kappa\coloneq d-\sum_iw_i$.
The ring $A$ consists of polynomial functions on the affine cone, whereas $\MM_f\coloneq S/\cJ_f$ records its critical-point equations.
A graded $A$-module carries compatible weighted degrees; the shifted module $A(k)$, with $A(k)_j\coloneq A_{j+k}$, gives the line bundle $\cO_\cX(k)$ after restriction to $Y$ and descent to the orbifold.
Thus the grading is the algebraic counterpart of the line-bundle twist \cite[\S\S1.1,1.4]{Dolgachev1982}.
Lemmas~\ref{lem:lines} and~\ref{lem:forms} below do not involve the link: they hold for every $f\in U_{\bw,d}$ in $N=r+2\ge3$ variables, with $r$ in place of $n$, and they are used in this generality in Section~\ref{sec:conj8}.

\paragraph{Descent to the leaf space.}
We begin by identifying the basic Dolbeault complex with the Dolbeault complex of the quotient orbifold.
This makes it possible to compute the basic Hodge numbers by coherent sheaf cohomology.

\begin{lemma}[basic forms are forms on the leaf space]
\label{lem:basic}
Let $K_f$ be the quasi-regular Sasakian link of an admissible weighted homogeneous polynomial $f\in U_{\bw,d}$, and let $\cX\coloneq \cX_f$ be its quotient orbifold as in Section~\ref{sec:links}.
For all $p,q$ there are isomorphisms $H^{p,q}_B(K_f)\cong H^q(\cX,\Omega^p_\cX)$; in particular $h^{p,q}_B(K_f)=h^{p,q}(\cX)$.
\end{lemma}
\begin{proof}
The weighted radial flow identifies $Y$ circle-equivariantly with $K_f\times\R_{>0}$ as a smooth manifold, since the real weighted dilations $z\mapsto e^r\cdot z$ move $Y$ transversally to the sphere.
Basic forms on $K_f$ therefore correspond to forms on $Y$ that are invariant under $\C^*$ and annihilated by the real and imaginary parts of the Euler vector field $E_{\bw}$, which are exactly the differential forms on $\cX$.
The transverse complex structure of the Reeb foliation is the one induced from $Y$ \cite[Ch.~9]{BoyerGalicki2008}, so this correspondence respects types and $\bar\partial$.
In a uniformising chart the $\bar\partial$-Poincar\'e lemma holds for invariant forms, by averaging over the finite stabiliser, so smooth basic $(p,\bullet)$-forms form a fine resolution of $\Omega^p_\cX$; compare \cite[\S1.2.2]{GNT2016}.
\end{proof}

\paragraph{The Gysin sequence and the middle degree.}
To extract the primitive contribution, we next compare basic cohomology with the ordinary cohomology of the link.
The connecting map is the Lefschetz operator introduced in Section~\ref{sec:links}.

\begin{lemma}[Gysin sequence]
\label{lem:gysin}
Let $K_f$ be the Sasakian link of $f\in U_{\bw,d}$, with contact form $\eta$, and put $L\coloneq [d\eta]\wedge\,\cdot\,$ on its real basic cohomology.
There is a long exact sequence
\[\cdots\to H^{k-2}_B(K_f)\xrightarrow{\ L\ }H^k_B(K_f)\to H^k(K_f;\R)\to H^{k-1}_B(K_f)\xrightarrow{\ L\ }H^{k+1}_B(K_f)\to\cdots\]
\end{lemma}
\begin{proof}
Averaging over the circle shows that circle-invariant forms compute $H^*(K_f;\R)$.
Every invariant form decomposes uniquely as $\gamma=\alpha+\eta\wedge\beta$ with $\alpha,\beta$ basic, namely $\beta=\iota_\xi\gamma$, and
\[
d\gamma=d\alpha+d\eta\wedge\beta-\eta\wedge d\beta.
\]
Hence $\gamma\mapsto\iota_\xi\gamma$ is, up to sign, a surjective chain map onto the basic complex shifted by one, with kernel the basic complex.
The connecting homomorphism of the resulting short exact sequence of complexes sends $[\beta]$ to $[d(\eta\wedge\beta)]=[d\eta\wedge\beta]$.
\end{proof}

Milnor's connectivity theorem now determines all basic cohomology groups away from the middle degree.
This is the reason that the Euler characteristics computed below suffice to recover the remaining Hodge numbers.

\begin{lemma}[basic cohomology away from the middle]
\label{lem:offmiddle}
Let $K_f$ be the $(2n+1)$-dimensional Sasakian link of an admissible polynomial $f\in U_{\bw,d}$, where $n\ge2$, with contact form $\eta$ and Lefschetz operator $L$ as in Lemma~\ref{lem:gysin}.
For $0\le k\le2n$ with $k\ne n$, the group $H^k_B(K_f)$ is spanned by $[d\eta]^{k/2}$ if $k$ is even and vanishes if $k$ is odd.
Consequently $h^{p,q}_B=\delta_{pq}$ whenever $p+q\neq n$, the map $L:H^{k-2}_B\to H^k_B$ is injective for $k\le n+1$, and $H^n(K_f;\R)\cong H^n_B/L\,H^{n-2}_B$.
In particular,
\begin{equation}
\label{eq:middle-betti}
b_n(K_f)=\sum_{p+q=n}h^{p,q}_S(K_f).
\end{equation}
\end{lemma}
\begin{proof}
Since $K_f$ is $(n-1)$-connected \cite{Milnor1968}, Poincar\'e duality gives $H^k(K_f;\R)=0$ unless $k\in\{0,n,n+1,2n+1\}$.
Inserting this into Lemma~\ref{lem:gysin}: for $1\le k\le n-1$ the map $L:H^{k-2}_B\to H^k_B$ is surjective and $L:H^{k-1}_B\to H^{k+1}_B$ is injective, so
\[
H^1_B=0,\qquad H^k_B=L^{\lfloor k/2\rfloor}H^{k-2\lfloor k/2\rfloor}_B\quad(k\le n-1),
\]
while $L$ is injective on $H^j_B$ for $j\le n-2$.
Basic forms of degree above $2n$ vanish, so the end of the sequence gives $H^{2n}_B\cong H^{2n+1}(K_f;\R)\cong\R$.
For $n+2\le k\le2n$ the same argument shows that $L:H^j_B\to H^{j+2}_B$ is bijective for $n+1\le j\le2n-2$ and that $H^{2n-1}_B$ injects into $H^{2n+1}_B=0$.
This determines the dimensions as claimed.
The class $[d\eta]^n$ is nonzero: if $(d\eta)^n=d\beta$ with $\beta$ basic, then Stokes' theorem gives
\[
\int_{K_f}\eta\wedge(d\eta)^n=\int_{K_f}d\eta\wedge\beta=0,
\]
since $d\eta\wedge\beta$ is a basic form of degree $2n+1$; but $\eta\wedge(d\eta)^n$ is a volume form.
Hence every $[d\eta]^j$ is nonzero, and it spans $H^{2j}_B$ when $2j\neq n$.
The basic Hodge decomposition \cite[Thm.~3.4.6]{ElKacimi1990} is
\[
H^k_B\otimes\C=\bigoplus_{p+q=k}H^{p,q}_B.
\]
Since $d\eta$ has type $(1,1)$, this gives $h^{p,q}_B=\delta_{pq}$ for $p+q\ne n$.
Finally, $L:H^{n-1}_B\to H^{n+1}_B$ is injective: for $n$ even its source vanishes, and for $n$ odd its source is spanned by $[d\eta]^{(n-1)/2}$, whose image $[d\eta]^{(n+1)/2}$ is nonzero.
The exact sequence of Lemma~\ref{lem:gysin} at $k=n$ therefore gives $H^n(K_f;\R)\cong H^n_B/L\,H^{n-2}_B$.
Here $L\,H^{n-2}_B$ is spanned by $[d\eta]^{n/2}$ for $n$ even and vanishes for $n$ odd, so $b_n(K_f)=\sum_{p+q=n}(h^{p,q}_B-\delta_{pq})$, which is $\sum_{p+q=n}h^{p,q}_S$ by \eqref{eq:hS-basic}.
\end{proof}

\begin{example}[the primitive correction for a quartic K3 surface]
\label{ex:quartic-gysin}
Take $f=z_1^4+z_2^4+z_3^4+z_4^4$, with all weights equal to $1$.
The circle action on the link is free, and its quotient is a smooth quartic K3 surface $X\subset\WP^3$.
Thus Lemma~\ref{lem:basic} is the usual descent of horizontal invariant forms along a principal circle bundle.
The relevant part of Lemma~\ref{lem:gysin} reads
\[
H^0(X;\R)\xrightarrow{\ L\ }H^2(X;\R)
\longrightarrow H^2(K_f;\R)\longrightarrow H^1(X;\R)=0.
\]
The image of $L$ is the one-dimensional span of the hyperplane class.
Consequently $b_2(K_f)=22-1=21$ and $h^{1,1}_S=20-1=19$; the $(2,0)$ and $(0,2)$ classes survive unchanged.
This explains concretely the diagonal correction in Lemma~\ref{lem:offmiddle} and Proposition~\ref{prop:k3}.
\end{example}

The remaining calculation takes place on the quotient orbifold.
For a coherent sheaf $\cF$, its holomorphic Euler characteristic is
\[
\chi(\cX,\cF)\coloneq \sum_i(-1)^i\dim H^i(\cX,\cF).
\]
Lemma~\ref{lem:offmiddle} leaves only one unknown Hodge number in each row, so this alternating sum determines it.
For fixed $p$, the only potentially unknown term of $\chi(\cX,\Omega^p_\cX)$ has $q=n-p$; every other term is the known diagonal contribution $\delta_{pq}$.
Thus taking an Euler characteristic here does not lose the individual Hodge number we seek.
We first compute it for line bundles, then express differential forms in terms of those bundles.
The next two lemmas also cover curves, since the twisted loci in Section~\ref{sec:conj8} can have complex dimension one.

\paragraph{Line-bundle cohomology.}
In the proof, $H^j_\fm(A)$ denotes \emph{local cohomology supported at the cone vertex} $\fm\coloneq (z_1,\dots,z_N)$.
It relates cohomology on the affine cone to that on the punctured cone $Y$, where the orbifold quotient is formed.
It is computed by the augmented \v{C}ech complex
\[
0\longrightarrow A\longrightarrow\bigoplus_i A_{z_i}
\longrightarrow\bigoplus_{i<j} A_{z_iz_j}\longrightarrow\cdots,
\]
where $A_{z_i}$ allows powers of $z_i$ in the denominator.
The term $A$ is in cohomological degree zero.
Deleting it gives the \v{C}ech complex of the cover of $Y$ by $z_i\ne0$, with its degrees shifted by one.
This explains the passage from $H^i(Y,\cO_Y)$ to $H^{i+1}_\fm(A)$ in the proof \cite[Ch.~3]{BrunsHerzog1993}.
More explicitly, the beginning of the comparison sequence is
\begin{equation}
\label{eq:punctured-cone-comparison}
0\longrightarrow H^0_\fm(A)\longrightarrow A
\longrightarrow H^0(Y,\cO_Y)
\longrightarrow H^1_\fm(A)\longrightarrow0.
\end{equation}
Thus unique extension of functions across the vertex follows from the vanishing of both local cohomology groups displayed here.
The input is Proposition~\ref{prop:algebra-inputs}(iii), the Cohen--Macaulay case of depth-vanishing and dimension-vanishing for local cohomology.
For a Cohen--Macaulay ring of dimension $r+1$, depth equals dimension (Remark~\ref{rem:cohen-macaulay}), so only degree $r+1$ can survive.
This is an algebraic extension argument on the possibly singular cone; it does not require the vertex to be smooth.

\begin{lemma}[line bundles on the leaf space]
\label{lem:lines}
Let $f$ be weighted homogeneous of degree $d$ with an isolated critical point at the origin, in $N=r+2\ge3$ variables of weights $\bw$, and let $\cX\coloneq [(f^{-1}(0)\setminus0)/\C^*]$.
Set $A\coloneq \C[z_1,\dots,z_N]/(f)$, $h_A(k)\coloneq \dim A_k$ and $\kappa\coloneq d-\sum_iw_i$.
For every $k\in\Z$,
\begin{equation}
\label{eq:line-cohomology}
\begin{aligned}
H^0(\cX,\cO_\cX(k))&\cong A_k,\\
H^i(\cX,\cO_\cX(k))&=0\qquad(0<i<r),\\
\dim H^r(\cX,\cO_\cX(k))&=h_A(\kappa-k).
\end{aligned}
\end{equation}
Hence
\[\chi\big(\cX,\cO_\cX(k)\big)=h_A(k)+(-1)^rh_A(\kappa-k).\]
\end{lemma}
\begin{proof}
We compute in the algebraic category, using the comparison in Section~\ref{sec:links}.
Algebraic $\C^*$-representations decompose into integer weight spaces, and equivariant maps preserve that decomposition.
Taking a fixed weight is therefore exact; this is the \emph{linear reductivity} used here.
Applied to the \v{C}ech complex, it identifies $H^i(\cX,\cO_\cX(k))$ with the weight-$k$ part of $H^i(Y,\cO_Y)$, where $Y$ is the punctured spectrum of $A$; compare \cite[\S1.4]{Dolgachev1982}.
Indeed, the differential preserves weight, so a weight-$k$ cocycle that is a boundary has a weight-$k$ primitive: take that homogeneous component of any primitive.
This explains why taking cohomology and selecting the line-bundle twist commute.
The ring $A$ is a hypersurface ring of dimension $r+1\ge2$, hence Cohen--Macaulay in the sense recalled in Remark~\ref{rem:cohen-macaulay}, so
\[
H^0(Y,\cO_Y)=A,\qquad
H^i(Y,\cO_Y)\cong H^{i+1}_\fm(A)\qforq i\ge1,
\]
where the latter group vanishes unless $i=r$ \cite[Ch.~3]{BrunsHerzog1993}.
Here $f$ is a non-zero-divisor in $S$; quotienting a Cohen--Macaulay ring by such an element lowers depth and dimension by one, preserving the equality in Remark~\ref{rem:cohen-macaulay} by Proposition~\ref{prop:algebra-inputs}(iii).
Since $r+1\ge2$, the two local cohomology groups in \eqref{eq:punctured-cone-comparison} vanish.
Taking weight $k$ gives $H^0(\cX,\cO_\cX(k))=A_k$; the same depth calculation gives the intermediate vanishings.
It remains to compute the one possible higher group, in degree $r$.
The sequence
\[
0\to S(-d)\xrightarrow{f}S\to A\to0,
\]
together with $H^j_\fm(S)=0$ for $j\neq N$ and $H^N_\fm(A)=0$, identifies $H^{r+1}_\fm(A)$ with the kernel of the surjection $f:H^N_\fm(S)(-d)\to H^N_\fm(S)$.
The module $H^N_\fm(S)$ has a basis of inverse monomials $z^{-\beta}$ with all $\beta_i\ge1$, so its degree-$k$ part has dimension $p_{\bw}(-k-\sum_iw_i)$.
Indeed, in the last term of the \v{C}ech complex all variables are inverted; the preceding term removes precisely the monomials with at least one non-negative exponent.
Hence
\[
\begin{aligned}
\dim H^{r+1}_\fm(A)_k
&=p_{\bw}(d-k-\sum_iw_i)-p_{\bw}(-k-\sum_iw_i)\\
&=h_A(\kappa-k).
\end{aligned}
\]
\end{proof}

\paragraph{Differential forms and adjunction.}
To pass from line bundles to differential forms, the conormal sequence imposes the relation $df=0$ on ambient one-forms, and the Euler sequence selects forms that annihilate the orbit direction $E_{\bw}$.
Together they express the classes of differential forms in terms of the line bundles covered by Lemma~\ref{lem:lines}.
Only this identity of classes, together with the canonical bundle, is needed in the proof of Theorem~\ref{thm:sheaf}.

The \emph{Grothendieck group} used here has a class $[V]$ for each vector bundle, with $[V]=[V']+[V'']$ for every short exact sequence $0\to V'\to V\to V''\to0$; multiplication is induced by tensor product.
The additivity of the holomorphic Euler characteristic follows from the long exact sequence in cohomology.
It allows us to use an identity of classes without choosing a splitting of either bundle sequence.
Taking determinants in the conormal sequence is the orbifold form of \emph{adjunction}, the rule expressing the canonical bundle of a hypersurface in terms of its ambient canonical bundle and normal line bundle; compare \cite[\S3.3]{Dolgachev1982}.

\begin{lemma}[differential forms on the leaf space]
\label{lem:forms}
Let $\cX$ be the $r$-dimensional quotient orbifold of the quasi-smooth hypersurface in Lemma~\ref{lem:lines}, with $\kappa\coloneq d-\sum_iw_i$, and write $\cE\coloneq \bigoplus_{i=1}^N\cO_\cX(-w_i)$.
In the Grothendieck group of vector bundles on $\cX$,
\[[\Omega^p_\cX]=\sum_{i+j+m=p}(-1)^{i+j}\,[\cO_\cX(-jd)]\cdot\big[{\textstyle\bigwedge^m}\cE\big],\qandq\Omega^r_\cX\cong\cO_\cX(\kappa).\]
\end{lemma}
\begin{proof}
The $\C^*$-equivariant vector bundles on $Y$ below descend to $\cX$.
The conormal sequence of $Y\subset\C^N\setminus0$,
\[0\to\cO_Y(-d)\xrightarrow{\ df\ }\bigoplus_{i=1}^N\cO_Y(-w_i)\,dz_i\to\Omega^1_Y\to0,\]
is exact on the left because $df\neq0$ on $Y$; here $dz_i$ has weight $w_i$ and $df$ has weight $d$.
It descends to an exact sequence
\begin{equation}
\label{eq:conormal-quotient}
0\to\cO_\cX(-d)\to\cE\to\cB\to0,
\end{equation}
where $\cB$ is the bundle on $\cX$ obtained from $\Omega^1_Y$.
The pullbacks of forms on $\cX$ are the forms on $Y$ annihilated by the nowhere-vanishing vector field $E_{\bw}$, so contraction with $E_{\bw}$ gives a second exact sequence
\begin{equation}
\label{eq:euler-quotient}
0\to\Omega^1_\cX\to\cB\to\cO_\cX\to0.
\end{equation}
Surjectivity is local: a nowhere-vanishing vector field can locally be completed to a frame, whose dual contains a one-form evaluating to $1$ on that field.
The kernel consists of horizontal forms, with equivariance already built into $\cB$.
Thus the two sequences remove, successively, the normal direction to the hypersurface and the tangent direction to the quotient orbit.
Their ranks decrease from $N$ to $N-1$ to $N-2=r$, as expected for the ambient space, the punctured cone and its leaf space.
In terms of $\lambda_u(V)\coloneq \sum_m[\bigwedge^mV]u^m$, which is multiplicative in short exact sequences, the sequences \eqref{eq:conormal-quotient} and~\eqref{eq:euler-quotient} give
\[\lambda_u(\Omega^1_\cX)=\frac{\lambda_u(\cE)}{(1+u)(1+[\cO_\cX(-d)]u)}=\frac{\prod_i\big(1+[\cO_\cX(-w_i)]u\big)}{(1+u)(1+[\cO_\cX(-d)]u)},\]
whose coefficient of $u^p$ is the stated identity.
Here $u$ is a formal variable, and each denominator is inverted as a power series with constant term $1$; the exterior-power filtration of an exact sequence gives the multiplicativity of $\lambda_u$.
Taking determinants gives
\[
\Omega^r_\cX\cong\det\cB\cong\det\cE\otimes\cO_\cX(d)=\cO_\cX(d-\sum_iw_i). \qedhere
\]
\end{proof}

For $p=1$, Lemma~\ref{lem:forms} reduces to the easily checked identity
\[
[\Omega^1_\cX]=\sum_i[\cO_\cX(-w_i)]-[\cO_\cX(-d)]-[\cO_\cX].
\]
The two subtractions encode the conormal line and the orbit direction, respectively.
This is an equality of classes for computing Euler characteristics, not a direct-sum decomposition of the cotangent bundle.

\begin{proof}[Proof of Theorem~\ref{thm:sheaf}]
Fix $0\le p\le n$.
We compute $\chi(\cX,\Omega^p_\cX)$ in two ways: from the known off-middle cohomology, and from the line-bundle formula.
In the first calculation there is an ambient diagonal class in bidegree $(p,p)$, and the only additional contribution lies in bidegree $(p,n-p)$.
When $2p=n$ these occupy the same group, whose dimension is $1+h^{p,p}_S$; this accounts for the diagonal correction below.
By Definition~\ref{def:sasakian} and Lemma~\ref{lem:offmiddle}, $$h^{0,0}_S=1, \quad h^{p,q}_S=0 \qforq 0<p+q<n, \qandq h^{p,n-p}_S=h^{p,n-p}_B-\delta_{2p,n}.
$$
Lemmas~\ref{lem:basic} and~\ref{lem:offmiddle} then give
\[\chi(\cX,\Omega^p_\cX)=\sum_q(-1)^qh^{p,q}_B=(-1)^p+(-1)^{n-p}\,h^{p,n-p}_S .\]
On the other hand, by Lemmas~\ref{lem:forms} and~\ref{lem:lines} with $r=n$,
\[\chi(\cX,\Omega^p_\cX)=\sum_{m=0}^p(-1)^{p-m}\sum_{|I|=m}\ \sum_{j=0}^{p-m}\chi\Big(\cO_\cX\big(-jd-\textstyle\sum_{i\in I}w_i\big)\Big).\]
The twist is negative except for $(I,j)=(\emptyset,0)$, so each term equals $(-1)^nh_A(\kappa+jd+w_I)$, with $w_I\coloneq \sum_{i\in I}w_i$, apart from the extra summand $1$ coming from $h_A(0)$.
Comparing the two expressions,
\[h^{p,n-p}_S=\sum_{m=0}^p(-1)^m\sum_{|I|=m}\ \sum_{j=0}^{p-m}h_A(\kappa+jd+w_I).\]
Since $h_A(k)=p_{\bw}(k)-p_{\bw}(k-d)$, the sum over $j$ telescopes to
\[
p_{\bw}(\kappa+(p-m)d+w_I)-p_{\bw}(w_I-\sum_iw_i),
\]
and the last term vanishes because $|I|\le p<N$.
Since
\[
\kappa+(p-|I|)d+w_I=(p+1)d-\sum_iw_i-\sum_{i\in I}(d-w_i),
\]
this is the first equality in \eqref{eq:sheaf}.
For the second, in the formal expansion \eqref{eq:incl-excl} with $\ell\coloneq (p+1)d-\sum_iw_i$ the terms with $|I|\ge p+1$ have argument
\[
(p+1-|I|)d-\sum_{i\notin I}w_i<0,
\]
so the two sums coincide; this step does not use Theorem~\ref{thm:hilbert}.
\end{proof}

\begin{remark}
The proof of Theorem~\ref{thm:sheaf} uses Milnor's connectivity theorem in Lemma~\ref{lem:offmiddle}, basic Hodge theory in Definition~\ref{def:sasakian}, and local cohomology of hypersurface rings in Lemma~\ref{lem:lines}.
It also produces the index convention of Theorem~\ref{thm:itoh} directly, and the case $p=n$ reproduces the familiar count $h^{n,0}=p_{\bw}(d-\sum_iw_i)$ of holomorphic volume forms.
Term-by-term evaluations of \eqref{eq:sheaf} agree with classical examples, Hilbert-series coefficients, independent Gr\"obner-basis computations and the datasets below.
The checks, including their scope in dimension one, are recorded in Appendix~\ref{app:comp} under \emph{Checks of the sheaf-cohomological formula}.
\end{remark}

\begin{remark}[the Sasakian Hodge diamond]
\label{rem:diamond}
Theorem~\ref{thm:sheaf} determines every Sasakian Hodge number of $K_f$, but all except those in the middle degree take trivial values.
Below the middle degree the link has no real cohomology, since it is $(n-1)$-connected \cite{Milnor1968}, and the only basic classes there are the powers $[d\eta]^j$ of the transverse K\"ahler class (Lemma~\ref{lem:offmiddle}).
For $j\ge1$ the class $[d\eta]^j$ is not primitive, because $L^{n-2j+1}[d\eta]^j=[d\eta]^{n-j+1}$ is nonzero, so these classes contribute nothing to $h^{p,q}_S$, while the constant functions give $h^{0,0}_S=1$.
Above the middle degree no new numbers arise: the Hodge star maps the harmonic $k$-forms on $K_f$ isomorphically onto the harmonic $(2n+1-k)$-forms, which is why Definition~\ref{def:sasakian} only concerns $p+q\le n$.
In the middle degree $h^{p,q}_S=h^{q,p}_S$ by complex conjugation, and in the Calabi--Yau case $h^{n,0}_S=p_{\bw}(0)=1$.
The Sasakian Hodge numbers in degrees $k\le n$ therefore consist of $h^{0,0}_S=1$, zeros in every degree $0<k<n$, and the middle row $(1,h^{n-1,1}_S,\dots,h^{1,n-1}_S,1)$; for a $9$-dimensional link, for instance, $h^{1,0}_S$, $h^{1,1}_S$, $h^{2,1}_S$ and every other $h^{p,q}_S$ with $0<p+q<4$ vanish.
In the Calabi--Yau case the only non-trivial numbers are thus $h^{1,1}_S$ for $n=2$, $h^{2,1}_S$ for $n=3$ and the pair $(h^{3,1}_S,h^{2,2}_S)$ for $n=4$, and these are the Sasakian Hodge numbers listed in Sections~\ref{sec:six}, \ref{sec:5d} and~\ref{sec:9d}.
They determine the real cohomology of the link:
\begin{equation}
\label{eq:link-betti}
b_0=b_{2n+1}=1,\qquad b_n=b_{n+1}=\sum_{p+q=n}h^{p,q}_S,
\end{equation}
and all other Betti numbers vanish.
\end{remark}

\begin{remark}[moduli of transverse Calabi--Yau structures]
\label{rem:moriyama}
For a compact manifold carrying a transverse Calabi--Yau structure with a fixed taut Reeb foliation, Moriyama showed that the moduli space $\mathfrak M_{CY}$ is smooth, with real tangent space \cite[Thm.~6.5]{Moriyama2011}
\[
H^{n,0}_B\oplus H^{n-1,1}_B\oplus P^{1,1}_{B,\R}.
\]
The quotient is by the identity component of the group of foliation-preserving diffeomorphisms; the basic K\"ahler class is not fixed.
The first two summands are regarded as real vector spaces and describe the transverse complex structure together with its holomorphic volume form, while the last describes primitive K\"ahler-class variations; for $n=3$ the dimension is $b_3+b_1+h^{1,1}_B-1$ \cite[eq.~(25)]{Moriyama2011}.
In \cite{AHHHSS2024} the term $h^{1,1}$ of this formula was read as a Sasakian Hodge number not accessible through the Milnor algebra, and mirror symmetry was suggested as a way to obtain it.
It is the basic Hodge number, however, and $h^{1,1}_B=1$ for every link with $n\ge3$ by Lemma~\ref{lem:offmiddle}, so that $P^{1,1}_B=0$ and
\[\dim\mathfrak M_{CY}(K_f)=2+2h^{n-1,1}_S\qquad(n\ge3),\]
which is $b_3(K_f)$ for $7$-dimensional links; for K3 links, \cite[eq.~(22)]{Moriyama2011} gives
\[
\dim\mathfrak M_{CY}(K_f)=3b_2-2=3h^{1,1}_S+4.
\]
These dimensions are therefore known for every link in the datasets considered here.
For $n\ge3$ there are no primitive K\"ahler-class directions, but the summand $H^{n,0}_B\cong\C$ still supplies two real directions.
Writing the transverse Calabi--Yau structure as $(\omega,\varepsilon)$, these are its common scale and the phase of the volume form:
\[
(\omega,\varepsilon)\longmapsto
\big(\lambda\omega,\lambda^{n/2}e^{i\theta}\varepsilon\big),
\qquad \lambda>0,\quad \theta\in\R.
\]
The remaining first-order directions are the complex deformations identified in Remark~\ref{rem:moduli}.
Thus the displayed dimensions refer to Moriyama's moduli problem, not to one with fixed basic K\"ahler class and a normalised volume form.
The extra classes counted by $h^{1,1}_{CY}-1$ belong to twisted sectors, rather than to primitive basic K\"ahler variations on the link (Remark~\ref{rem:mirror}).
\end{remark}

\subsection{The six remaining weight systems}
\label{sec:six}
\begin{figure}[ht]
\centering
\includegraphics[width=0.85\textwidth]{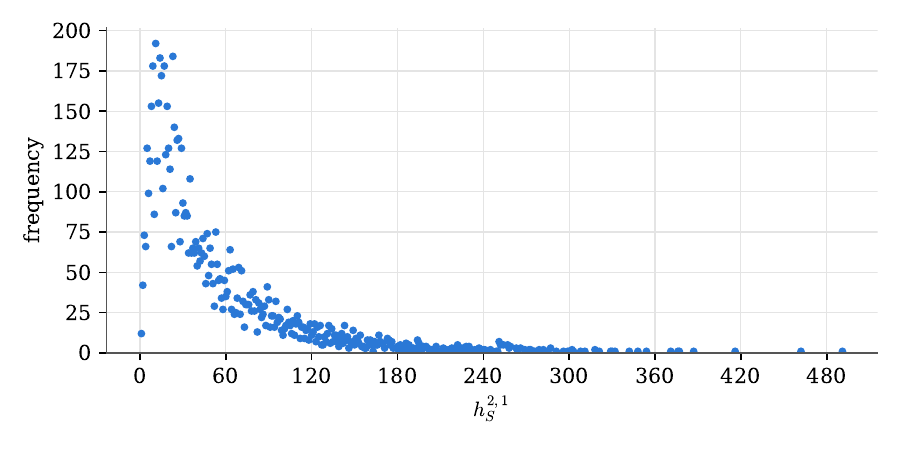}
\caption{Distribution of the Sasakian Hodge number $h^{2,1}_S$ of the $7$-dimensional links of the $7555$ weighted Calabi--Yau threefolds: each point gives a value of $h^{2,1}_S$ and the number of weight systems attaining it.
The remaining Sasakian Hodge numbers are $h^{0,0}_S=h^{3,0}_S=1$ and zero (Remark~\ref{rem:diamond}).}
\label{fig:cy3-dist}
\end{figure}

In \cite{AHHHSS2024} the invariants of $6$ of the $7555$ weight systems were left uncomputed, since the Gr\"obner-basis computations for them repeatedly exceeded the available walltime, and machine-learning predictions were given in their place.
By Theorem~\ref{cor:conj7} and Corollary~\ref{cor:closed} these invariants are determined by the weights, and they are obtained from \eqref{eq:incl-excl} in a fraction of a second.
They are listed in Table~\ref{tab:six}, alongside the predictions of \cite{AHHHSS2024}, and they complete the dataset of Sasakian Hodge numbers of $7$-dimensional Calabi--Yau links.
The distribution of $h^{2,1}_S$ over all $7555$ weight systems is shown in Figure~\ref{fig:cy3-dist}.
It takes $287$ distinct values between $1$ and $491$, with median $33$; the maximum is attained by $(1,1,12,28,42)$, one of the six systems of Table~\ref{tab:six}, and exceeds the largest value $416$ among the $7549$ systems of \cite{AHHHSS2024}.

\begin{table}[H]
\centering
\begin{tabular}{lrrrrrr}\toprule
$\bw$ & $h^{2,1}_S$ & $h^{2,1}_{CY}$ & $\mu(\bw)$ & $\nu$ & NN prediction & SR prediction\\\midrule
$(1,1,8,19,28)^\dagger$ & 348 & 348 & 39788 & 27 & 257 & 338\\
$(1,1,9,21,32)$ & 387 & 387 & 49665 & 31 & 272 & 377\\
$(1,1,11,26,39)^\dagger$ & 462 & 462 & 72226 & 35 & 293 & 447\\
$(1,1,12,28,42)^\dagger$ & 491 & 491 & 82668 & 37 & 299 & 476\\
$(1,6,34,81,122)^\dagger$ & 245 & 246 & 119805 & 19 & 243 & 243\\
$(1,6,40,93,140)$ & 274 & 275 & 153714 & 19 & 265 & 272\\\bottomrule
\end{tabular}
\caption{The invariants of the six $7$-dimensional Calabi--Yau links left open in \cite{AHHHSS2024}, computed from Corollary~\ref{cor:closed} and \eqref{eq:nuSt}; in all cases $h^{3,0}_S=1$.
The last two columns are the neural-network and symbolic-regression predictions of $h^{2,1}_S$ from \cite{AHHHSS2024}.
The rows marked $\dagger$ have since been confirmed by Gr\"obner bases \cite{MLcCY7github}; every row has been confirmed by the exact certificates described in the text.
The two cases with $h^{2,1}_S<h^{2,1}_{CY}$ each contain a curve of genus one whose generic points are $A_1$ singularities, in accordance with Theorem~\ref{thm:conj8}.}
\label{tab:six}
\end{table}

\paragraph{Admissible representatives.}
The exact checks require explicit polynomials whose critical point is isolated.
The following elementary criterion supplies the representatives used here and in Subsection~\ref{sec:9d}.

\begin{lemma}[Fermat, chain and loop polynomials]
\label{lem:invertible}
Consider polynomials in $N$ complex variables $z_1,\dots,z_N$.
Let $j$ be a map from $\{1,\dots,N\}$ to itself such that every $k$ has at most one preimage other than $k$ itself, and let $a_1,\dots,a_N\ge2$.
Then
\[f=\sum_{i:\,j(i)=i}z_i^{a_i}+\sum_{i:\,j(i)\neq i}z_i^{a_i}z_{j(i)}\]
has an isolated critical point at the origin.
If $f$ is $\bw$-homogeneous of degree $d=\sum_iw_i$ with $N\ge3$, the condition $a_i\ge2$ holds automatically.
\end{lemma}
\begin{proof}
Every component of the graph of $j$ contains one cycle.
If the cycle is a fixed point $k$, the condition on preimages makes the component a path $i_1\to\dots\to i_m=k$; otherwise every vertex of the cycle already has a preimage on the cycle, and the component is the cycle.
Hence $f$ is a sum, in disjoint sets of variables, of terms of the following forms:
\[
\begin{array}{ll}
\text{Fermat:}&x^a,\\
\text{chain:}&x_1^{a_1}x_2+\dots+x_{m-1}^{a_{m-1}}x_m+x_m^{a_m},\\
\text{loop:}&x_1^{a_1}x_2+\dots+x_{m-1}^{a_{m-1}}x_m+x_m^{a_m}x_1.
\end{array}
\]
It suffices to treat each summand.
For a chain put $M_i\coloneq x_i^{a_i}x_{i+1}$ for $i<m$ and $M_m\coloneq x_m^{a_m}$; then
\[
x_1\partial_1f=a_1M_1,\qquad x_i\partial_if=M_{i-1}+a_iM_i\quad(i\ge2).
\]
At a critical point all $M_i$ vanish, so $x_m=0$; then $\partial_mf=x_{m-1}^{a_{m-1}}$, since $a_m\ge2$, gives $x_{m-1}=0$, and in the same way $\partial_if=x_{i-1}^{a_{i-1}}+a_ix_i^{a_i-1}x_{i+1}$ gives $x_{i-1}=0$ for $i=m-1,\dots,2$.
For a loop put $M_i\coloneq x_i^{a_i}x_{i+1}$ with indices modulo $m$; then $x_i\partial_if=M_{i-1}+a_iM_i$, so at a critical point
\[
M_i=(-1)^m\big(\prod_ja_j\big)M_i,
\]
which forces $M_i=0$ as $\prod_ja_j>1$.
Hence some $x_i$ vanishes, and $\partial_if=x_{i-1}^{a_{i-1}}+a_ix_i^{a_i-1}x_{i+1}$ gives $x_{i-1}=0$, and so on around the loop.
Finally, if $d=\sum_iw_i$ and $N\ge3$, then $a_iw_i=d-w_{j(i)}>w_i$ for $j(i)\ne i$, and $a_i=d/w_i>1$ for $j(i)=i$.
\end{proof}

These are the Fermat, chain and loop building blocks of the \emph{invertible} polynomials of Kreuzer and Skarke \cite{KreuzerSkarke1992}: an admissible weighted homogeneous polynomial with as many monomials as variables, whose square exponent matrix is nonsingular over $\Q$.
Here the exponent matrix has one row for each monomial and one column for each variable; the chain and loop forms in Lemma~\ref{lem:invertible} make its nonsingularity explicit.
Invertibility refers to this matrix, not to the polynomial as a function.
The invertible polynomials used here, such as the one displayed below, were found by a randomised search for a map $j$ compatible with $(\bw,d)$, and they have all coefficients equal to $1$.

\paragraph{Verification and comparison.}
The values have been confirmed in two further ways.
Firstly, Gr\"obner-basis results for four of the six systems were added to the public dataset of \cite{MLcCY7github} in December 2024, and they agree with Table~\ref{tab:six} in both $h^{2,1}_S$ and $\nu$.
That dataset contains no values for the remaining two systems, $(1,1,9,21,32)$ and $(1,6,40,93,140)$, for which the values in Table~\ref{tab:six} are new.
Secondly, for each of the six systems the Hilbert function of $\MM_f$ was certified exactly in every degree $0\le\ell\le3d$ for an explicit invertible polynomial, such as
\[f=z_1^{280}+z_1z_4^3+z_2^{40}z_3+z_3^7+z_5^2\qquad\text{for }\bw=(1,6,40,93,140),\]
which is admissible by Lemma~\ref{lem:invertible}.
The certificate, described in Appendix~\ref{app:certificates}, is exact linear algebra over $\Q$ on the degree-$\ell$ part of the Koszul complex of $\partial_1f,\dots,\partial_5f$, checked by ranks modulo a prime; it uses neither a Gr\"obner basis nor Theorem~\ref{thm:hilbert}.
The certified Hilbert functions agree with the coefficients of $P_{\bw}$, their sums equal $\mu(\bw)$, and through \eqref{eq:nuSt} they give the values of $\nu$ in Table~\ref{tab:six}.
In addition, $h^{2,1}_S=\dim(\MM_g)_d$ was certified in the same way for a second polynomial $g$ in each case, with all monomials of degree $d$ and coefficients drawn uniformly from $\{1,\dots,999\}$, after the singularity of $g$ had been certified to be isolated (Appendix~\ref{app:certificates}).
These computations take about $80$ seconds in total.

Two features of these systems help explain the earlier predictions.
They are not those of largest Milnor number, as $\mu(\bw)$ reaches $909720$ among the systems computed in \cite{AHHHSS2024}; however, they include the two largest values of $h^{2,1}_S$ among all $7555$ weight systems, $491$ for $(1,1,12,28,42)$ and $462$ for $(1,1,11,26,39)$, whereas the largest value in the dataset of \cite{AHHHSS2024} is $416$.
The neural network underestimates all six values, by as much as $192$, and two of them lie outside the range of the data on which it was trained.
The symbolic-regression formula, which was observed in \cite{AHHHSS2024} to predict values very close to the bound $h^{2,1}_{CY}$, is within $15$ of the exact value in every case.
This behaviour is explained by the exact values: the bound $h^{2,1}_S\le h^{2,1}_{CY}$ of Conjecture~\ref{conj:cyb}, proved in Section~\ref{sec:conj8}, is attained in four of the cases and missed by one in the other two, where the defect comes from the genus-one curves $X_f\cap\WP(6,34,122)$ and $X_f\cap\WP(6,40,140)$, whose generic points have stabiliser $\mu_2$.

\section{Proof of the Calabi--Yau Bound Conjecture}
\label{sec:conj8}

Throughout this section $(\bw,d)$ is Calabi--Yau, $d=\sum_iw_i$, $f\in U_{\bw}$, $\cX\coloneq \cX_f$ and $X\coloneq X_f$.
Building on Lemmas~\ref{lem:basic} and~\ref{lem:offmiddle}, we compare the three Hodge theories appearing in \eqref{eq:intro-chain}.
The idea of the proof is that the Calabi--Yau Hodge numbers of the Kreuzer--Skarke list are \emph{stringy} Hodge numbers of $X$, which decompose into the Hodge numbers of the orbifold $\cX$ and contributions from its loci of orbifold points; the former are the Sasakian Hodge numbers, up to the diagonal correction of \eqref{eq:intro-chain}.

\subsection{Orbifold and stringy Hodge numbers}
We use the fixed loci $S(t)$, $Y^t$, $\cX_t$ and $X^t$ of Subsection~\ref{sec:orbifold-background}, and the age convention fixed there.
Theorem~\ref{thm:orbifold-stringy} will identify the sector sum with the stringy Hodge numbers once the local tangent actions have been checked.
\paragraph{Ages and fixed loci.}
The \emph{age} is the sum of the fractional rotation angles of the stabiliser action on the tangent space.
It shifts both Hodge indices of a sector in \eqref{eq:orbifold}; in particular, an age-one sector contributes its $(p-1,q-1)$-forms to bidegree $(p,q)$.
The following tangent-space calculation fixes this shift directly from the weights.

\begin{lemma}[ages]
\label{lem:age}
Let $f\in U_{\bw}$ have Calabi--Yau degree $d=\sum_iw_i$, and write $Y\coloneq f^{-1}(0)\setminus0$ and $\cX\coloneq [Y/\C^*]$.
Let $z\in Y$ and let $t=e^{2\pi i\theta}$ belong to the stabiliser of $z$.
The eigenvalues of $t$ on the tangent space of $\cX$ at $z$ are the multiset $\{t^{w_i}\}_{i=1}^N$ with one eigenvalue $t^d$ and one eigenvalue $1$ removed.
Hence the age of $t$, the sum of the numbers $\vartheta\in[0,1)$ such that $e^{2\pi i\vartheta}$ is an eigenvalue, is
\[\age(t)=\sum_{i=1}^N\langle\theta w_i\rangle-\langle\theta d\rangle,\qquad\langle x\rangle\coloneq x-\lfloor x\rfloor,\]
which depends only on $t$, and equals $\sum_{i\notin S(t)}\langle\theta w_i\rangle$ if $t^d=1$.
\end{lemma}
\begin{proof}
The tangent space of $\cX$ at $z$ is $\ker(df_z)/\C E_{\bw}(z)$, with $E_{\bw}(z)=(w_iz_i)_i$.
Differentiating $f(t\cdot z)=t^df(z)$ at the fixed point $z$ gives
\[
df_z\circ t=t^d\,df_z.
\]
The map $df_z$ is surjective because $Y$ is smooth, and its one-dimensional target has character $t^d$ under the finite cyclic group generated by $t$.
Representations of this group over $\C$ split by averaging: choose a linear splitting and average it over the group to make it equivariant.
Hence the eigenvalue multiset of the ambient tangent space is the union of those of the kernel and the target.
Also $E_{\bw}(z)$ is nonzero and fixed, because $z$ is fixed and the weighted action commutes with its Euler field.
These facts justify the following two removals, including when an eigenvalue has multiplicity greater than one.
Thus passing to $\ker df_z$ removes one eigenvalue $t^d$ from those of $t$ on $\C^N$, and quotienting by the fixed line $\C E_{\bw}(z)$ removes one eigenvalue $1$.
\end{proof}

\begin{example}[an age-one curve]
\label{ex:octic-age}
Consider the admissible polynomial
\[
f=z_1^8+z_2^8+z_3^4+z_4^4+z_5^4
\quad\text{of degree }8\text{ and weights }(1,1,2,2,2).
\]
The element $t=-1$ fixes the last three coordinates.
Its ambient eigenvalues are $(-1,-1,1,1,1)$; Lemma~\ref{lem:age} removes two eigenvalues $1$, leaving $(-1,-1,1)$ and age $\tfrac12+\tfrac12=1$.
The fixed curve has coarse space the smooth plane quartic
$z_3^4+z_4^4+z_5^4=0$, after dividing its weights by $2$.
Its three holomorphic one-forms therefore contribute to bidegree $(2,1)$ in the sector formula \eqref{eq:orbifold}.
Equivalently, the genus count is $p_{(2,2,2)}(2)=3$, with monomials $z_3,z_4,z_5$.
This is the defect of $3$ discussed in Remark~\ref{rem:classical}; it is a contribution from the quotient singularities, additional to the untwisted forms detected by the link.
\end{example}

For $F\subseteq\{1,\dots,N\}$ write $g_F\coloneq \gcd(w_i:i\in F)$.
Let $J\coloneq e^{2\pi i/d}\in\C^*$, and for $k\in\Z/d$ put $F_k\coloneq S(J^k)=\{i:d\mid kw_i\}$.

We next determine which fixed loci can occur.
Part~(iv) will count the sectors supported on each codimension-two locus in the defect formula.

\begin{lemma}[the twisted sectors]
\label{lem:twisted}
Let $f\in U_{\bw}$ have degree $d=\sum_iw_i$ in $N=n+2$ variables, and write $Y\coloneq f^{-1}(0)\setminus0$ and $X\coloneq X_f$.
For a finite-order element $t\in\C^*$, set $$S(t)\coloneq \{i:t^{w_i}=1\}, \quad Y^t\coloneq Y\cap\C^{S(t)} \qandq X^t\coloneq X\cap\WP(\bw_{S(t)}).$$
Assume $t\neq1$ and $Y^t\neq\emptyset$.
As above, set $J\coloneq e^{2\pi i/d}$, $F_k\coloneq S(J^k)$ and $g_F\coloneq \gcd(w_i:i\in F)$.
\begin{enumerate}
\item[(i)] If $t^d\neq1$, then $f$ vanishes on $\C^{S(t)}$, $|S(t)|\le N/2$, and $X^t=\WP(\bw_{S(t)})$ is a weighted projective space.
\item[(ii)] If $t^d=1$, then $f|_{\C^{S(t)}}$ has an isolated critical point at the origin, so $X^t$ is a quasi-smooth hypersurface of dimension $|S(t)|-2$ in $\WP(\bw_{S(t)})$.
\item[(iii)] If $t=J^k$, then $|F_k|\neq N-1$, and $\age(t)=1$ whenever $|F_k|=n$.
\item[(iv)] If $n\ge3$, $|F|=n$ and $g_F>1$, then $g_F\mid d$, and exactly $g_F-1$ classes $k\in\Z/d$ satisfy $F_k=F$.
\end{enumerate}
\end{lemma}
\begin{proof}
\begin{enumerate}
\item[(i)] Every monomial in the variables indexed by $S(t)$ is fixed by $t$; one of degree $d$ would force $t^d=1$, so $f|_{\C^{S(t)}}\equiv0$.
Then $\C^{S(t)}\setminus0\subset Y$, and smoothness of $Y$ requires the $N-|S(t)|$ polynomials $\partial_jf|_{\C^{S(t)}}$, $j\notin S(t)$, to have no common zero on $\C^{S(t)}\setminus0$; as they are homogeneous of positive degree, this forces $N-|S(t)|\ge|S(t)|$.
Indeed, Proposition~\ref{prop:algebra-inputs}(i) says that $r$ equations in an affine space of dimension $s$ have common zero set of dimension at least $s-r$; if $r<s$, positive homogeneity then leaves nonzero common zeros.
Here $s=|S(t)|$ and $r=N-|S(t)|$.
\item[(ii)] If $t^d=1$ and $z\in Y^t$, then $df_z\circ t=df_z$, so $df_z$ lies in the span of the $dz_i$ with $i\in S(t)$, and the differential of $f|_{\C^{S(t)}}$ at $z$ is nonzero.
By Euler's identity \eqref{eq:euler} this excludes critical points other than the origin; compare \cite[Lemma~2.1.10]{FanJarvisRuan2013}.
\item[(iii)] By Lemma~\ref{lem:age} the age is $\sum_{i\notin F_k}\langle kw_i/d\rangle$, an integer since $\sum_ikw_i/d=k$.
If $|F_k|=N-1$ the single remaining fractional part would be a nonzero integer.
If $|F_k|=n$ there are two fractional parts in $(0,1)$ with integer sum, so the age is $1$.
\item[(iv)] If $g_F\nmid d$, no monomial of degree $d$ involves only the variables in $F$, so $f|_{\C^F}\equiv0$ and the argument of (i) gives $n\le N/2$, that is $n\le2$.
For $n\ge3$ therefore $g_F\mid d$.
Now $F_k\supseteq F$ if and only if $k\in\frac d{g_F}\Z$, which gives $g_F$ classes; for $k\neq0$ among them $|F_k|\in\{n,n+1,n+2\}$, the value $n+1$ is excluded by (iii), and $n+2$ forces $k=0$ since $\gcd(\bw)=1$.
\end{enumerate}
\end{proof}

Lemmas~\ref{lem:age} and~\ref{lem:twisted} determine the local structure of $\cX$.
Every $t\neq1$ in the stabiliser of a point $z\in Y$ acts non-trivially on the tangent space of $\cX$ at $z$.
Otherwise, by Lemma~\ref{lem:age}, either $t^d\neq1$ and $|S(t)|=N-1$, which contradicts Lemma~\ref{lem:twisted}(i), or $t^d=1$ and $|S(t)|=N$, so that $t=1$ as $\gcd(\bw)=1$.
The determinant of $t$ on the tangent space is $t^{\sum_iw_i-d}=1$, so no $t$ acts as a reflection, which would have exactly one eigenvalue different from $1$, equal to its determinant.
Hence $\cX$ is an effective orbifold without reflections, and every $t\neq1$ with $Y^t\neq\emptyset$ has positive integer age.
Theorem~\ref{thm:quotient-criterion} therefore gives Gorenstein canonical cyclic quotient singularities on its coarse space $X$.

We have now checked all hypotheses of Theorem~\ref{thm:orbifold-stringy}.
Thus $h^{p,q}_{st}(X)=h^{p,q}_{\rm orb}(\cX)$, with the explicit sector sum \eqref{eq:orbifold}; no crepant resolution is needed for this equality.

For $n=3$ the Kreuzer--Skarke values were computed with Vafa's Landau--Ginzburg orbifold formula \cite{Vafa1989,KlemmSchimmrigk1994,Kreuzer:1992np}, and for $n=4$ in \cite{LSW1999} by the same method.
By \cite[Prop.~4.3, Cor.~4.4]{Batyrev2020} the generating function of that formula is
\[(-u)^nE_{\rm orb}(\cX;u^{-1},v)=\sum_{p,q}(-1)^{p+q}h^{n-p,q}_{st}(X)\,u^pv^q,\]
in which the Hodge diamond is rotated; the tabulated $h^{p,1}_{CY}$ are read off from the coefficients of $u^{n-p}v$, and they are the stringy Hodge numbers $h^{p,1}_{st}(X)$.
For $n=3$ the same values are the Hodge numbers of Batyrev's toric crepant resolution \cite{Batyrev1994,Batyrev1998}.
In what follows $h^{p,q}_{CY}\coloneq h^{p,q}_{st}(X)$.

\needspace{6\baselineskip}
The untwisted term can now be read in the Sasakian notation of Definition~\ref{def:sasakian}.

\begin{lemma}[the untwisted sector]
\label{lem:untwisted}
Let $K_f$ be the link of an admissible Calabi--Yau polynomial in $N=n+2$ variables, and let $\cX\coloneq \cX_f$ be its quotient orbifold.
For $p+q=n$, the untwisted contribution to \eqref{eq:orbifold} is
\begin{equation}
\label{eq:untwisted-hodge}
h^{p,q}(\cX)=h^{p,q}_S(K_f)+\delta_{pq}.
\end{equation}
\end{lemma}
\begin{proof}
By Lemma~\ref{lem:basic}, $h^{p,q}(\cX)=h^{p,q}_B(K_f)$, and by Definition~\ref{def:sasakian} and Lemma~\ref{lem:offmiddle}, $h^{p,q}_B=h^{p,q}_S+\delta_{pq}$ for $p+q=n$.
\end{proof}

\subsection{The bound and its defect}
Lemma~\ref{lem:untwisted} identifies the part of the stringy Hodge numbers detected by the link.
The other sectors give a non-negative correction, and Lemma~\ref{lem:twisted} identifies exactly those that can contribute to $h^{n-1,1}$.

\begin{proof}[Proof of Theorem~\ref{thm:defect}]
\begin{enumerate}
\item[(i)] In \eqref{eq:orbifold}, the untwisted term is $h^{p,q}_S+\delta_{pq}$ by Lemma~\ref{lem:untwisted}, and every twisted term is a Hodge number of a compact orbifold and hence non-negative.

\item[(ii)] A twisted sector with $\age(t)=c\ge1$ contributes $h^{n-1-c,1-c}(\cX_t)$ to the $(n-1,1)$ entry, which vanishes unless $c=1$; so only the numbers $h^{n-2,0}(\cX_t)$ of sectors of age one and dimension at least $n-2$ contribute.
If $t^d\neq1$, then $X^t$ is a weighted projective space by Lemma~\ref{lem:twisted}(i), for which $h^{n-2,0}=0$ as $n\ge3$, by Theorem~\ref{thm:vmanifold-hodge}.
If $t^d=1$, then $t=J^k$ for some $k$, and the dimension constraint
\[
\dim\cX_t=|F_k|-2\ge n-2
\]
forces $|F_k|=n$, since $|F_k|=n+1$ is excluded by Lemma~\ref{lem:twisted}(iii) and $|F_k|=N$ would force $t=1$ as $\gcd(\bw)=1$; the age is then $1$.
Fix such an $F=F_k$, and put $r\coloneq n-2\ge1$ and $g\coloneq g_F$, so that $\cX_t=\mathcal C_F$.
By Lemma~\ref{lem:twisted}(ii) and~(iv), $\mathcal C_F$ is a $\mu_g$-gerbe over the quotient of the quasi-smooth hypersurface $\{f|_{\C^F}=0\}\setminus0$ by $\C^*$ acting with the weights $\bw_F/g$, which have greatest common divisor $1$, and $f|_{\C^F}$ has degree $d/g$ for these weights.
Neither the gerbe nor the rescaling changes Hodge numbers, and Lemmas~\ref{lem:lines} and~\ref{lem:forms} apply to the rescaled hypersurface with $r$ in place of $n$.
Its holomorphic $r$-forms are therefore the sections of $\cO(\kappa_F)$ with $\kappa_F=(d-\sum_{i\in F}w_i)/g$, and since $\kappa_F<d/g$ there is no correction from the equation:
\[h^{r,0}(\mathcal C_F)=p_{\bw_F/g}(\kappa_F)=p_{\bw_F}\Big(d-\sum_{i\in F}w_i\Big).\]
Summing over the $g_F-1$ elements $t=J^k$ with $F_k=F$ of Lemma~\ref{lem:twisted}(iv) gives \eqref{eq:defect}.
Finally, a point $z\in Y$ has stabiliser containing $\mu_{g_F}$ exactly when $z_i=0$ for every $i$ with $g_F\nmid w_i$, and the set of such $i$ is the complement of $F$, since by Lemma~\ref{lem:twisted}(iii) and $\gcd(\bw)=1$ no further weight is divisible by $g_F$.
\end{enumerate}
\end{proof}

The \emph{geometric genus} of an $r$-dimensional sector is $p_g=h^{r,0}$, the number of independent holomorphic top-degree forms; the preceding argument computes it by adjunction.
Thus the defect formula has two inputs: the number $g_F-1$ of age-one sectors on a locus, and the number of top-degree forms on that locus.
For $n=3$ these are curves, whose holomorphic $1$-forms contribute to stringy bidegree $(2,1)$ after the age shift.
The summand of \eqref{eq:defect} is therefore $g_F-1$ times their genus.
Here \emph{rationality} means that the smooth coarse curve is isomorphic to $\WP^1$; it places no restriction on the stabilisers retained by the orbifold curve.
This proves the Calabi--Yau Bound Conjecture and describes exactly when it is sharp.

\begin{proof}[Proof of Theorem~\ref{thm:conj8}]
The identity is Theorem~\ref{thm:defect}(ii) with $n=3$, using $h^{2,1}_{CY}=h^{2,1}_{st}(X_f)$.
The coarse space of a one-dimensional orbifold is a smooth curve, $C_F$ is connected because $h^0(\cO)=1$ by Lemma~\ref{lem:lines}, and its genus $h^{1,0}$ is $p_g(C_F)$ by the proof of Theorem~\ref{thm:defect}; hence $p_g(C_F)=0$ exactly when $C_F\cong\WP^1$.
The restricted polynomial $f|_{\C^F}$ is irreducible: if it factored into positive-degree weighted homogeneous factors, their common zero set in $\C^3$ would have dimension at least one and would consist of critical points, contradicting Lemma~\ref{lem:twisted}(ii).
It is not a coordinate variable, since its degree exceeds every weight by Lemma~\ref{lem:VJ}.
Its zero set is therefore not contained in any coordinate hyperplane.
Consequently the points of $C_F$ with all coordinates in $F$ nonzero form a dense open subset, and their stabiliser is $\mu_{g_F}$.
At such a point, write $\{1,\dots,5\}\setminus F=\{j,j'\}$ and let $t$ generate $\mu_{g_F}$.
By Lemma~\ref{lem:age} the eigenvalues of $t$ on the tangent space are $1$, along the curve, and $t^{w_j}$, $t^{w_{j'}}$, whose product is $t^{d-\sum_{i\in F}w_i}=1$.
A prime dividing $g_F$ and $w_j$ would also divide $w_{j'}$, and hence every weight, so $t^{w_j}$ is a primitive $g_F$-th root of unity and the transversal singularity is of type $A_{g_F-1}$.
Indeed, for the action $(u,v)\mapsto(\zeta u,\zeta^{-1}v)$ the invariant coordinates $U=u^{g_F}$, $V=v^{g_F}$ and $Z=uv$ satisfy $UV=Z^{g_F}$.
Example~\ref{ex:cyclic-surface} recalls why its minimal resolution has a chain of $g_F-1$ rational $(-2)$-curves, explaining the multiplicity in \eqref{eq:defect}.
\end{proof}

The identity was evaluated for all $7555$ weight systems (Appendix~\ref{app:comp}): equality holds for exactly $4202$ of them, and for $4198$ of the $7549$ considered in \cite{AHHHSS2024}.

\begin{example}
For $\bw=(1,6,34,81,122)$, of degree $244$, the only index set $F$ with $|F|=3$ and $g_F>1$ is the one of the weights $6,34,122$, with $g_F=2$.
The curve $C_F\subset\WP(6,34,122)$ has genus $p_{(6,34,122)}(82)=1$, since $82=8\cdot6+34$ is the only representation, and therefore $h^{2,1}_{CY}-h^{2,1}_S=1$, in agreement with Table~\ref{tab:six}.
The curve passes through the coordinate point of weight $34$, which lies on $X_f$ for every $f$ because $34\nmid244$, and whose stabiliser is $\mu_{34}$; the genus in Theorem~\ref{thm:conj8} is that of the whole curve.
\end{example}

\begin{remark}
\label{rem:classical}
The identity of Theorem~\ref{thm:conj8} was also verified numerically for all $7555$ weight systems, using the Kreuzer--Skarke values of $h^{2,1}_{CY}$ \cite{KSwebsite}, the Gr\"obner-basis values of $h^{2,1}_S$ from \cite{AHHHSS2024} for the $7549$ systems computed there, and Table~\ref{tab:six} for the other six.
The defect is related to the classical distinction between polynomial and non-polynomial deformations \cite{GreenHubsch1987,CdOK1995}: $h^{2,1}_S$ counts the polynomial deformations of $X_f$ (Remark~\ref{rem:moduli}), so the right-hand side of \eqref{eq:defect} is the number of the remaining ones.
For instance, for $\WP(1,1,2,2,2)[8]$ the value $h^{2,1}_{CY}=86$ consists of $83$ polynomial deformations and $3$ attached to the curve of genus $3$ whose generic points are $A_1$ singularities.
Candelas, de la Ossa and Katz \cite[\S5.2]{CdOK1995} also identify $15$ non-polynomial deformations of $\WP(1,2,2,2,7)[14]$ with the genus of its singular curve.
Theorem~\ref{thm:conj8} recovers that contribution as $(2-1)p_{(2,2,2)}(8)=15$.
More generally, a correction from codimension-two strata already appears in Batyrev's toric formula \cite[Thm.~4.3.7]{Batyrev1994}, under its reflexive-polytope hypotheses.
Theorem~\ref{thm:defect} gives the Sasakian interpretation for admissible weighted hypersurfaces: it identifies the polynomial contribution with $h^{n-1,1}_S$, expresses the correction directly in terms of weights and stabilisers, and applies without choosing a crepant resolution.

In \cite[\S5.1]{Hirst:2023kdl} an approximation to Vafa's formula, in the form of \cite{Batyrev2020}, was studied that keeps the term of age zero and the terms of maximal size; it gives the exact value of $h^{2,1}$ for $56\%$ of the $7555$ threefolds \cite[Table~5.1]{Hirst:2023kdl}.
The term of age zero equals $\sum_ka_{kd}(uv)^k$, so in the rotated grading above it contributes $a_d=h^{n-1,1}_S$ to $h^{n-1,1}$, while the terms of maximal size contribute only to diagonal entries.
For $n\ge3$ that approximation to $h^{n-1,1}$ is therefore exactly $h^{n-1,1}_S$, and the omitted terms that change it are those of the sectors $\mathcal C_F$ of Theorem~\ref{thm:defect}.
Consequently, Theorem~\ref{thm:defect} recovers the $h^{1,2}$ and $h^{1,3}$ lower bounds in \cite[eqs.~(5.6)--(5.7)]{Hirst:2023kdl}, using Hodge symmetry, and gives their exact error and equality criterion; it also gives the corresponding bound for every $n\ge3$.
The observation of \cite{Hirst:2023kdl} is thus the equality case of Theorem~\ref{thm:conj8}, with proportion $4202/7555\approx55.6\%$, and for fourfolds the proportion $48.6\%$ of Subsection~\ref{sec:9d} matches the $48.7\%$ found in \cite[Table~5.2]{Hirst:2023kdl} for a sample of about $220\,000$ weight systems.
\end{remark}

For fourfolds Theorem~\ref{thm:defect} gives the analogous statement, with surfaces in place of curves.

\begin{corollary}[$n=4$]
\label{cor:n4}
Let $K_f$ be the nine-dimensional link of an admissible Calabi--Yau polynomial in six variables, and let $X_f$ be the corresponding fourfold.
With the stringy Hodge numbers and fixed-locus notation of Theorem~\ref{thm:defect},
\[h^{3,1}_{CY}-h^{3,1}_S=\sum_{|F|=4,\ g_F>1}(g_F-1)\,p_g(\mathcal C_F)\qandq h^{2,2}_S\le h^{2,2}_{CY}-1,\]
where the surfaces $\mathcal C_F$ have coarse spaces $X_f\cap\WP(\bw_F)$ and $p_g(\mathcal C_F)=p_{\bw_F}(d-\sum_{i\in F}w_i)$.
Here
\begin{equation}
\label{eq:cy4-lw}
h^{2,2}_{CY}=44+4h^{1,1}_{CY}-2h^{2,1}_{CY}+4h^{3,1}_{CY}.
\end{equation}
\end{corollary}
\begin{proof}
The first two statements are Theorem~\ref{thm:defect}(ii) and~(i).
For the last, $X_f$ is projective with Gorenstein canonical singularities and trivial canonical class, and its stringy Hodge numbers exist, so Theorem~\ref{thm:background-lw}, with $m=4$, gives:
\begin{equation}
\label{eq:stringy-lw}
\sum_{p,q}(-1)^{p+q}(p-2)^2h^{p,q}_{st}=\frac13e_{st},\qquad
 e_{st}\coloneq \sum_{p,q}(-1)^{p+q}h^{p,q}_{st}.
\end{equation}
Twisted sectors have positive age and therefore do not contribute to $h^{p,0}_{st}$, so $h^{0,0}_{st}=h^{4,0}_{st}=1$ and $h^{p,0}_{st}=0$ for $0<p<4$ by Lemma~\ref{lem:offmiddle} and Corollary~\ref{cor:closed}.
With the symmetries $h^{p,q}_{st}=h^{q,p}_{st}=h^{4-p,4-q}_{st}$, the identity becomes the stated formula, exactly as for smooth Calabi--Yau fourfolds.
\end{proof}

Restricted to the untwisted sector, the same identity relates the two Sasakian Hodge numbers of a $9$-dimensional link, up to a correction from the twisted sectors.

\begin{proof}[Proof of Theorem~\ref{prop:lw}]
Write the decomposition \eqref{eq:orbifold} into the untwisted sector and the sum of the twisted ones as
\begin{equation}
\label{eq:twisted-hodge-sum}
\begin{aligned}
h^{p,q}_{CY}&=h^{p,q}_{\rm untw}+T^{p,q},\\
T^{p,q}&\coloneq \sum_{t\ne1}h^{p-\age(t),q-\age(t)}(\cX_t).
\end{aligned}
\end{equation}
By Lemmas~\ref{lem:untwisted} and~\ref{lem:offmiddle} the untwisted parts of $h^{1,1}$, $h^{2,1}$, $h^{3,1}$ and $h^{2,2}$ are $1$, $0$, $h^{3,1}_S$ and $h^{2,2}_S+1$.
Inserting this into \eqref{eq:cy4-lw} gives \eqref{eq:link-lw} with
\[
R=4T^{1,1}-2T^{2,1}+4T^{3,1}-T^{2,2},
\]
so that a sector of age $a$ contributes
\[\rho(t)=4h^{1-a,1-a}-2h^{2-a,1-a}+4h^{3-a,1-a}-h^{2-a,2-a}\]
of its own Hodge numbers.
By Lemma~\ref{lem:twisted}(i) and~(iii), $Y^t=Y\cap\C^{S(t)}$ with $|S(t)|\le3$ if $t^d\ne1$ and $|S(t)|\le4$ if $t^d=1$, so $m\le2$.
The eigenvalues of $t^{-1}$ on the normal space of $\cX_t$ are the inverses of those of $t$, none equal to $1$, so
\begin{equation}
\label{eq:inverse-ages}
\age(t)+\age(t^{-1})=4-m,
\end{equation}
and both ages are at least $1$.
Hence $a=1$ if $m=2$, $\{\age(t),\age(t^{-1})\}=\{1,2\}$ if $m=1$, and $a\in\{1,2,3\}$ if $m=0$.
For $m=2$, Hodge symmetry and Poincar\'e duality on the coarse space give
\[
\begin{aligned}
\rho&=4h^{0,0}-2h^{1,0}+4h^{2,0}-h^{1,1}\\
&=6(h^{0,0}-h^{0,1}+h^{0,2})-\sum_{p,q}(-1)^{p+q}h^{p,q}\\
&=6\chi(\cO)-e.
\end{aligned}
\]
For $m=1$ the sector is connected (Lemmas~\ref{lem:twisted} and~\ref{lem:lines}), so $\rho=4-2g-1=e+1$ for age one, where $g=h^{1,0}$, and $\rho=-h^{0,0}=-1$ for age two.
For $m=0$, $\rho=4h^{0,0}$, $-h^{0,0}$ or $0$ for ages $1$, $2$ or $3$.
If $X_f$ has no orbifold points there are no twisted sectors and $R=0$.
\end{proof}

The identity has been checked independently of the Kreuzer--Skarke data for $2068$ weight systems: $2000$ drawn at random and the $25$ largest values of $h^{3,1}_S$, of $R$ and of $-R$.
For each of them every twisted sector was constructed from the weights by Lemmas~\ref{lem:age} and~\ref{lem:twisted}, and its Hodge numbers were computed by Lemma~\ref{lem:lines} and Theorem~\ref{thm:sheaf} for surfaces and curves, and by counting points for sectors of dimension zero.
In every case the sum of the $\rho(t)$ equals $h^{2,2}_S-4h^{3,1}_S-47$ computed from \eqref{eq:9d}, and the resulting twisted parts reproduce the Kreuzer--Skarke values of $h^{1,1}_{CY}$, $h^{2,1}_{CY}$ and $h^{3,1}_{CY}$ (script \texttt{lw\_identity\_9d.py} of the \href{https://github.com/GeometrIALab/Link_SasakianHodgeNumbers}{GitHub} repository).

\begin{remark}[mirror symmetry]
\label{rem:mirror}
The Sasakian Hodge numbers of Calabi--Yau links do not exhibit mirror symmetry, and Theorem~\ref{thm:defect} shows why.
For Calabi--Yau $n$-folds mirror symmetry exchanges $h^{p,q}$ with $h^{n-p,q}$, and in particular $h^{n-1,1}$ with $h^{1,1}$.
For $n\ge3$, however, $h^{1,1}_S=0$ for every link by Remark~\ref{rem:diamond}, whereas $h^{n-1,1}_S\ge1$ for each of the $7555$ links with $n=3$ and the $1\,100\,055$ links with $n=4$ considered here; hence the Sasakian Hodge numbers of no link in these datasets are the mirror image of those of another link, or of its own.
The underlying reason is that the link sees only the untwisted sector of \eqref{eq:orbifold} (Lemma~\ref{lem:untwisted}), and mirror symmetry does not preserve the splitting into untwisted and twisted sectors.
To $h^{n-1,1}$ the untwisted sector contributes exactly the polynomial deformations, counted by $h^{n-1,1}_S$ (Remark~\ref{rem:moduli}), while by Lemma~\ref{lem:offmiddle} its contribution to $h^{1,1}$ is only the class $[d\eta]$, which is proportional to the hyperplane class; every further class of type $(1,1)$ comes from a twisted sector of age one, since a sector of age $c$ contributes $h^{1-c,1-c}$ to $h^{1,1}$ in \eqref{eq:orbifold}.
Suppose, for instance, that the mirror of $X_f$ is itself a quasi-smooth Calabi--Yau hypersurface $X_g$ in a weighted projective space, so that $h^{p,q}_{CY}(X_g)=h^{n-p,q}_{CY}(X_f)$.
By \eqref{eq:orbifold} and Lemma~\ref{lem:offmiddle}, $h^{1,1}_{CY}(X_g)=1+\tau(\cX_g)$, where $\tau(\cX_g)$ is the number of connected components of the twisted sectors of age one of $\cX_g$, and Theorem~\ref{thm:defect} becomes
\[h^{n-1,1}_S(K_f)+\sum_{|F|=n,\ g_F>1}(g_F-1)\,p_{\bw_F}\Big(d-\sum_{i\in F}w_i\Big)=1+\tau(\cX_g).\]
The left-hand side is computed on the link of $X_f$, while the right-hand side, apart from the hyperplane class, is invisible to the link of $X_g$, whose own $h^{1,1}_S$ vanishes.
Mirror symmetry is therefore a statement about the stringy Hodge numbers $h^{p,q}_{CY}$, not about the Sasakian Hodge numbers of the links.
At that level every quasi-smooth Calabi--Yau hypersurface $X_f$ has a mirror $X^\vee$ with $h^{p,q}_{st}(X^\vee)=h^{n-p,q}_{CY}(X_f)$ \cite{Batyrev1994,Batyrev2020}, but $X^\vee$ is a compactified toric hypersurface that is in general not a hypersurface in a weighted projective space, so the lists of weight systems are not closed under the exchange \cite{CandelasLynkerSchimmrigk1990,Kreuzer:1992np}, and in general the mirror has no link of the kind considered here.
\end{remark}

\section{Applications in dimensions five and nine}
\label{sec:applications}
The preceding results apply uniformly in the dimension of the quotient.
We now compare their geometric consequences for K3 surfaces and Calabi--Yau fourfolds: in the first case they recover the topology of the links, and in the second they give a complete dataset of middle Hodge numbers.

\subsection{Five-dimensional links}
\label{sec:5d}

We apply Corollary~\ref{cor:closed} and Theorem~\ref{thm:defect} to links whose quotient orbifold has complex dimension $n=2$.
The Hodge numbers determine the second Betti number, while the twisted sectors record the exceptional curves of the K3 resolution.
The complete list is given in Appendix~\ref{app:k3data}.

\paragraph{Hodge numbers and exceptional curves.}
For $n=2$ the $95$ Calabi--Yau weight systems of quasi-smooth K3 surfaces \cite{Yonemura1990,IanoFletcher2000} give $5$-dimensional links $K_f$.
They are simply connected Sasakian $5$-manifolds whose leaf space is a K3 orbifold.
By Corollary~\ref{cor:closed} and Lemma~\ref{lem:offmiddle} their Sasakian Hodge numbers and second Betti number are
\begin{equation}
\label{eq:k3}
h^{2,0}_S=h^{0,2}_S=1,\qquad h^{1,1}_S=p_{\bw}(d)-\sum_{i=1}^4p_{\bw}(w_i),\qquad b_2(K_f)=2+h^{1,1}_S .
\end{equation}
Together with $h^{0,0}_S=1$ and $h^{1,0}_S=h^{0,1}_S=0$ these are all the Sasakian Hodge numbers of the link (Remark~\ref{rem:diamond}).
The number $h^{1,1}_S$ is related to the singularities of the K3 surface $X_f$ by a simple count.
Every singular point $x$ of $X_f$ is a cyclic quotient singularity $\C^2/\mu_{k_x}$ with $\mu_{k_x}\subset\SL_2(\C)$, by Lemma~\ref{lem:age} and the local-structure discussion following Lemma~\ref{lem:twisted}, that is, a Du Val singularity of type $A_{k_x-1}$, whose minimal resolution has $k_x-1$ exceptional curves.

\begin{proposition}
\label{prop:k3}
Let $K_f$ be the five-dimensional link of an admissible Calabi--Yau polynomial in four variables, and let $X_f$ be its K3 coarse space.
For each singular point $x\in X_f$, write $A_{k_x-1}$ for its Du Val type.
Then
\begin{equation}
\label{eq:k3-exceptional}
20=1+h^{1,1}_S(K_f)+\sum_{x\in\mathrm{Sing}X_f}(k_x-1).
\end{equation}
\end{proposition}
\begin{proof}
The minimal resolution of $X_f$ is crepant and is a K3 surface, so $h^{1,1}_{st}(X_f)=20$ \cite[Thm.~3.12]{Batyrev1998}.
In \eqref{eq:orbifold} the untwisted sector contributes $h^{1,1}(\cX_f)=1+h^{1,1}_S$, as in Lemma~\ref{lem:untwisted}.
Each $t\neq1$ in the stabiliser $\mu_{k_x}$ of a point $x$ acts on the tangent plane non-trivially and with determinant $1$, so without eigenvalue $1$; its fixed points are therefore isolated, and it has age $1$.
The twisted sectors are thus points, one for each pair $(x,t)$ with $t\in\mu_{k_x}\setminus\{1\}$, and each contributes $h^{0,0}=1$ to $h^{1,1}_{st}$.
Summing over all singular points gives $\sum_x(k_x-1)$.
\end{proof}

\begin{figure}[ht]
\centering
\includegraphics[width=0.6\textwidth]{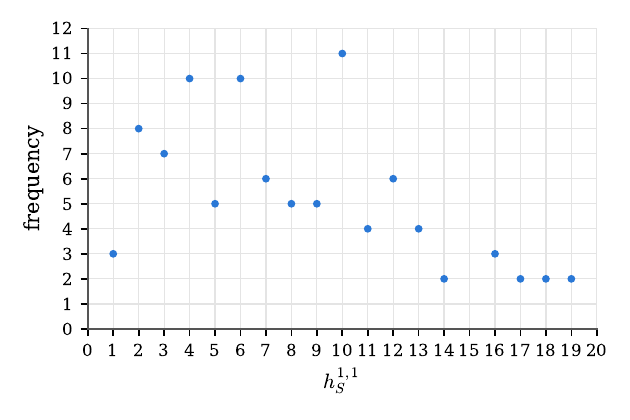}
\caption{Distribution of the Sasakian Hodge number $h^{1,1}_S$ of the $5$-dimensional links of the $95$ weighted K3 surfaces: each point gives a value of $h^{1,1}_S$ and the number of weight systems attaining it.
The remaining Sasakian Hodge numbers are $h^{0,0}_S=h^{2,0}_S=1$ and $h^{1,0}_S=0$ (Remark~\ref{rem:diamond}).}
\label{fig:k3-dist}
\end{figure}

\paragraph{The dataset and its topology.}
The data are collected in Table~\ref{tab:k3} in Appendix~\ref{app:k3data}.
In the authors' follow-up computations $(h^{2,0}_S,h^{1,1}_S)$ was computed with Gr\"obner bases for one explicit polynomial for each of the $95$ weight systems (file \texttt{K3Topological\_Data.txt} of the \href{https://github.com/GeometrIALab/Link_SasakianHodgeNumbers}{GitHub} repository), and these values agree with \eqref{eq:k3} in all $95$ cases.
Each of these polynomials was certified to have an isolated singularity (Appendix~\ref{app:certificates}), and by Theorem~\ref{cor:conj7} the values hold for every admissible polynomial.
The singular points were read off from the weights.
For generic $f$ and each set $S$ of coordinates with $g_S>1$, the points of $X_f$ whose nonzero coordinates are exactly those in $S$ number $p_{\bw_S}(d)-1$ when $|S|=2$, and one or zero when $|S|=1$ according to whether $w_i\nmid d$ or $w_i\mid d$; the positions of these points depend on $f$.
The identity of Proposition~\ref{prop:k3} holds for all $95$ systems.
Across the dataset $h^{1,1}_S$ takes every value from $1$ to $19$ except $15$, with the distribution shown in Figure~\ref{fig:k3-dist}.

For example the K3 surface $X_{27}\subset\WP(5,6,7,9)$ has singularities $A_6+A_5+A_4+A_2$, hence $17$ exceptional curves, and its link has $h^{1,1}_S=2$.

\begin{remark}[topology of K3 links]
\label{rem:k3top}
The topology of these links is known \cite{BGM2006}.
A K3 link is spin, since it is the boundary of the Milnor fibre, which is parallelisable, so that $K_f$ is stably parallelisable.
Its second homology is torsion-free \cite[Thm.~39]{BGM2006}, so by Smale's classification of simply connected spin $5$-manifolds \cite{Smale1962}
\[H_2(K_f;\Z)\cong\Z^{\,2+h^{1,1}_S},\qquad K_f\cong\#_{\,2+h^{1,1}_S}(S^2\times S^3),\]
as observed in \cite[Ex.~41]{BGM2006}.
Table~\ref{tab:k3} recovers the range $3\le b_2\le21$, $b_2\ne17$, recorded in \cite[Cor.~42]{BGM2006}; the absent value $b_2=17$ corresponds to the absent value $h^{1,1}_S=15$.
\end{remark}

\subsection{Nine-dimensional links}
\label{sec:9d}

We now use Corollary~\ref{cor:closed} to compute the middle Sasakian Hodge numbers for the complete fourfold list.
Corollary~\ref{cor:n4} provides the comparison with stringy Hodge numbers, and the independent checks below test the calculation on weights and on explicit polynomials.

For $n=4$ the link $K_f$ is a $3$-connected $9$-manifold whose leaf space is a Calabi--Yau fourfold orbifold $\cX_f$, with coarse space $X_f\subset\WP^5(\bw)$.
By Corollary~\ref{cor:closed} and \eqref{eq:middle-betti}, its Sasakian Hodge numbers and middle Betti number are
\begin{equation}
\label{eq:9d}
\begin{aligned}
h^{4,0}_S&=1,\\
h^{3,1}_S&=p_{\bw}(d)-\sum_ip_{\bw}(w_i),\\
h^{2,2}_S&=p_{\bw}(2d)-\sum_ip_{\bw}(d+w_i)+\sum_{i<j}p_{\bw}(w_i+w_j).
\end{aligned}
\end{equation}
\[b_4(K_f)=2+2h^{3,1}_S+h^{2,2}_S .\]
By Remark~\ref{rem:diamond} these, together with $h^{1,3}_S=h^{3,1}_S$ and $h^{0,4}_S=1$, are all the non-trivial Sasakian Hodge numbers: $h^{0,0}_S=1$, every $h^{p,q}_S$ with $0<p+q<4$, such as $h^{1,1}_S$ and $h^{2,1}_S$, vanishes, and the harmonic forms of degree $5$ to $9$ are the Hodge duals of those of degree $4$ to $0$.

\paragraph{Computation.}
These numbers were computed for all $1\,100\,055$ non-degenerate Calabi--Yau weight systems of \cite{LSW1999}, as listed by Kreuzer and Skarke \cite{KSfourfolds}, whose degrees range up to $d=6\,521\,466$.
Each Hodge number requires at most $22$ denumerants $p_{\bw}(k)$.
These were computed exactly by enumerating the multiplicities of the four largest weights and counting the solutions for the two smallest weights in closed form.
For positive integers $a,b$, the number of $(x,y)\in\Z_{\ge0}^2$ with $ax+by=r$ is zero if $r<0$ or $c\coloneq \gcd(a,b)$ does not divide $r$.
Otherwise, put
\[
a'\coloneq a/c,\qquad b'\coloneq b/c,\qquad r'\coloneq r/c,
\]
and choose $0\le x_0<b'$ satisfying
\[
a'x_0\equiv r'\pmod{b'}.
\]
When $a'x_0\le r'$, the number of solutions is
\begin{equation}
\label{eq:two-weight-count}
\left\lfloor\frac{r'-a'x_0}{a'b'}\right\rfloor+1.
\end{equation}
It is zero when $a'x_0>r'$; when $b'=1$, the congruence simply gives $x_0=0$.
The computation of all invariants, together with the twisted surfaces and the comparison (a) below, took $322$\,s of wall-clock time on a laptop with eight processes (Appendix~\ref{app:certificates}).
By contrast, the Gr\"obner-basis approach of \cite{AHHHSS2024} enumerates a monomial basis of the whole Milnor algebra, whose dimension $\mu(\bw)$ has median $1.1\times10^7$ over this dataset, exceeds $10^8$ for $23\%$ and $10^9$ for $7\%$ of the weight systems, and reaches $3.0\times10^{12}$.
Gr\"obner-basis computations with the code of \cite{AHHHSS2024} for these links had not finished after a month on a computing cluster, which these sizes explain.

\paragraph{Verification.}
Since no Gr\"obner-basis values were available in this dimension, the results were checked by four methods, which test different parts of the computation.
\begin{enumerate}
\item[(a)] \emph{The Kreuzer--Skarke data.}
The defect formula of Corollary~\ref{cor:n4}, evaluated on the twisted surfaces, reproduces $h^{3,1}_{CY}-h^{3,1}_S$ exactly for all $1\,100\,055$ systems, where $h^{3,1}_{CY}$ is the value computed in \cite{LSW1999} with Vafa's formula and tabulated in \cite{KSfourfolds}.
The bound $h^{2,2}_S\le h^{2,2}_{CY}-1$ also holds in every case, with $h^{2,2}_{CY}$ obtained from \eqref{eq:cy4-lw}.
This compares $h^{3,1}_S$ with independently computed data, but it uses the same denumerant code as the computation.
\item[(b)] \emph{Topology.}
Milnor and Orlik computed the middle Betti number of a link from the monodromy of the Milnor fibration \cite{MilnorOrlik1970}: writing $d/w_i=u_i/v_i$ in lowest terms,
\[b_{N-2}(K_f)=\sum_{I\subseteq\{1,\dots,N\}}(-1)^{N-|I|}\,\frac{\prod_{i\in I}u_i/v_i}{\mathrm{lcm}(u_i:i\in I)} .\]
This formula involves no denumerants and was implemented separately.
It agrees with $b_4=2+2h^{3,1}_S+h^{2,2}_S$ for all $1\,100\,055$ systems, and it reproduces $b_3=2+2h^{2,1}_S$ for all $7549$ Gr\"obner-basis computations of \cite{AHHHSS2024}.
It checks the sum of the Hodge numbers rather than each of them; together with (a), it determines $h^{2,2}_S$ independently of its own formula.
\item[(c)] \emph{Explicit polynomials.}
These checks use neither the denumerant code nor Theorem~\ref{thm:hilbert}.
For each of the $42$ weight systems of the authors' follow-up computations (file \texttt{CY4Polynomials\_42.txt} of the \href{https://github.com/GeometrIALab/Link_SasakianHodgeNumbers}{GitHub} repository), a Gr\"obner basis of $\cJ_f$ over $\Q$ was computed for an invertible polynomial (Lemma~\ref{lem:invertible}), and all standard monomials were enumerated.
In every case their number is $\mu(\bw)$, the Hilbert function equals $P_{\bw}(t)$ in every degree, and $(h^{4,0}_S,h^{3,1}_S,h^{2,2}_S)$ agrees with \eqref{eq:9d}.
For larger cases, $\dim(\MM_f)_\ell$ was certified in single degrees by the exact certificate of Appendix~\ref{app:certificates}.
In this way $h^{3,1}_S$ was certified for $900$ weight systems that admit an invertible polynomial, $100$ drawn at random from each range $10^k\le\mu(\bw)<10^{k+1}$ with $4\le k\le12$, and $h^{2,2}_S$ for the $655$ of them with $\dim S_{2d}\le4\times10^5$.
Both numbers were also certified for $25$ randomly drawn weight systems that admit no invertible polynomial and have $\dim S_{2d}\le2500$ and $\mu(\bw)\le10^7$, using polynomials with all monomials of degree $d$ and coefficients drawn uniformly from $\{1,\dots,999\}$, and $h^{3,1}_S$ was certified for the $42$ dense polynomials of the follow-up computations, which contain all monomials of degree $d$ with coefficients between $1$ and $4$.
The invertible polynomials are admissible by Lemma~\ref{lem:invertible}, and the dense ones were certified separately to have isolated singularities (Appendix~\ref{app:certificates}).
Every certified value agrees with \eqref{eq:9d}.
\item[(d)] \emph{The sheaf-cohomology formula.}
Formula \eqref{eq:sheaf} of Theorem~\ref{thm:sheaf}, evaluated term by term, reproduces both $h^{3,1}_S$ and $h^{2,2}_S$ for all $1\,100\,055$ systems.
This detects slips in the algebra of Corollary~\ref{cor:closed} and of the telescoping, but it shares the denumerant code with the computation.
\end{enumerate}

\begin{figure}[ht]
\centering
\includegraphics[width=0.7\textwidth]{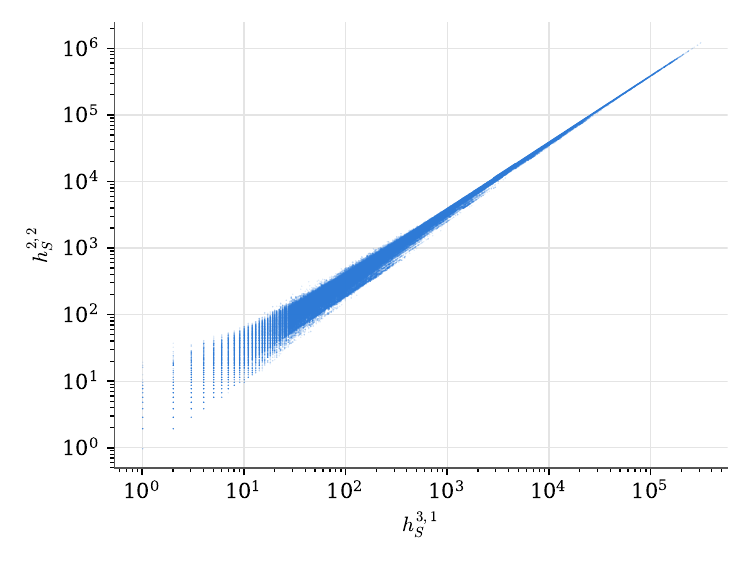}
\caption{The Sasakian Hodge numbers $(h^{3,1}_S,h^{2,2}_S)$ of the $9$-dimensional links of the $1\,100\,055$ weighted Calabi--Yau fourfolds, on logarithmic axes; each point is a weight system, drawn semi-transparently so that darker regions contain more of them.
The remaining Sasakian Hodge numbers are $h^{0,0}_S=h^{4,0}_S=1$ and zero (Remark~\ref{rem:diamond}).}
\label{fig:cy4-dist}
\end{figure}

\paragraph{Results.}
Table~\ref{tab:9d-stats} summarises the distribution of the invariants, and Table~\ref{tab:9d-examples} lists a selection of links, including the extreme cases.
Equality $h^{3,1}_S=h^{3,1}_{CY}$ holds for $534\,935$ weight systems ($48.6\%$), which are exactly those with no twisted surface of positive geometric genus.
Twisted surfaces occur in $930\,088$ systems.
The defect is positive in the remaining $565\,120$ cases, with maximum $8500$ for $\WP(2,2,70,294,367,735)$, where the surface $X_f\cap\WP(2,2,70,294)$, whose generic points have stabiliser $\mu_2$, has $p_g=8500$.
Equality $h^{2,2}_S=h^{2,2}_{CY}-1$ holds only for the sextic $\WP(1,1,1,1,1,1)[6]$ and for $\WP(1,1,1,1,1,5)[10]$; in both cases $X_f$ is smooth, since it avoids the orbifold points of the ambient space, so that there are no twisted sectors.
The pair $(h^{3,1}_S,h^{2,2}_S)$ takes $477\,048$ distinct values, plotted in Figure~\ref{fig:cy4-dist}.
The two numbers are strongly correlated: for the $364\,766$ weight systems with $h^{3,1}_S\ge1000$ the ratio $h^{2,2}_S/h^{3,1}_S$ lies between $2.59$ and $4.16$, and it narrows to $[3.99,4.01]$ for the $1273$ systems with $h^{3,1}_S\ge10^5$, so that for large $h^{3,1}_S$ the points approach the line $h^{2,2}_S=4h^{3,1}_S$.
This is explained by the identity \eqref{eq:link-lw} of Theorem~\ref{prop:lw}.
The slope $4$ is the coefficient of $h^{3,1}$ in the Libgober--Wood identity, a consequence of Riemann--Roch for $c_1=0$, and
\[\frac{h^{2,2}_S}{h^{3,1}_S}=4+\frac{47+R}{h^{3,1}_S},\]
where $R$ is a sum of Euler-characteristic-type invariants of the twisted sectors, which have dimension at most two, whereas $h^{3,1}_S$ is a Hodge number of the fourfold.
Over the whole dataset $-10\,177\le R\le533$.
The large negative values come from twisted surfaces, for which $6\chi(\cO)-e$ is negative and of the order of their geometric genus.
For $\WP(2,2,90,281,750,1125)$, for instance, the surface $X_f\cap\WP(2,2,90,750)$ has $p_g=7236$ and contributes $-10\,485$ to $R=-10\,177$, so that $h^{2,2}_S/h^{3,1}_S=140\,494/37\,656\approx3.73$.
A twisted surface of positive genus lies in a weighted projective space $\WP(\bw_F)$ that omits two weights $w_a,w_b$ not divisible by $g_F$, and its genus $p_g=p_{\bw_F}(w_a+w_b)$ counts monomials of the small degree $w_a+w_b$.
For the $1273$ links with $h^{3,1}_S\ge10^5$, whose smallest weight is $1$ or $2$, these genera are small, with $\sum_F(g_F-1)p_g\le116$, and $-146\le47+R\le556$, which is less than $0.5\%$ of $h^{3,1}_S$.
Since there are only finitely many weight systems, this approach to the line is a property of the dataset rather than a limit theorem.

The full dataset, containing for each weight system its degree, the Kreuzer--Skarke Hodge numbers, $h^{3,1}_S$, $h^{2,2}_S$, $b_4(K_f)$, the exact value of $\mu(\bw)$ and the twisted surfaces, is available in the \href{https://github.com/GeometrIALab/Link_SasakianHodgeNumbers}{GitHub} repository.

\begin{table}[ht]
\centering
\begin{tabular}{lrrrr}\toprule
Quantity & Minimum & Median & Maximum & Mean\\\midrule
degree $d$ & 6 & 4480 & 6\,521\,466 & 30\,866\\
$h^{3,1}_S$ & 1 & 448 & 303\,148 & 2268\\
$h^{2,2}_S$ & 1 & 1715 & 1\,213\,142 & 8989\\
$b_4(K_f)$ & 5 & 2610 & 1\,819\,440 & 13\,526\\
$h^{3,1}_{CY}-h^{3,1}_S$ & 0 & 1 & 8500 & 32.6\\
$\mu(\bw)$ & $1.6\times10^4$ & $1.1\times10^7$ & $3.0\times10^{12}$ & $1.2\times10^9$\\\bottomrule
\end{tabular}
\caption{Summary of the Sasakian Hodge numbers of the $1\,100\,055$ nine-dimensional Calabi--Yau links, together with the degree, the defect of Corollary~\ref{cor:n4} and the Milnor number.
Means are rounded to the nearest integer, except for the defect (one decimal) and $\mu(\bw)$, which like its other entries is given to two significant figures.
The Hodge numbers are spread over several orders of magnitude, while the defect is usually small: its median is $1$ and it vanishes for $48.6\%$ of the links.}
\label{tab:9d-stats}
\end{table}

\begin{table}[ht]
\centering\footnotesize\setlength{\tabcolsep}{4.5pt}
\begin{tabular}{lrrrrrrr}\toprule
$\bw$ & $d$ & $h^{3,1}_S$ & $h^{2,2}_S$ & $b_4(K_f)$ & $h^{3,1}_{CY}$ & $h^{2,2}_{CY}$ & $\mu(\bw)$\\\midrule
$(1,1,1,1,1,1)$ & 6 & 426 & 1751 & 2605 & 426 & 1752 & $1.6\times10^4$\\
$(1,1,1,1,1,2)$ & 7 & 455 & 1866 & 2778 & 455 & 1868 & $1.9\times10^4$\\
$(1,1,1,1,2,2)$ & 8 & 443 & 1815 & 2703 & 443 & 1820 & $2.2\times10^4$\\
$(1,1,2,2,2,2)$ & 10 & 346 & 1406 & 2100 & 350 & 1452 & $2.1\times10^4$\\
$(2,2,2,3,3,3)$ & 15 & 187 & 788 & 1164 & 187 & 792 & $1.8\times10^4$\\
$(2,2,70,294,367,735)$ & 1470 & 20804 & 75830 & 117440 & 29304 & 117772 & $1.3\times10^8$\\
$(1,1,84,516,1204,1806)$ & 3612 & 303148 & 1213142 & 1819440 & 303148 & 1213644 & $6.6\times10^9$\\
$(105,230,279,462,534,805)$ & 2415 & 1 & 18 & 22 & 1 & 1140 & $4.8\times10^4$\\
$(521,600,624,625,630,750)$ & 3750 & 1 & 1 & 5 & 1 & 1752 & $1.6\times10^4$\\
$(1805,1806,151662,931638,2173822,3260733)$ & 6521466 & 252 & 502 & 1008 & 252 & 1213644 & $6.6\times10^9$\\
$(1,1806,75894,466206,1087814,1631721)$ & 3263442 & 151700 & 606822 & 910224 & 151700 & 1213644 & $3.0\times10^{12}$\\\bottomrule
\end{tabular}
\caption{Selected nine-dimensional Calabi--Yau links, with the stringy Hodge numbers of the fourfold and the Milnor number.
The rows comprise the sextic and four further links of small degree, among them $(1,1,2,2,2,2)$, whose twisted quintic surface $X_f\cap\WP(2,2,2,2)$ has $p_g=4$ and accounts for the defect of $4$; the link with the largest defect; the maxima of $h^{3,1}_S$, $h^{2,2}_S$, $d$ and $\mu$; and links with $h^{3,1}_S=1$ and with $h^{2,2}_S=1$, the latter over the mirror of the sextic, with $(h^{1,1}_{CY},h^{3,1}_{CY})=(426,1)$.}
\label{tab:9d-examples}
\end{table}

\newpage
\section*{Afterword}
\phantomsection
\addcontentsline{toc}{section}{Afterword}
\label{sec:conclusions}

We return to the relation between the structural results and the computational observations that motivated them, before indicating the remaining topological questions.
The two conjectures of \cite{AHHHSS2024} are proved in Theorem~\ref{cor:conj7} and Theorem~\ref{thm:conj8}.
The Weight Determinacy Theorem makes the isolated-singularity hypothesis explicit and holds in every degree of the Milnor algebra, because the Hilbert series of a weighted homogeneous isolated singularity is determined by its weights.
Theorem~\ref{thm:sheaf} gives an alternative proof of the Milnor-algebra description used in \cite{AHHHSS2024}, through sheaf cohomology on the leaf space.
Theorem~\ref{thm:defect} bounds the middle Sasakian Hodge numbers by the stringy Hodge numbers in every dimension $n\ge2$, with a correction on the diagonal; for $n\ge3$ the gap in $h^{n-1,1}$ is an explicit sum of geometric genera of twisted sectors, and for threefolds the bound is sharp exactly when all curves of orbifold points are rational.

\paragraph{Interpretation of the experiments.}
The machine-learning results of \cite{AHHHSS2024} can be read in this light.
The Sasakian Hodge number $h^{2,1}_S$ is an alternating sum of six lattice-point counts in the weights, which may be part of the reason why it could be learnt accurately, although such counts depend on arithmetic properties of the weights and are not smooth functions of them.
Its correlation with $h^{2,1}_{CY}$ is explained by the defect, which vanishes for $56\%$ of the weight systems.
The same formulae give complete datasets of Sasakian Hodge numbers, here for all $95$ K3 links and all $1\,100\,055$ Calabi--Yau fourfold links, the latter beyond the reach of the Gr\"obner-basis pipeline of \cite{AHHHSS2024}.
These datasets are complete in a strong sense, since all Sasakian Hodge numbers outside the middle degree are trivial.
For the same reason they cannot exhibit mirror symmetry: the links see only the untwisted sector of the stringy Hodge numbers, and mirror symmetry exchanges the polynomial deformations counted by $h^{n-1,1}_S$ with classes of type $(1,1)$ on the mirror that come from its twisted sectors (Remark~\ref{rem:mirror}).

\paragraph{Further questions.}
Natural next steps are to determine the torsion and diffeomorphism types of the $7$- and $9$-dimensional links, the $5$-dimensional ones being classified \cite{BGM2006,Smale1962}; to extend the Sasakian Hodge number database to higher-dimensional links; to look for analogues of the Crowley--Nordstr\"om invariant for the $\SU(n)$-structures of contact Calabi--Yau manifolds of higher dimension; to compare the invariants of the links, including $\nu$, for pairs of weight systems whose Calabi--Yau hypersurfaces are mirror to each other; and to use the new datasets as testing grounds for machine learning in geometry.

\appendix
\section{Certificates and computations}
\label{app:comp}
This appendix records the computational checks used in Sections~\ref{sec:sheaf}, \ref{sec:six}, \ref{sec:5d} and~\ref{sec:9d}, together with the computing environment and the complete K3 dataset.
We retain the weighted grading $\deg z_i=w_i$ and the Jacobian ideal $\cJ_f$ of Section~\ref{sec:weights}; admissibility of the invertible representatives is established in Lemma~\ref{lem:invertible}.
The code, the input data with their sources and SHA-256 hashes, the complete output datasets, the test polynomials with their random seeds, and the logs of the computations described here are available in the \href{https://github.com/GeometrIALab/Link_SasakianHodgeNumbers}{GitHub} repository, from which every table and count in this paper can be regenerated.
The inputs are the weight systems and Kreuzer--Skarke Hodge numbers of the $7555$ threefolds as used in \cite{AHHHSS2024} (files \texttt{WP4s.txt} and \texttt{WP4\_Hodges.txt} of \cite{MLcCY7github}), the $7549$ Gr\"obner-basis computations of \cite{AHHHSS2024}, and the fourfold list \texttt{5dTransWH.all.gz} of \cite{KSfourfolds}, whose Hodge numbers are the chiral-ring data computed with Vafa's formula.

\paragraph{Checks of the sheaf-cohomological formula.}
To check the algebra in the proof of Theorem~\ref{thm:sheaf}, \eqref{eq:sheaf} was evaluated term by term, without using the telescoping.
The results were compared with:
\begin{enumerate}
\item the Hodge numbers of nine smooth hypersurfaces in projective space, including the cubic fourfold with $(h^{3,1},h^{2,2}_{\rm prim})=(1,20)$;
\item the coefficients of $P_{\bw}$ for $3000$ Hodge numbers of random non-degenerate weight systems with $N=3,\dots,6$, most of them not Calabi--Yau;
\item Gr\"obner-basis computations over $\Q$ of the Milnor algebras of invertible polynomials (Lemma~\ref{lem:invertible}) for $300$ of these systems;
\item all Sasakian Hodge numbers in Subsections~\ref{sec:5d} and~\ref{sec:9d}.
\end{enumerate}
For $N=3$, only the formula is tested: Theorem~\ref{thm:sheaf} assumes $n\ge2$ and therefore $N\ge4$.
No discrepancy was found.

\subsection{Certificates for explicit polynomials}
\label{app:certificates}
The checks below address two separate requirements: the critical point must be isolated, and the relevant graded piece must have the stated dimension.
The first is automatic for the representatives of Lemma~\ref{lem:invertible}; dense polynomials require a separate isolation test.

\paragraph{Isolation certificates.}
The remaining test polynomials contain all monomials of degree $d$: the $6+25$ polynomials with coefficients drawn uniformly from $\{1,\dots,999\}$ in Sections~\ref{sec:six} and~\ref{sec:9d}, and the $95+42$ polynomials of the follow-up computations, whose coefficients lie between $1$ and $4$ and are mostly equal to $1$.
For such a polynomial $f$ a Gr\"obner basis of the Jacobian ideal of the reduction $\bar f$ of $f$ modulo $p=32003$ was computed with Singular \cite{Singular}, and the quotient was checked to have Krull dimension $0$ and vector-space dimension $\mu(\bw)$.
In this setting, \emph{Krull dimension zero} means that the quotient is finite-dimensional; its dimension is the number of standard monomials recalled in Section~\ref{sec:weights}.
This implies that $f$ itself has an isolated critical point at the origin.
Indeed, in every degree $\ell$,
\begin{equation}
\label{eq:reduction-bound}
\dim_\Q(\MM_f)_\ell\le\dim_{\bF_p}(\MM_{\bar f})_\ell,
\end{equation}
because the rank of the integer matrix of the map
\[
\bigoplus_iS_{\ell-d+w_i}\longrightarrow S_\ell,\qquad
(g_i)\longmapsto\sum_ig_i\partial_if,
\]
can only drop modulo $p$.
Thus \eqref{eq:reduction-bound} shows that $\MM_f$ is finite-dimensional if $\MM_{\bar f}$ is.
Isolation was certified in this way for the six polynomials of Section~\ref{sec:six}, the $95$ K3 polynomials of Subsection~\ref{sec:5d}, the $25$ random fourfold polynomials of Subsection~\ref{sec:9d}(c), with $\mu(\bw)$ up to $9.2\times10^6$, and the $42$ fourfold polynomials of the follow-up computations.
No polynomial failed the test, and none was discarded or resampled.

\paragraph{Exact certificates for graded pieces.}
A \emph{strand} of a graded complex is obtained by taking one fixed weighted degree in every term.
This turns the Koszul complex from Theorem~\ref{thm:hilbert} into finite-dimensional linear algebra.
For an integer polynomial $f$ and a degree $\ell$, the relevant part of its degree-$\ell$ strand is
\begin{equation}
\label{eq:koszul-strand}
K_{2,\ell}\xrightarrow{\ D_2\ }K_{1,\ell}\xrightarrow{\ D_1\ }K_{0,\ell}=S_\ell,
\end{equation}
with
\[
K_{1,\ell}=\bigoplus_iS_{\ell-(d-w_i)},\qquad K_{2,\ell}=\bigoplus_{i<j}S_{\ell-(2d-w_i-w_j)},
\]
where $D_1(g_i)_i=\sum_ig_i\partial_if$ and $D_2$ maps the generator $e_{ij}$ to $\partial_if\,e_j-\partial_jf\,e_i$.
In the monomial bases these are integer matrices with $D_1D_2=0$, and $(\MM_f)_\ell=S_\ell/\mathrm{im}\,D_1$.
Let $r_k$ be the rank of $D_k$ modulo a prime $p$, and $c_1\coloneq \dim K_{1,\ell}$.
Since $r_k\le\rank_\Q D_k$ and $\rank_\Q D_1+\rank_\Q D_2\le c_1$, the condition
\begin{equation}
\label{eq:rank-certificate}
r_1+r_2=c_1
\end{equation}
forces $\rank_\Q D_k=r_k$ for both $k$, and then
\[\dim_\Q(\MM_f)_\ell=\dim S_\ell-r_1\]
exactly.
For example, if $c_1=10$ and the ranks modulo $p$ are $r_1=7$ and $r_2=3$, neither rational rank can be larger: increasing either would make their sum exceed $10$, contradicting $\mathrm{im}\,D_2\subseteq\ker D_1$.
This is why the pair of finite-field computations certifies an exact rational result; a single rank modulo $p$ would only give a lower bound.
Thus a \emph{rank certificate} is the checkable equality \eqref{eq:rank-certificate}, which forces the finite-field ranks to equal the rational ranks.
For an admissible $f$ the Koszul complex is exact by Theorem~\ref{thm:hilbert}, so \eqref{eq:rank-certificate} holds for all but finitely many primes.
It was tested with exact Gaussian elimination modulo $p=2^{31}-1$, without floating-point arithmetic, and it held at this prime in every case reported here, so that no other prime was needed.
For $\ell=d$ the space $K_{2,\ell}$ is zero in the Calabi--Yau case, and the condition is that $D_1$ has full rank.
The largest strands certified have $\dim S_\ell=1.8\times10^5$ for $\ell=d$ and $3.8\times10^5$ for $\ell=2d$; the certificate for $h^{2,2}_S$ was attempted when $\dim S_{2d}\le4\times10^5$.
The certificate does not by itself show that $f$ is admissible, which is established separately by Lemma~\ref{lem:invertible} or by the isolation certificates above.

Table~\ref{tab:checks} lists the checks with explicit polynomials.
The $42$ fourfold weight systems of the follow-up computations appear twice, with an invertible and with a dense polynomial; the other samples are disjoint.

\begin{table}[ht]
\centering\footnotesize
\begin{tabular}{llllll}\toprule
Section & $n$ & Polynomials & Isolation & Quantity certified & Method\\\midrule
\ref{sec:sheaf} & $1$--$4$ & $300$ invertible, $\mu<4000$ & Lemma~\ref{lem:invertible} & all middle $h^{p,q}_S$ & Gr\"obner basis over $\Q$\\
\ref{sec:six} & $3$ & $6$ invertible & Lemma~\ref{lem:invertible} & $\dim(\MM_f)_\ell$, $0\le\ell\le3d$ & rank certificate\\
\ref{sec:six} & $3$ & $6$ dense, random & Singular & $h^{2,1}_S$ & rank certificate\\
\ref{sec:5d} & $2$ & $95$ dense & Singular & $h^{1,1}_S$ & Gr\"obner basis\\
\ref{sec:9d} & $4$ & $42$ invertible & Lemma~\ref{lem:invertible} & Hilbert function & Gr\"obner basis over $\Q$\\
\ref{sec:9d} & $4$ & $900$ invertible & Lemma~\ref{lem:invertible} & $h^{3,1}_S$ ($900$), $h^{2,2}_S$ ($655$) & rank certificate\\
\ref{sec:9d} & $4$ & $25$ dense, random & Singular & $h^{3,1}_S$, $h^{2,2}_S$ & rank certificate\\
\ref{sec:9d} & $4$ & $42$ dense & Singular & $h^{3,1}_S$ & rank certificate\\\bottomrule
\end{tabular}
\caption{Checks with explicit polynomials.
All certified values agree with the closed formulae.}
\label{tab:checks}
\end{table}

\paragraph{Computing environment.}
The symbolic and numerical verifications reported in the manuscript used Python (\texttt{sympy}, \texttt{numpy}), Singular \cite{Singular} and Macaulay2 \cite{M2}.
Unless stated otherwise, the computations were run with Python~3.10.13 and \texttt{numpy}~1.26.4, and the reruns recorded in the \href{https://github.com/GeometrIALab/Link_SasakianHodgeNumbers}{GitHub} repository with \texttt{numpy}~2.2.6, on an Apple M1 laptop with eight cores and $16$\,GB of memory under macOS~26, using exact integer arithmetic throughout.
The $322$\,s for the fourfold dataset is the wall-clock time of eight worker processes that compute all invariants of Subsection~\ref{sec:9d}, the twisted surfaces and the comparison (a), including writing the output file and excluding reading the input file.
The certificates of Section~\ref{sec:six} took about $80$\,s on one core, and those of Subsection~\ref{sec:9d}(c) took $69$\,s for the $900$ systems with invertible polynomials and $15$\,s for the $25$ systems without, with eight processes.
The isolation certificates were computed with Singular~4.4.1 \cite{Singular} on the same machine while it was shared with other computations; individual polynomials took between a fraction of a second and $86$ minutes.

\subsection{Complete K3 dataset}
\label{app:k3data}
Table~\ref{tab:k3} contains the $95$ K3 links discussed in Subsection~\ref{sec:5d}.
The Hodge and Betti numbers use \eqref{eq:k3}, and the singularity column
records the Du Val types entering Proposition~\ref{prop:k3}.
The table also recovers the topological classification recalled in
Remark~\ref{rem:k3top}.

The equivalent full dataset for the nine-dimensional links is available in this work's respective \href{https://github.com/GeometrIALab/Link_SasakianHodgeNumbers}{GitHub} repository.

{\small
\begin{longtable}{lrrrlr}
\caption{The $95$ five-dimensional K3 links, ordered by degree and then lexicographically by weights: weights $\bw$, degree $d$, Sasakian Hodge number $h^{1,1}_S$ from \eqref{eq:k3}, second Betti number $b_2(K_f)$, singularities of the K3 surface $X_f$, and Milnor number $\mu(\bw)$.
In every case $h^{2,0}_S=1$, and the singularities have $19-h^{1,1}_S$ exceptional curves in total (Proposition~\ref{prop:k3}).}
\label{tab:k3}
\\
\toprule $\bw$ & $d$ & $h^{1,1}_S$ & $b_2$ & $\mathrm{Sing}\,X_f$ & $\mu$\\\midrule\endfirsthead
\multicolumn{6}{l}{\footnotesize Table~\ref{tab:k3} continued.}\\
\toprule $\bw$ & $d$ & $h^{1,1}_S$ & $b_2$ & $\mathrm{Sing}\,X_f$ & $\mu$\\\midrule\endhead
\bottomrule\endfoot
$(1,1,1,1)$ & 4 & 19 & 21 & -- & 81 \\
$(1,1,1,2)$ & 5 & 18 & 20 & $A_{1}$ & 96 \\
$(1,1,1,3)$ & 6 & 19 & 21 & -- & 125 \\
$(1,1,2,2)$ & 6 & 16 & 18 & $3A_{1}$ & 100 \\
$(1,1,2,3)$ & 7 & 16 & 18 & $A_{2}+A_{1}$ & 120 \\
$(1,1,2,4)$ & 8 & 17 & 19 & $2A_{1}$ & 147 \\
$(1,2,2,3)$ & 8 & 13 & 15 & $A_{2}+4A_{1}$ & 105 \\
$(1,1,3,4)$ & 9 & 16 & 18 & $A_{3}$ & 160 \\
$(1,2,3,3)$ & 9 & 12 & 14 & $3A_{2}+A_{1}$ & 112 \\
$(1,1,3,5)$ & 10 & 17 & 19 & $A_{2}$ & 189 \\
$(1,2,2,5)$ & 10 & 14 & 16 & $5A_{1}$ & 144 \\
$(1,2,3,4)$ & 10 & 12 & 14 & $A_{3}+A_{2}+2A_{1}$ & 126 \\
$(1,2,3,5)$ & 11 & 12 & 14 & $A_{4}+A_{2}+A_{1}$ & 144 \\
$(1,1,4,6)$ & 12 & 18 & 20 & $A_{1}$ & 242 \\
$(1,2,3,6)$ & 12 & 13 & 15 & $2A_{2}+2A_{1}$ & 165 \\
$(1,2,4,5)$ & 12 & 12 & 14 & $A_{4}+3A_{1}$ & 154 \\
$(1,3,4,4)$ & 12 & 10 & 12 & $3A_{3}$ & 132 \\
$(2,2,3,5)$ & 12 & 9 & 11 & $A_{4}+6A_{1}$ & 105 \\
$(2,3,3,4)$ & 12 & 8 & 10 & $4A_{2}+3A_{1}$ & 90 \\
$(1,3,4,5)$ & 13 & 10 & 12 & $A_{4}+A_{3}+A_{2}$ & 144 \\
$(1,2,4,7)$ & 14 & 13 & 15 & $A_{3}+3A_{1}$ & 195 \\
$(2,2,3,7)$ & 14 & 10 & 12 & $A_{2}+7A_{1}$ & 132 \\
$(2,3,4,5)$ & 14 & 7 & 9 & $A_{4}+A_{3}+A_{2}+3A_{1}$ & 99 \\
$(1,2,5,7)$ & 15 & 12 & 14 & $A_{6}+A_{1}$ & 208 \\
$(1,3,4,7)$ & 15 & 10 & 12 & $A_{6}+A_{3}$ & 176 \\
$(1,3,5,6)$ & 15 & 10 & 12 & $A_{5}+2A_{2}$ & 168 \\
$(2,3,5,5)$ & 15 & 6 & 8 & $3A_{4}+A_{1}$ & 104 \\
$(3,3,4,5)$ & 15 & 6 & 8 & $A_{3}+5A_{2}$ & 88 \\
$(1,2,5,8)$ & 16 & 13 & 15 & $A_{4}+2A_{1}$ & 231 \\
$(1,3,4,8)$ & 16 & 11 & 13 & $2A_{3}+A_{2}$ & 195 \\
$(1,4,5,6)$ & 16 & 9 & 11 & $A_{5}+A_{4}+A_{1}$ & 165 \\
$(2,3,4,7)$ & 16 & 7 & 9 & $A_{6}+A_{2}+4A_{1}$ & 117 \\
$(2,3,5,7)$ & 17 & 6 & 8 & $A_{6}+A_{4}+A_{2}+A_{1}$ & 120 \\
$(1,2,6,9)$ & 18 & 14 & 16 & $A_{2}+3A_{1}$ & 272 \\
$(1,3,5,9)$ & 18 & 11 & 13 & $A_{4}+2A_{2}$ & 221 \\
$(1,4,6,7)$ & 18 & 9 & 11 & $A_{6}+A_{3}+A_{1}$ & 187 \\
$(2,3,4,9)$ & 18 & 8 & 10 & $A_{3}+2A_{2}+4A_{1}$ & 140 \\
$(2,3,5,8)$ & 18 & 6 & 8 & $A_{7}+A_{4}+2A_{1}$ & 130 \\
$(3,4,5,6)$ & 18 & 5 & 7 & $A_{4}+A_{3}+3A_{2}+A_{1}$ & 91 \\
$(3,4,5,7)$ & 19 & 4 & 6 & $A_{6}+A_{4}+A_{3}+A_{2}$ & 96 \\
$(1,4,5,10)$ & 20 & 10 & 12 & $2A_{4}+A_{1}$ & 228 \\
$(2,3,5,10)$ & 20 & 7 & 9 & $2A_{4}+A_{2}+2A_{1}$ & 153 \\
$(2,4,5,9)$ & 20 & 6 & 8 & $A_{8}+5A_{1}$ & 132 \\
$(2,5,6,7)$ & 20 & 5 & 7 & $A_{6}+A_{5}+3A_{1}$ & 117 \\
$(3,4,5,8)$ & 20 & 4 & 6 & $A_{7}+2A_{3}+A_{2}$ & 102 \\
$(1,3,7,10)$ & 21 & 10 & 12 & $A_{9}$ & 264 \\
$(1,5,7,8)$ & 21 & 8 & 10 & $A_{7}+A_{4}$ & 208 \\
$(2,3,7,9)$ & 21 & 6 & 8 & $A_{8}+2A_{2}+A_{1}$ & 152 \\
$(3,5,6,7)$ & 21 & 4 & 6 & $A_{5}+A_{4}+3A_{2}$ & 96 \\
$(1,3,7,11)$ & 22 & 11 & 13 & $A_{6}+A_{2}$ & 285 \\
$(1,4,6,11)$ & 22 & 10 & 12 & $A_{5}+A_{3}+A_{1}$ & 252 \\
$(2,4,5,11)$ & 22 & 7 & 9 & $A_{4}+A_{3}+5A_{1}$ & 153 \\
$(1,3,8,12)$ & 24 & 12 & 14 & $A_{3}+2A_{2}$ & 322 \\
$(1,6,8,9)$ & 24 & 8 & 10 & $A_{8}+A_{2}+A_{1}$ & 230 \\
$(2,3,7,12)$ & 24 & 7 & 9 & $A_{6}+2A_{2}+2A_{1}$ & 187 \\
$(2,3,8,11)$ & 24 & 6 & 8 & $A_{10}+3A_{1}$ & 182 \\
$(3,4,5,12)$ & 24 & 5 & 7 & $A_{4}+2A_{3}+2A_{2}$ & 133 \\
$(3,4,7,10)$ & 24 & 3 & 5 & $A_{9}+A_{6}+A_{1}$ & 119 \\
$(3,6,7,8)$ & 24 & 4 & 6 & $A_{6}+4A_{2}+A_{1}$ & 102 \\
$(4,5,6,9)$ & 24 & 3 & 5 & $A_{8}+A_{4}+A_{2}+2A_{1}$ & 95 \\
$(4,5,7,9)$ & 25 & 2 & 4 & $A_{8}+A_{6}+A_{3}$ & 96 \\
$(1,5,7,13)$ & 26 & 9 & 11 & $A_{6}+A_{4}$ & 285 \\
$(2,3,8,13)$ & 26 & 7 & 9 & $A_{7}+A_{2}+3A_{1}$ & 207 \\
$(2,5,6,13)$ & 26 & 6 & 8 & $A_{5}+A_{4}+4A_{1}$ & 168 \\
$(2,5,9,11)$ & 27 & 4 & 6 & $A_{10}+A_{4}+A_{1}$ & 160 \\
$(5,6,7,9)$ & 27 & 2 & 4 & $A_{6}+A_{5}+A_{4}+A_{2}$ & 88 \\
$(1,4,9,14)$ & 28 & 10 & 12 & $A_{8}+A_{1}$ & 342 \\
$(3,4,7,14)$ & 28 & 4 & 6 & $2A_{6}+A_{2}+A_{1}$ & 150 \\
$(4,6,7,11)$ & 28 & 2 & 4 & $A_{10}+A_{5}+2A_{1}$ & 102 \\
$(1,4,10,15)$ & 30 & 11 & 13 & $A_{4}+A_{3}+A_{1}$ & 377 \\
$(1,6,8,15)$ & 30 & 9 & 11 & $A_{7}+A_{2}+A_{1}$ & 319 \\
$(2,3,10,15)$ & 30 & 8 & 10 & $A_{4}+2A_{2}+3A_{1}$ & 252 \\
$(2,6,7,15)$ & 30 & 6 & 8 & $A_{6}+A_{2}+5A_{1}$ & 184 \\
$(3,4,10,13)$ & 30 & 3 & 5 & $A_{12}+A_{3}+A_{1}$ & 153 \\
$(4,5,6,15)$ & 30 & 4 & 6 & $2A_{4}+A_{3}+A_{2}+2A_{1}$ & 130 \\
$(5,6,8,11)$ & 30 & 1 & 3 & $A_{10}+A_{7}+A_{1}$ & 95 \\
$(2,5,9,16)$ & 32 & 5 & 7 & $A_{8}+A_{4}+2A_{1}$ & 207 \\
$(4,5,7,16)$ & 32 & 3 & 5 & $A_{6}+A_{4}+2A_{3}$ & 135 \\
$(3,5,11,14)$ & 33 & 2 & 4 & $A_{13}+A_{4}$ & 152 \\
$(3,4,10,17)$ & 34 & 4 & 6 & $A_{9}+A_{3}+A_{2}+A_{1}$ & 186 \\
$(4,6,7,17)$ & 34 & 3 & 5 & $A_{6}+A_{5}+A_{3}+2A_{1}$ & 135 \\
$(1,5,12,18)$ & 36 & 10 & 12 & $A_{5}+A_{4}$ & 434 \\
$(3,4,11,18)$ & 36 & 4 & 6 & $A_{10}+2A_{2}+A_{1}$ & 200 \\
$(7,8,9,12)$ & 36 & 1 & 3 & $A_{7}+A_{6}+A_{3}+A_{2}$ & 87 \\
$(3,5,11,19)$ & 38 & 3 & 5 & $A_{10}+A_{4}+A_{2}$ & 189 \\
$(5,6,8,19)$ & 38 & 2 & 4 & $A_{7}+A_{5}+A_{4}+A_{1}$ & 132 \\
$(5,7,8,20)$ & 40 & 2 & 4 & $A_{6}+2A_{4}+A_{3}$ & 132 \\
$(1,6,14,21)$ & 42 & 10 & 12 & $A_{6}+A_{2}+A_{1}$ & 492 \\
$(2,5,14,21)$ & 42 & 6 & 8 & $A_{6}+A_{4}+3A_{1}$ & 296 \\
$(3,4,14,21)$ & 42 & 5 & 7 & $A_{6}+A_{3}+2A_{2}+A_{1}$ & 247 \\
$(4,5,13,22)$ & 44 & 2 & 4 & $A_{12}+A_{4}+A_{1}$ & 186 \\
$(3,5,16,24)$ & 48 & 4 & 6 & $A_{7}+A_{4}+2A_{2}$ & 258 \\
$(7,8,10,25)$ & 50 & 1 & 3 & $A_{7}+A_{6}+A_{4}+A_{1}$ & 129 \\
$(4,5,18,27)$ & 54 & 3 & 5 & $A_{8}+A_{4}+A_{3}+A_{1}$ & 245 \\
$(5,6,22,33)$ & 66 & 2 & 4 & $A_{10}+A_{4}+A_{2}+A_{1}$ & 244 \\
\end{longtable}
}

\phantomsection
\addcontentsline{toc}{section}{References}
\bibliographystyle{utphys}
\bibliography{references}{}

\end{document}